\documentclass[referee,pdflatex,sn-basic]{sn-jnl}% referee option is meant for double line spacing

\usepackage{endnotes}
\usepackage{natbib}
\usepackage{url}          
\usepackage{enumitem} 
\usepackage{bm}        
\usepackage{setspace} 
\usepackage[ruled]{algorithm2e}
\usepackage{pdflscape}
\usepackage{multirow}
\usepackage[dvipsnames]{xcolor}
\usepackage{tikz}
\usetikzlibrary{shapes.geometric, arrows, automata, positioning}
\usepackage{graphicx}%
\usepackage{multirow}%
\usepackage{amsmath,amssymb,amsfonts}%
\usepackage{amsthm}%
\usepackage{mathrsfs}%
\usepackage[title]{appendix}%
\usepackage{textcomp}%
\usepackage{manyfoot}%
\usepackage{booktabs}%
\usepackage{algpseudocode}%
\usepackage{listings}%
\usepackage{footnote}%
\usepackage{pgfplots}%
\theoremstyle{plain}%
\newtheorem{theorem}{Theorem}%  meant for continuous numbers
\newtheorem{proposition}[theorem]{Proposition}% 
\newtheorem{lemma}{Lemma}

\theoremstyle{remark}%
\newtheorem{remark}{Remark}%

\theoremstyle{definition}%

\DeclareRobustCommand{\cpluspluslogo}{\hbox{C\hspace{-0.5ex}
                       \protect\raisebox{0.5ex}
                       {\protect\scalebox{0.67}{++}}}}

\def\ROCPP{RO\cpluspluslogo}

\newcommand{\field}[1]{\ensuremath{\mathbb{#1}}}
\newcommand{\sets}[1]{\ensuremath{\mathcal{#1}}}

\newcommand{\reals}{\ensuremath{\field{R}}} % real numbers
\newcommand{\1}{\ensuremath{{\rm \mathbf e}}} % vector of all 1's
\definecolor{DarkGreen}{rgb}{0, 0.5, 0.1}
\DeclareMathOperator{\st}{s.t.}
\DeclareMathOperator{\interior}{int}
\DeclareMathOperator*{\argmin}{\mathrm{argmin}}

\newcommand{\newqj}[1]{{\color{black}{#1}}}

\usepackage{xr}
\makeatletter
\newcommand*{\addFileDependency}[1]{
	\typeout{(#1)}
	\@addtofilelist{#1}
	\IfFileExists{#1}{}{\typeout{No file #1.}}
}
\makeatother

\newcommand*{\myexternaldocument}[1]{
	\externaldocument{#1}
	\addFileDependency{#1.tex}
	\addFileDependency{#1.aux}
}
\myexternaldocument{response-to-reviewers}

\begin{document}
%%%%%%%%%%%%%%%%

\title[Article Title]{Distributionally Robust Optimization with Decision-Dependent Information Discovery}


\author*[1,2]{\fnm{Qing} \sur{Jin}}\email{qingjin@usc.edu}

\author[3]{\fnm{Angelos} \sur{Georghiou}}\email{georghiou.angelos@ucy.ac.cy}

\author[1,2,4]{\fnm{Phebe} \sur{Vayanos}}\email{phebe.vayanos@usc.edu}

\author[5]{\fnm{Grani A.} \sur{Hanasusanto}}\email{gah@illinois.edu}

\affil[1]{\orgdiv{Department of Industrial and Systems Engineering}, \orgname{University of Southern California}}

\affil[2]{\orgdiv{Center for Artificial Intelligence in Society}, \orgname{University of Southern California}}

\affil[3]{\orgdiv{Department of Business and Public Administration}, \orgname{University of Cyprus}}

\affil[4]{\orgdiv{Department of Computer Science}, \orgname{University of Southern California}}

\affil[5]{\orgdiv{Department of Industrial \& Enterprise Systems Engineering}, \orgname{University of Illinois Urbana-Champaign}}

\abstract{%
%\begin{spacing}{1.0}
We study two-stage distributionally robust optimization (DRO) problems with decision-dependent information discovery (DDID) wherein (a portion of) the uncertain parameters are revealed only if an (often costly) investment is made at the first stage. This class of problems finds many important applications in selection problems (e.g., in hiring, project portfolio optimization, or optimal sensor location). Despite the \newqj{wide applicability of the problem}, it has not been previously studied. We propose a framework for modeling and approximately solving DRO problems with DDID. We formulate the problem as a min-max-min-max problem and adopt the popular~\unboldmath{$K$}-adaptability approximation scheme, which chooses~$K$ candidate recourse actions here-and-now and implements the best of those actions after the uncertain parameters that were chosen to be observed are revealed. We then present a decomposition algorithm that solves the~$K$-adaptable formulation exactly. In particular, we devise a cutting plane algorithm that iteratively solves a relaxed version of the problem, evaluates the true objective value of the corresponding solution, generates valid cuts, and imposes them on the relaxed problem. For the evaluation problem, we develop a branch-and-cut algorithm that provably converges to an optimal solution. We showcase the effectiveness of our framework on the R\&D project portfolio optimization problem and the best box problem. 
%\end{spacing}
% Enter your abstract
}%


\maketitle

%\newpage
%%%%%%%%%%%%%%%%%%%%%%%%%%%%%%%%%%%%%%%%%%%%%%%%%%%%%%%%%%%%%%%%%%%%%%

%%%%%%%%%%%%%%%%%%%%%%%%%%%%%%%%%%%%%%%%%%%%%%%%%%%%%%%%%%%%%%%%%%%%%%%%%%%%%%%%%%%%%%%%%%%%%%%%
%%%%%%%%%%%%%%%%%%%%%%%%%%%%%%%%%%%%%%%%%%%%%%%%%%%%%%%%%%%%%%%%%%%%%%%%%%%%%%%%%%%%%%%%%%%%%%%%
%%%%%%%%%%%%%%%%%%%%%%%%%%%%%%%%%%%% Introductions %%%%%%%%%%%%%%%%%%%%%%%%%%%%%%%%%%%%%%%
%%%%%%%%%%%%%%%%%%%%%%%%%%%%%%%%%%%%%%%%%%%%%%%%%%%%%%%%%%%%%%%%%%%%%%%%%%%%%%%%%%%%%%%%%%%%%%%%
%%%%%%%%%%%%%%%%%%%%%%%%%%%%%%%%%%%%%%%%%%%%%%%%%%%%%%%%%%%%%%%%%%%%%%%%%%%%%%%%%%%%%%%%%%%%%%%%
\section{Introduction}
\label{sec:introductions}

\subsection{Background \& Motivation}

Distributionally robust optimization, which optimizes for a safe action that performs best under the most adverse distribution in an \emph{ambiguity set} of distributions consistent with known information, is a popular modeling and solution paradigm for decision-making under uncertainty. It can be used to tackle both single-stage (see, e.g., \citet{Delage2010DistributionallyProblems}, \citet{Ben-Tal2013RobustProbabilities}, \citet{Wiesemann2014}, \citet{Jiang2016Data-drivenProgram}, \citet{gao2023distributionally}, \citet{MohajerinEsfahani2018Data-drivenReformulations}) and multi-stage (see e.g., \citet{Goh2010DistributionallyApproximations}, \citet{Hanasusanto2018ConicBalls}, \citet{Bertsimas2019AdaptiveOptimization}, \citet{Chen2020RobustRSOME}, \citet{Yu2020MultistageSets}) problems. In multi-stage problems, uncertain parameters are revealed sequentially over time, and decisions are allowed to \emph{adapt} to the history of observations. Mathematically, decisions must be modeled as functions of the uncertain parameters that have been revealed by the time the decision is made.

To the best of our knowledge, all models and solution approaches assume that the sequence in which uncertain parameters are revealed is \emph{exogenous}, i.e., \emph{independent} of the decision-maker's actions. However, this assumption fails to hold in many real-world applications where uncertain parameters (\emph{information}) only become observable following an often costly investment and \emph{measurement decisions} control the \emph{time of information discovery}. We now describe some important examples of such problems \newqj{that involve} \emph{decision-dependent information discovery}.

\textit{R\&D Project Portfolio Optimization.} Companies often maintain long pipelines of projects that are candidates to be undertaken~\citep{Solak2010OptimizationUncertainty}. The return of each project is uncertain and will only be revealed if the project is completed, which often necessitates substantial resources to be expanded. Thus, the costly decisions to undertake and complete a project are measurement decisions \newqj{that} control the time of information discovery in this problem.

\textit{Optimal Sensor Placement for Emergency Warning.} Emergency warning systems (that issue warnings about, e.g., landslides~\citep{Chu2021SitkaNet:Study}) rely on information collected from sensors (about, e.g., soil moisture) to make predictions. Sensors are often costly or difficult to install (due to e.g., difficulty in accessing remote locations). The information available to help inform the warning system depends on the location of the sensors. Thus, the costly sensor placement decisions are measurement decisions that control the time of information discovery in this problem.

\textit{Selection Problems.} Selection problems\newqj{,} such as those arising in company hiring, loan approval, or bail decisions\newqj{,} involve selecting, over time, a subset of available candidates to move to the next stage until a final selection decision is made in the last stage. The personal characteristics (e.g., qualifications, gender, race) of individuals are uncertain and only become observed progressively over time as they get selected into subsequent stages. Accordingly, the outcomes  (e.g., job performance, loan repayment, recidivism) are also uncertain and only become observed for individuals who make it into the final stage. Thus, the candidate selection decisions are measurement decisions that control the time of information discovery in this problem. 

In most real-world situations, the precise distribution of the uncertain parameters in such problems is unknown (e.g., the joint distribution of project returns in R\&D, \newqj{the} joint distribution of soil characteristics at different remote locations, \newqj{the} joint distribution of job candidates' characteristics, and corresponding job performance). Yet, basic characteristics of this joint distribution, such as first-order moments, may be available (e.g., through past project performance, current weather conditions, and past worker performance data). In our work, we propose to leverage such information in \emph{distributionally robust optimization problems with decision-dependent information discovery} where such moment information is added explicitly as constraints in the ambiguity set.

\subsection{Literature Review}
\label{subsec:lit-review}

{Our work most closely relates to the literature on optimization under decision-dependent information discovery, to distributionally robust optimization with decision-dependent uncertainty sets, and to the works on the $K$-adaptability approximation to (distributionally) robust optimization. We review all of these in turn.}

{Decision-making problems under decision-dependent information discovery have mostly been studied in the stochastic programming literature. The majority of works assume the distribution of the uncertain parameters is discrete or require the distribution to be discretized before they can be applied, see \citet{colvin2010modeling,goel2004stochastic,goel2006class,gupta2011solution,gupta2014new,APAP2017233}. \citet{Vayanos2011DecisionProgramming} investigate a problem with continuous uncertain parameters and propose a conservative solution approach based on piecewise constant and piecewise linear decision rule approximations. In the robust optimization literature, \citet{Zhang2020AUncertainty} derive a \newqj{decision rule} approximation for multi\newqj{-}stage robust optimization with endogenous uncertainty. \citet{paradiso2022exact} develop an exact solution scheme based on a nested decomposition algorithm for two-stage robust optimization problems with DDID and objective uncertainty. \citet{omer:hal-04097679} provide a mixed-integer linear programming formulation for robust combinatorial problems with DDID under polyhedral uncertainty through column generation and constraint generation algorithms.

{Distributionally robust optimization with a decision-dependent ambiguity set has been studied in both static and adaptive settings. Such a framework can model the case where the distribution of the underlying random variables depends on decisions. One class of decision-dependent ambiguity set\newqj{s} allows moments of the distributions to depend on the decisions, see \citet{zhang2016quantitative,Luo2020DistributionallySets,ryu2019nurse,Basciftci2021DistributionallyDemand}, and \citet{Yu2020MultistageSets}. Another class allows decisions to influence the nominal or marginal distributions; see \citet{Luo2020DistributionallySets,Noyan2022DistributionallyLogistics}, and \citet{Doan2021DistributionallyPlanning}. None of the models considers the case where decisions affect the time of information discovery.}

{The $K$-adaptability approximation scheme, originally introduced by \citet{Bertsimas2010FiniteOptimization} for two-stage robust optimization, involves selecting $K$ candidate policies here-and-now and implementing the best of these policies after the uncertain parameters are revealed. This method naturally accommodates discrete adaptive variables. \citet{Bertsimas2010FiniteOptimization} study the complexity of the problem and provide an exact finite\newqj{-}dimensional bilinear formulation for the case where~$K=2$. \citet{Subramanyam2020} propose a branch-and-bound algorithm to speed-up computation. For the class of problems involving only binary adaptive variables, \citet{Hanasusanto2015} characterize the problem’s complexity in terms of the number of second-stage policies~$K$ needed to obtain an optimal solution to the original, fully adaptive problem and derive practical explicit mixed-integer linear optimization (MILO) reformulations. \citet{buchheim2017min,buchheim2018complexity} study the complexity of two-stage robust combinatorial optimization problems under objective uncertainty and propose efficient algorithms for the cases of polyhedral and discrete uncertainty sets, respectively. \citet{Chassein2019} propose a faster algorithm for problems with \newqj{a} budgeted uncertainty set. \newqj{\citet{kurtz2024many} derives bounds on the number of policies~$K$ which guarantee optimality for non-linear two-stage robust optimization problems with integer decisions, and extends the results to the case with DDID.} \citet{Hanasusanto2016} adapt the $K$-adaptability approximation for binary decision variables to the case of two-stage distributionally robust optimization problems, derive explicit mixed-integer linear programming reformulations, and provide efficient methods for bounding the selection probabilities of each of the $K$ second-stage policies. \citet{han2023finite} provide a single-stage robust optimization reformulation for computing optimal $K$-adaptable policies with optimal partitions for a two-stage distributionally robust optimization problem with constraint uncertainty, show the reformulated problem is~$\mathcal{NP}$-complete, and propose a partitioning framework that results in a mixed-integer optimization problem of moderate size. Our work most closely relates to the paper of \citet{Vayanos2020}, wherein the authors study a robust optimization problem with DDID and propose a conservative approximation based on~$K$-adaptability. We study the more general class of \emph{distributionally} robust problems with DDID.}

% \notepv{what is the role of this paragraph?} Decomposition methods have been well studied in stochastic programming problems (see~\cite{Laporte1993TheRecourse},~\cite{Angulo2016ImprovingMethod},~\cite{Zou2019StochasticProgramming}) and robust optimization problems (see~\cite{Zhao2012AnProblems},~\cite{Bertsimas2016MultistagePartitions},~\cite{Subramanyam2020}). Recently, many researchers have leveraged its flexibility and extended it to the DRO setting.~\cite{Chen2021DecompositionProblems} propose an algorithmic framework based on Lagrangian decomposition for solving the two-stage DRO problem with finite support in the ambiguity set.~\cite{Luo2019DecompositionModels} study DRO problem with Wasserstein ambiguity set and reformulate it into a decomposable semi-infinite program. They then solve the reformulated problem with a cutting-surface algorithm. To solve a two-stage mixed-binary DRO problem,~\cite{Bansal2018DecompositionPrograms} present a decomposition algorithm utilizing a distribution separation procedure and forcing parametric cuts.~\cite{Yu2020MultistageSets} extend the general moment-based ambiguity set to the multistage decision-dependent setting and derive MISDP reformulations of stage-wise subproblems. They deploy the Stochastic Dual Dynamic integer Programming method for solving the problem under the three ambiguity sets with risk-neutral or risk-averse objective functions.  
    
\subsection{{Proposed Approach \&} Contributions}
\label{subsec:contributions} 

We now summarize our proposed approach and main contributions in this paper.

\begin{enumerate}

\item We study two-stage distributionally robust optimization problems with decision-dependent information discovery. We formulate this class of problems as adaptive distributionally robust optimization problems that optimize over \newqj{decision rules} on which we impose decision-dependent non-anticipativity constraints. We establish that this model admits an equivalent two-and-a-half\newqj{-}stage (min-max-min-max) robust optimization reformulation and adapt the popular~$K$-adaptability approach to obtain a conservative approximation. The parameter~$K$, which denotes the number of candidate policies, enables us to conveniently trade off computational complexity with solution quality.

\item We show that applying the~$K$-adaptability approximation to the reformulated robust optimization problems leads to a bilinear optimization problem {involving products of \emph{real-valued} decision variables} (see Theorem~\ref{thm:cstrunc-bilinear}). For moderately sized problems, this can be solved with off-the-shelf solvers. For larger instances, we propose a decomposition algorithm that solves the problem exactly. The scheme iteratively solves \newqj{the} main problem with only the measurement variables and adds optimality, feasibility, and Benders-like cuts generated by an evaluation problem. We develop a provably convergent branch-and-cut algorithm to evaluate the objective with fixed measurement variables.

% We present a decomposition algorithm that solves the~$K$-adaptable formulation exactly. The algorithm iteratively finds the current optimal information discovery solution, evaluates the objective value with the optimal solution, and adds cuts. For the evaluation problem, where the measurement variables are fixed, we develop a branch-and-cut algorithm that is provably convergent. We further enhance the decomposition algorithm by enforcing information cut derived from the information discovery nature. 

\item We conduct numerical experiments on stylized instances of the R\&D project portfolio optimization and best box problems that showcase the effectiveness of our approach. We show that the solutions to the $K$-adaptability approximation can be obtained by either solving the bilinear optimization problem directly or by employing our decomposition algorithm. In the R\&D project portfolio optimization problem, the decomposition algorithm outperforms the state-of-the-art solver, finding better solutions in up to 100\% more cases.  We assess the value of incorporating distributional information by comparing our method to that returned by a purely robust solution. Our results show that the distributionally robust optimization solution improves the robust solution by up to 51\% on average in terms of the relative difference in the objective value under the worst-case distribution and by 49\% in out-of-sample performance.
\end{enumerate}

\subsection{Organization of the Paper and Notations}

The remainder of the paper is organized as follows. Section~\ref{sec:DRO-DDDID} introduces the two-stage distributionally robust optimization problem with the decision-dependent information discovery and derives an equivalent min-max-min-max robust optimization counterpart. Section~\ref{sec:k-adapt} presents the~$K$-adaptability approximation scheme, and Section~\ref{sec:l-shaped} describes the decomposition algorithm. Section~\ref{sec:numerical} conducts numerical experiments to test the performance of our approach. All proofs are relegated to the appendix.

\textit{Notation.} Throughout the paper, vectors (matrices) are denoted by boldface lowercase (uppercase) letters. Uncertainty is modeled by the probability space $\left(\mathbb{R}^k,\mathcal{B}\left(\mathbb{R}^k\right), \mathbb{P}\right)$, which consists of the sample space $\mathbb{R}^k$, the Borel $\sigma$-algebra $\mathcal{B}\left(\mathbb{R}^k\right)$ and the probability measure~$\mathbb{P}$, whose support we denote by $\Xi$. We let~$\sets L_{N_\xi}^{N_{y}}$ be the space of all measurable functions from~$\reals^{N_{\xi}}$ to~$\reals^{N_{y}}$ that are bounded on compact sets. Given two vectors~$\bm x$ and $\bm y$ of the same dimension, $\bm x \circ \bm y$ denotes their Hadamard product. Given two sets~$\sets X$ and~$\sets Y,\; \sets X \xrightarrow{L} \sets Y$ denotes a linear mapping from~$\sets X$ to~$\sets Y$. We use~$\bm A_i$ and $[\bm A]_i$ to denote the~$i$th column and the~$i$th row of matrix~$\bm A$, respectively. With a slight abuse of notations, we use the maximum and minimum operators even when the optimum may not be achieved. In such instances, these operators should be interpreted as suprema and infima, respectively.

%%%%%%%%%%%%%%%%%%%%%%%%%%%%%%%%%%%%%%%%%%%%%%%%%%%%%%%%%%%%%%%%%%%%%%%%%%%%%%%%%%%%%%%%%%%%%%%%
%%%%%%%%%%%%%%%%%%%%%%%%%%%%%%%%%%%%%%%%%%%%%%%%%%%%%%%%%%%%%%%%%%%%%%%%%%%%%%%%%%%%%%%%%%%%%%%%
%%%%%%%%%%%%%%%%%%%%%%%%%%%%%%%%%%%% DRO with Decision-Dependent Information Discovery %%%%%%%%%%%%%%%%%%%%%%%%%%%%%%%%%%%%%%%
%%%%%%%%%%%%%%%%%%%%%%%%%%%%%%%%%%%%%%%%%%%%%%%%%%%%%%%%%%%%%%%%%%%%%%%%%%%%%%%%%%%%%%%%%%%%%%%%
%%%%%%%%%%%%%%%%%%%%%%%%%%%%%%%%%%%%%%%%%%%%%%%%%%%%%%%%%%%%%%%%%%%%%%%%%%%%%%%%%%%%%%%%%%%%%%%%
\section{DRO with Decision-Dependent Information Discovery}
\label{sec:DRO-DDDID}
In this section, we describe the two-stage distributionally robust optimization problems with decision-dependent information discovery. We propose a decision rule-based formulation for this problem and demonstrate that it is equivalent to a two-and-a-half\newqj{-}stage min-max-min-max robust problem.

\subsection{Decision Rule Formulation}
\label{dr-formulation}

The DRO-DDID problem minimizes the worst-case expected objective value over all distributions~$\mathbb P$ of the uncertain parameters $\bm\xi\in\mathbb{R}^{N_{\xi}}$, supported on $\Xi$, and belonging to the ambiguity set $\mathcal{P}$. The ambiguity set contains plausible distributions that share certain properties that are known about the true distribution. The first-stage decisions~$\bm x \in \sets X \subseteq\reals^{N_{x}}$, and $\bm w \in \sets W \subseteq \{0,1\}^{N_\xi}$ are selected here-and-now, before any of the uncertain parameters~$\bm \xi \in \Xi \subseteq \mathbb{R}^{N_{\xi}}$ are realized. The wait-and-see decisions~$\bm y(\bm \xi) \in \sets Y \subseteq \reals^{N_y}$ are dependent on the observed portion of~$\bm \xi$. The sets $\mathcal X$ and $\mathcal Y$ may involve integrality constraints, implying that both the first and second stage decisions may involve a mix of real- and discrete-valued decisions. The binary vector~$\bm w$ collects the measurement decisions determining which components of the uncertain vector~$\bm \xi$ to observe. Specifically,~$\xi_i$ is revealed between the first and second stages if and only if~$w_i = 1$.

The two-stage distributionally robust problem with decision-dependent information discovery can be written mathematically as
\begin{equation}
    \begin{array}{cl}
         \min & \quad \displaystyle \sup_{\mathbb{P} \in \mathcal{P}} \;\;\;\mathbb{E}_{\mathbb{P}}\Big(  {\bm \xi}^\top {\bm C} \; {\bm x} + {\bm \xi}^\top {\bm D} \; {\bm w}  + {\bm \xi}^\top {\bm Q} \; {\bm y}({\bm \xi})\Big) \\
         \st & \quad {\bm x} \in \sets X, \; {\bm w} \in \sets W, \; {\bm y} \in \mathcal L_{N_\xi}^{N_{y}} \\
         & \quad \!\! \left. \begin{array}{l} 
         {\bm y}({\bm \xi}) \in \sets Y  \\
         {\bm T} ({\bm \xi}){\bm x} + {\bm V} ({\bm \xi}){\bm w} + {\bm W} ({\bm \xi}){\bm y}({\bm \xi}) \leq {\bm H}{\bm \xi} 
         \end{array} \quad \right\} \quad \forall {\bm \xi} \in \Xi \\
         & \quad {\bm y}({\bm \xi}) = {\bm y}({\bm \xi}') \quad \forall {\bm \xi}, \bm \xi' \in \Xi \; : \; {\bm w} \circ {\bm \xi} = {\bm w} \circ {\bm \xi}',
    \end{array}
\label{eq:endo_whole}
\end{equation}
where~$\boldsymbol{C} \in \mathbb{R}^{N_{\xi} \times N_{x}}, \; \boldsymbol{D} \in \mathbb{R}^{N_{\xi} \times N_{\xi}}, \; \boldsymbol{Q} \in \mathbb{R}^{N_{\xi} \times N_{y}}$, and~$\boldsymbol{H} \in \mathbb{R}^{L \times N_{\xi}}$. This problem optimizes over the here-and-now decisions~$\bm x \in \sets X$,~$\bm w \in \sets W$, and over the \newqj{generic} \emph{functional} decision rules $\bm y \in \mathcal L_{N_\xi}^{N_{y}}$, \newqj{which include any (measurable) function form, such as linear, piecewise linear, etc}. Modeling $\bm y$ as (measurable) functions of ${\bm \xi}$ ensures that these recourse decisions enjoy a lot of modeling flexibility, being able to adapt \emph{fully} to the (portion of) uncertain parameters that have been revealed between the first and second decision-stages. We assume that the left-hand-side of the constraints is linear in~$\bm \xi$, specifically,~$\; \boldsymbol{T}(\bm \xi):\Xi \xrightarrow{L} \mathbb{R}^{L \times N_x}, \; \boldsymbol{V}(\bm \xi):\Xi \xrightarrow{L} \mathbb{R}^{L \times N_{\xi}}, \; \boldsymbol{W}(\bm \xi):\Xi \xrightarrow{L} \mathbb{R}^{L \times N_y}$. We can account for affine dependencies on~$\bm \xi$ by introducing an auxiliary uncertain parameter~$ \xi_{N_\xi+1}$ restricted to equal one.

% The coefficient vectors~$[\bm T(\bm \xi)]_l,\;[\bm V(\bm \xi)]_l,\;[\bm W(\bm \xi)]_l$ in the~$l$-th constraint takes the form~$\bm \xi^\top \bm T_l,\;\bm \xi^\top \bm V_l$, and $\bm \xi^\top \bm W_l$, where~$\bm T_l \in \reals^{N_\xi \times N_x},\,\bm V_l \in \reals^{N_\xi \time N_w}$, and~$\bm W_l \in \reals^{N_\xi \times N_y}$.

The last set of constraints, which is a decision-dependent non-anticipativity constraint, enforces that~$\bm y(\cdot)$ can only depend on the observed uncertain parameters. Given two realizations~$\bm \xi$ and~$\bm \xi^\prime$, the equality~$\bm w \circ \bm \xi = \bm w \circ \bm \xi^\prime$ holds if and only if the observed portions of the two realizations are the same. In this case, the corresponding wait-and-see decisions~$\bm y(\bm \xi)$ and~$\bm y(\bm \xi^\prime)$ should also be consistent with one another. The uncertainty set~$\Xi$ in \eqref{eq:endo_whole} is a nonempty bounded polyhedron expressible as~$\Xi=\left\{\boldsymbol{\xi} \in \mathbb{R}^{N_{\xi}}: \bm A \bm{\xi} \leq \bm b\right\}$ for some matrix~$\boldsymbol{A} \in \mathbb{R}^{R \times N_{\xi}}$ and vector~$\boldsymbol{b} \in \mathbb{R}^{R}$. 

The ambiguity set~$\mathcal P$ takes the form
\begin{equation}
    \mathcal{P} = \Big\{\mathbb{P}\in\mathcal{M}_+(\mathbb{R}^{N_\xi})\; :\; \mathbb{P}(\bm \xi\in\Xi) = 1,\;\;\mathbb{E}_{\mathbb{P}}(\bm g(\bm \xi))\leq \bm c \Big\},
\label{eq:ambi_set}
\end{equation}
where $\mathcal{M}_{+}\left(\mathbb{R}^{N_{\xi}}\right)$ denotes the cone of nonnegative Borel measures supported on $\mathbb{R}^{N_{\xi}}$. The uncertainty set~$\Xi$ is also the support set, being defined as the smallest set that is known to satisfy $\boldsymbol{\xi} \in \Xi$ with probability 1 for some distribution $\mathbb{P}\in\mathcal{P}$. We assume that $\bm c \in \mathbb{R}^{N_g}$ and that $\bm g: \mathbb{R}^{N_{\xi}} \rightarrow \mathbb{R}^{N_g}$ is a convex piecewise linear multifunction defined as
$$
g_{s}(\boldsymbol{\xi})=\max _{t \in \mathcal{T}} \boldsymbol{g}_{s t}^{\top} \boldsymbol{\xi} \quad \forall s \in \mathcal{S}=\{1, \ldots, N_g\}.
$$
The vector~$\bm g_{st} \in \reals^{N_{\xi}}$ collects the coefficients of the~$t$-th piece of the function~$g_s(\bm \xi)$. Without loss of generality, the index $t$ of the linear pieces ranges over the same index set $\mathcal{T}=\{1, \ldots, T\}$ for each $g_s(\bm \xi)$, where $s \in \sets S$. Finally, we assume that the ambiguity set $\mathcal{P}$ contains a Slater point in the sense that there is a distribution $\mathbb{P} \in \mathcal{P}$ such that $\mathbb{E}_{\mathbb{P}}\left[ g_{s}({\boldsymbol{\xi}})\right]<c_s$ for all $s \in \mathcal{S}$ for which $ g_{s}(\boldsymbol{\xi})$ is nonlinear. 

In the DDID problem, since the information is discovered after implementing the measurement decision, we do not have empirical data to build the data-driven ambiguity set. For example, in the landslide prediction and warning problem, we may have information from other sites that could be used to obtain statistics on moments, but if we never installed sensors on the slope, no data was collected. As such, we consider using the moment-based ambiguity set~\eqref{eq:ambi_set}, which can be built judiciously using domain knowledge and expert assessments. The ambiguity set~\eqref{eq:ambi_set} has the flexibility to characterize many parameters of the unknown true distribution~$\mathbb P_0$, e.g., mean, mean absolute deviation. It can also approximate the nonlinear parameters, e.g., variance and standard deviation. \newqj{All of these parameters can be calibrated using domain expertise or historical data (e.g., data from related projects in R\&D project portfolio optimization can be used to calibrate the mean return and cost of each project as well as the mean deviation).} Another benefit of~\eqref{eq:ambi_set} is tractability, as it usually leads to mixed integer linear programming or mixed integer nonlinear programming formulations. 

{\remark When~$\bm w = \bm 0$, no uncertain parameters are revealed between the first and second decision stages. In this case, problem~\eqref{eq:endo_whole} reduces to a single-stage DRO problem, where the decision-maker chooses~$\bm x,\;\bm w,\;\bm y$ here-and-now to optimize in view of the worst-case distribution. When~$\bm w = \mathbf e$, the problem~\eqref{eq:endo_whole} reduces to a two-stage adaptive DRO problem. Finally, when the ambiguity set only contains the support information, we recover the robust optimization problem studied in~\citet{Vayanos2020}.}

Solving the \newqj{decision rule-}based formulation~\eqref{eq:endo_whole} is computationally challenging since it optimizes over a function space~$\mathcal L_{N_\xi}^{N_{y}}$. To the best of our knowledge, no approximation schemes for the adaptive binary variables can be directly applied. The prepartioning approach proposed in \citet{Vayanos2011DecisionProgramming} is tailored to stochastic and pure robust optimization problems and does not generalize to cases where the distribution is unknown. The approaches proposed in \citet{han2023finite} and \citet{Hanasusanto2016} can only solve the special instance of~\eqref{eq:endo_whole} where~$\bm w = \mathbf e$. The lack of approximation schemes applicable to formulation~\eqref{eq:endo_whole} motivates us to derive an equivalent nested reformulation that gives us access to a broader set of approximation and solution schemes.

\subsection{An Equivalent min-max-min-max Reformulation}
\label{subsec:nested-formulation}
In this section, we derive an equivalent min-max-min-max reformulation of problem~\eqref{eq:endo_whole} that yields the same set of optimal first-stage decisions~$(\bm x,\bm w)$ and corresponding optimal objective values. We formalize the equivalence in the following theorem.
\begin{theorem}
The objective value and set of optimal solutions~$(\bm x,\bm w)$ to problem~\eqref{eq:endo_whole} are equal to those of the min-max-min-max problem
\addtocounter{equation}{+1}
\begin{equation}
    % \tag{$\sets N$}
    \begin{array}{cl}
         \min & \;\; \displaystyle \max_{\overline{\bm \xi} \in \Xi} \displaystyle \left\{ 
         \begin{array}{cl}
            \displaystyle\min_{ {\bm y} \in \sets Y } \;\;
             &\displaystyle\max_{ {\bm \xi} \in \Xi({\bm w},\overline{\bm \xi}) } \; \bm c^\top \bm \psi + {\bm \xi}^\top {\bm C} \; {\bm x} + {\bm \xi}^\top {\bm D} \; {\bm w}  + {\bm \xi}^\top {\bm Q} \; {\bm y} - \bm\psi^\top \bm g(\bm\xi) \\
             \st & \;\; {\bm T} ({\bm \xi}){\bm x} + {\bm V} ({\bm \xi}){\bm w} + {\bm W} ({\bm \xi}){\bm y} \leq {\bm H}{\bm \xi} \;\;\; \forall {\bm \xi} \in \Xi({\bm w},\overline{\bm \xi})
         \end{array}
          \right\}  \\
         \st & \;\; \bm\psi\in\mathbb{R}_+^{N_g},\;{\bm x} \in \sets X,\; {\bm w} \in \sets W,
    \end{array}\label{eq:nested}
\end{equation}
\label{thm:min-max}
where~$\Xi(\boldsymbol{w}, \overline{\boldsymbol{\xi}}):=\{\boldsymbol{\xi} \in \Xi: \boldsymbol{w} \circ \boldsymbol{\xi}=\boldsymbol{w} \circ \overline{\boldsymbol{\xi}}\}.$
\end{theorem}

Problem~\eqref{eq:nested} can be interpreted as a sequential game against ``nature''. In the first stage, the decision-maker first decides on $\bm \psi,\,\bm x,\,\text{and }\bm w$. Then, nature chooses a realization~$\overline {\bm \xi}$ of the uncertain parameters, of which the decision-maker can only observe those elements $\xi_i$ such that $w_i=1$. In the second stage, the decision-maker selects ${\bm y}$ in a way that is robust to the remaining elements of $\xi_i$ (that have not been observed). Finally, nature selects a realization of ${\bm \xi}$ that coincides with $\overline{\bm \xi}$ for those elements that have been observed between the first and second decision stages, i.e., such that $\xi_i = \overline \xi_i$ for $i \in \mathcal I$ such that $w_i=1$.

The sketch of the proof \newqj{of} Theorem~\ref{thm:min-max} is as follows. For fixed~$\bm w, \;\bm x, \;\bm y(\cdot)$ we identify the inner maximization of problem~\eqref{eq:endo_whole} as a moment problem and dualize it.
% and combine the dualized problem with the outer minimization. 
Strong duality guarantees that the optimal solution~$(\bm x,\;\bm w,\; \bm y(\cdot))$ in the dualized problem is also optimal in the original problem~\eqref{eq:endo_whole} and that both problems have the same objective value. The dual problem then reduces to a robust optimization problem with decision-dependent information discovery. Applying \citet{Vayanos2020} 
Theorem 1 yields the equivalent nested formulation~\eqref{eq:nested}. The full proof is deferred to Appendix~\ref{ec:proofs}.

%%%%%%%%%%%%%%%%%%%%%%%%%%%%%%%%%%%%%%%%%%%%%%%%%%%%%%%%%%%%%%%%%%%%%%%%%%%%%%%%%%%%%%%%%%%%%%%%
%%%%%%%%%%%%%%%%%%%%%%%%%%%%%%%%%%%%%%%%%%%%%%%%%%%%%%%%%%%%%%%%%%%%%%%%%%%%%%%%%%%%%%%%%%%%%%%%
%%%%%%%%%%%%%%%%%%%%%%%%%%%%%%%% K-Adaptability counterpart %%%%%%%%%%%%%%%%%%%%%%%%%%%%%%%%%%%
%%%%%%%%%%%%%%%%%%%%%%%%%%%%%%%%%%%%%%%%%%%%%%%%%%%%%%%%%%%%%%%%%%%%%%%%%%%%%%%%%%%%%%%%%%%%%%%%
%%%%%%%%%%%%%%%%%%%%%%%%%%%%%%%%%%%%%%%%%%%%%%%%%%%%%%%%%%%%%%%%%%%%%%%%%%%%%%%%%%%%%%%%%%%%%%%%
\section{\texorpdfstring{$K$-Adaptability Counterpart}{}}
\label{sec:k-adapt}
The min-max-min-max reformulation~\eqref{eq:nested} enables us to propose a new approximate and potentially exact solution approach to solve the DRO problem with DDID. The new solution approach is based on the~$K$-adaptability approximation, which allows us to control the trade-off between complexity and solution quality. In the~$K$-adaptability approximation,~$K$ different policies are chosen here-and-now before any uncertain parameters are revealed. After a portion of uncertain parameters~$\{\xi_i:w_i=1\}$ are observed, the best candidate policy is implemented among those that are robustly feasible considering the uncertain parameters that remain unknown. If there exists some~$\bm \xi \in \Xi$ such that no policies are feasible, the~$K$-adaptability problem evaluates to~$+\infty$. While in problem~\eqref{eq:nested}, the wait-and-see decisions~$\bm y$ can take on any value that satisfies the constraint in~$\mathcal{Y}$, the~$K$-adaptability approximation only allows us to choose one of the $K$ candidate values identified in the first stage. It thus results in a conservative approximation to problem~\eqref{eq:nested}.

The~$K$-adaptability counterpart of problem~\eqref{eq:nested} can be written as
\begin{equation}
% \tag{$\sets N_K$}
    \begin{array}{cl}
         \min & \;\; \displaystyle \max_{\overline{\bm \xi} \in \Xi} \displaystyle \left\{ 
         \begin{array}{cl}
            \displaystyle\min_{ k \in \sets K } \;\;
             &\displaystyle\max_{ {\bm \xi} \in \Xi({\bm w},\overline{\bm \xi}) } \; \bm c^\top \bm \psi + {\bm \xi}^\top {\bm C} \; {\bm x} + {\bm \xi}^\top {\bm D} \; {\bm w}  + {\bm \xi}^\top {\bm Q} \; {\bm y^k} - \bm\psi^\top \bm g(\bm\xi) \\
             \st & \;\; {\bm T} ({\bm \xi}){\bm x} + {\bm V} ({\bm \xi}){\bm w} + {\bm W} ({\bm \xi}){\bm y^k} \leq {\bm H}{\bm \xi} \;\;\; \forall {\bm \xi} \in \Xi({\bm w},\overline{\bm \xi})
         \end{array}
          \right\}  \\
         \st & \;\; \bm\psi\in\mathbb{R}_+^{N_g},\;{\bm x} \in \sets X,\; {\bm w} \in \sets W, \; \bm y^k \in\mathcal{Y},\; \forall k\in\sets K,
    \end{array}\label{eq:k-nested}
\end{equation}
where~$\sets K = \{1,\ldots,K\}$. 

\newqj{
\remark{
The~$K$-adaptability counterpart~\eqref{eq:k-nested} has the following equivalent \newqj{decision rule} formulation
\begin{equation}
    \begin{array}{cl}
         \min & \quad \displaystyle \sup_{\mathbb{P} \in \mathcal{P}} \;\;\;\mathbb{E}_{\mathbb{P}}\Big(  {\bm \xi}^\top {\bm C} \; {\bm x} + {\bm \xi}^\top {\bm D} \; {\bm w}  + {\bm \xi}^\top {\bm Q} \; {\bm y}({\bm \xi})\Big) \\
         \st & \quad {\bm x} \in \sets X, \; {\bm w} \in \sets W, \; {\bm y} \in \mathcal L_{N_\xi}^{N_{y}}, \bm y^k \in \sets Y, \forall k \in \sets K \\
         & \quad \!\! \left. \begin{array}{l} 
         {\bm y}({\bm \xi}) \in \{\bm y^1,\,\ldots,\bm y^K\}  \\
         {\bm T} ({\bm \xi}){\bm x} + {\bm V} ({\bm \xi}){\bm w} + {\bm W} ({\bm \xi}){\bm y}({\bm \xi}) \leq {\bm H}{\bm \xi} 
         \end{array} \quad \right\} \quad \forall {\bm \xi} \in \Xi \\
         & \quad {\bm y}({\bm \xi}) = {\bm y}({\bm \xi}') \quad \forall {\bm \xi}, \bm \xi' \in \Xi \; : \; {\bm w} \circ {\bm \xi} = {\bm w} \circ {\bm \xi}'.
    \end{array}
\label{eq:k-adapt-dr-formulation}
\end{equation}
}
}
%%%%%%%%%%%%%%%%%%%%%%%%%%%%%%%%%%%%%%%%%%%%%%%%%%%%%%%%%%%%%%%%%%%%%%%%%%%%%%%%%%%%%%%%%%%%%%%%
%%%%%%%%%%%%%%%%%%%%%%%%%%%%%%%%%%%%%%%%%%%%%%%%%%%%%%%%%%%%%%%%%%%%%%%%%%%%%%%%%%%%%%%%%%%%%%%%
%%%%%%%%%%%%%%%%%%%%%%%%%%%%%%%% Mixed Binary Bilinear Formulation %%%%%%%%%%%%%%%%%%%%%%%%%%%%%%%%%%%
%%%%%%%%%%%%%%%%%%%%%%%%%%%%%%%%%%%%%%%%%%%%%%%%%%%%%%%%%%%%%%%%%%%%%%%%%%%%%%%%%%%%%%%%%%%%%%%%
%%%%%%%%%%%%%%%%%%%%%%%%%%%%%%%%%%%%%%%%%%%%%%%%%%%%%%%%%%%%%%%%%%%%%%%%%%%%%%%%%%%%%%%%%%%%%%%%
% \subsection{Mixed Binary Bilinear Formulation}
% \label{bilinear formulation}
We provide reformulations of the~$K$-adaptability problem~\eqref{eq:k-nested} that can be fed into off-the-shelf solvers, beginning with the case involving only objective uncertainty.  In the absence of uncertainty in the constraints, the problem is expressible as
\begin{equation}
    \begin{array}{cl}
         \min & \quad \displaystyle \max_{\overline{\bm \xi} \in \Xi} \;\;
            \displaystyle\min_{ { k} \in \sets K }
            \left\{
            \displaystyle\max_{ {\bm \xi} \in \Xi({\bm w},\overline{\bm \xi}) } \; \bm c^\top \bm \psi + {\bm \xi}^\top {\bm C} \; {\bm x} + {\bm \xi}^\top {\bm D} \; {\bm w}  + {\bm \xi}^\top {\bm Q} \; {\bm y^k} - \bm\psi^\top \bm g(\bm\xi) 
            \right\} \\
         \st & \quad \bm\psi\in\mathbb{R}_+^{N_g},\;{\bm x} \in \sets X,\; {\bm w} \in \sets W, \; \bm y^k \in\mathcal{Y},\; \forall k\in\sets K \\
         & \quad {\bm T} {\bm x} + {\bm V} {\bm w} + {\bm W} {\bm y^k} \leq {\bm h},\; \forall k\in\sets K,
    \end{array}\label{eq:k-objunc}
\end{equation}
where~$\bm h \in \reals^L$. We show the equivalent mixed binary bilinear reformulation of problem~\eqref{eq:k-objunc} in the following theorem. 

\begin{theorem}
The~$K$-adaptability problem with objective uncertainty~\eqref{eq:k-objunc} is equivalent to the following mixed binary bilinear problem
\begin{equation}
    \begin{array}{cl}
         \min & \quad  \bm c^\top \bm \psi + {\bm b}^\top {\bm \beta}  + \sum_{k\in \mathcal K} {\bm b}^\top {\bm \beta}^k\\
         \st & \quad {\bm x} \in \sets X, \; {\bm w} \in \sets W, \; {\bm y}^k \in \sets Y, \; \forall k \in \sets K \\
         & \quad {\bm \alpha} \in \mathbb R^K_+, \; {\bm \beta} \in \mathbb R^R_+, \; {\bm \beta}^k \in \mathbb R^R_+, \; \forall k \in \sets K, \; {\bm \gamma}^k \in \mathbb R^{N_\xi}, \; \forall k \in \mathcal K \\
        & \quad \bm\psi\in\mathbb{R}_+^{N_g},\, \delta_{s,t}^k\in \mathbb R_+, \; \forall s \in \sets S, \, \forall t\in\sets T,\;\forall k\in\sets K\\
        & \displaystyle\quad    {\bm A}^\top{\bm \beta}^k +  {\bm w} \circ {\bm \gamma}^k + \sum_{s\in \sets S}\sum_{t\in \sets T} \delta_{s,t}^k \bm g_{s,t} = {\alpha}^k \left( {\bm C} {\bm x} + {\bm D} {\bm w} + {\bm Q} {\bm y}^k \right) \quad \forall k \in \mathcal K \\
        & \displaystyle\quad \sum_{t\in\sets T}\bm \delta_{t}^k = \alpha^k \bm  \psi \quad \forall k\in\sets K\\
        & \quad \1^\top {\bm \alpha}  = 1  \\
        & \quad  {\bm A}^\top{\bm \beta}  =  \displaystyle \sum_{k \in \mathcal K} {\bm w} \circ  {\bm \gamma}^k \\
         & \quad {\bm T} {\bm x} + {\bm V} {\bm w} + {\bm W}{\bm y}^k \leq {\bm h} \quad \forall k \in \sets K. 
    \end{array}
    \label{eq:objunc-bilinear}
\end{equation}
\label{thm:objunc-bilinear}
\end{theorem}

For the case of constraint uncertainty, to keep the paper concise, we only present the main theorem and defer the full mixed integer nonlinear optimization (MINLO) reformulation to Appendix~\ref{ec:proofs}. The dual variables~${\bm \alpha}$ are associated with the set of epigraph constraints in problem~\eqref{subproblem1}. The variables ${\bm \beta}$ and $\{{\bm \beta}^k\}_{k \in \sets K}$ are associated with the constraints~$\bm A \bar{\bm \xi} \leq \bm b$ and~$\bm A \bm \xi^k \leq \bm b,\; \forall k \in \sets K$, respectively. The dual variables $\{{\bm \gamma}^k\}_{k \in \mathcal K}$ are associated with the constraint~$\bm w \circ \bm \xi^k = \bm w \circ \bar{\bm \xi},\; \forall k \in \sets K$. Finally, $\{\delta_{s,t}^k\}_{s \in \sets S, t \in \sets T, k \in \sets K}$ arises from dualizing the set of constraints~$\zeta_s^k \leq -\bm g_{s,t}^\top\bm \xi^k$ with respect to the moment information.
\begin{theorem}
\label{thm:cstrunc-bilinear}
The~$K$-adaptability problem~\eqref{eq:k-nested} with constraint uncertainty has an equivalent mixed binary bilinear reformulation.
\end{theorem}

In the case with only objective uncertainty, formulation \eqref{eq:objunc-bilinear} presents two sets of bilinear terms. The first set of constraints contain~$\alpha^k\bm x,\,\alpha^k\bm w,$ and~$\alpha^k \bm y^k$, which can be linearized if~$\bm x$ and~$\bm y^k$ are discrete variables. The second set of constraints involves the term~$\alpha^k\bm \psi$ which is a product between continuous variables (recall that~$\bm \psi$ is the dual variable of the constraint~$\mathbb E_{\mathbb P}(\bm g(\bm \xi)) \leq \bm c$ in the ambiguity set~\eqref{eq:ambi_set}). The number of nonlinear terms scale\newqj{s} as~$O( K \cdot N_g )$. In cases where problem~\eqref{eq:k-nested} involves constraint uncertainty, its MINLO reformulation~\eqref{eq:cstrunc-bilinear} also includes bilinear terms between continuous variables, scaling as~$O( (1 + L)^K \cdot K \cdot N_g )$. These terms do not appear if the ambiguity set~\eqref{eq:ambi_set} only has support information.}

Solving a mixed-integer bilinear problem to global optimality is difficult in general. Since version 9.0, Gurobi \newqj{has} added support for solving nonconvex MINLO problems using a spatial branching algorithm. The generic algorithm in Gurobi converges slowly when the number of bilinear terms is large (see Table~\ref{table:r&d-per}). Instead of solving the monolithic MINLO problems~\eqref{eq:objunc-bilinear} and~\eqref{eq:cstrunc-bilinear} in the case of objective uncertainty and constraint uncertainty, respectively, in the following section, we propose a decomposition algorithm that solves problem~\eqref{eq:k-nested} using a nested decomposition algorithm. 

%%%%%%%%%%%%%%%%%%%%%%%%%%%%%%%%%%%%%%%%%%%%%%%%%%%%%%%%%%%%%%%%%%%%%%%%%%%%%%%%%%%%%%%%%%%%%%%%
%%%%%%%%%%%%%%%%%%%%%%%%%%%%%%%%%%%%%%%%%%%%%%%%%%%%%%%%%%%%%%%%%%%%%%%%%%%%%%%%%%%%%%%%%%%%%%%%
%%%%%%%%%%%%%%%%%%%%%%%%%%%%%%%% L-shaped algorithm %%%%%%%%%%%%%%%%%%%%%%%%%%%%%%%%%%%
%%%%%%%%%%%%%%%%%%%%%%%%%%%%%%%%%%%%%%%%%%%%%%%%%%%%%%%%%%%%%%%%%%%%%%%%%%%%%%%%%%%%%%%%%%%%%%%%
%%%%%%%%%%%%%%%%%%%%%%%%%%%%%%%%%%%%%%%%%%%%%%%%%%%%%%%%%%%%%%%%%%%%%%%%%%%%%%%%%%%%%%%%%%%%%%%%

\section{Decomposition Algorithm}
\label{sec:l-shaped}

In this section, we leverage the nested structure of problem~\eqref{eq:k-nested} to design a two-layer decomposition algorithm that solves it to global optimality. A decomposition algorithm has been used to solve a two-stage robust optimization problem with DDID and objective uncertainty in \citet{paradiso2022exact}. We extend it to the DRO setting with both objective and constraint uncertainty. Figure \ref{fig:overview} provides an overview of our overall proposed solution procedure. We present an~$L$-shaped algorithm designed to optimize the measurement variables~$\bm w$, as introduced in Section~\ref{subsec:l-shaped-outter}. A key subroutine within the~$L$-shaped algorithm involves the evaluation of the true objective value of \eqref{eq:k-nested} for a given $\bm w$. We propose to address this evaluation problem through a branch-and-cut algorithm, the details of which are discussed in Section~\ref{subsec:b&c}.  In addition, we develop three cut families to strengthen the main problem in the~$L$-shaped algorithm. We discuss them in Section~\ref{subsec:tighter-cuts}. 
\begin{figure}[htbp!]
\centering
\tikzstyle{startstop} = [rectangle, rounded corners, minimum width=3cm, minimum height=1cm,text centered, draw=black]
\tikzstyle{io} = [trapezium, trapezium left angle=70, trapezium right angle=110, minimum width=3cm, minimum height=1cm, text centered, draw=black]
\tikzstyle{process} = [rectangle, minimum width=3cm, minimum height=1cm, text centered, draw=black, execute at begin node=\setlength{\baselineskip}{3.5ex}]
\tikzstyle{decision} = [diamond, minimum width=3cm, minimum height=1cm, text centered, draw=black]

\tikzstyle{arrow} = [thick,->,>=stealth]

\scalebox{0.87}{
\begin{tikzpicture}[align=center,scale=5]
\small
\node (start) [startstop] {Start};

\node (solve-out) [process, below of=start, node distance=1.75cm] {Solve relaxation \\problem \eqref{eq:mins-w-master} from Section~\ref{subsec:l-shaped-outter}};
\draw [arrow] (start) -- (solve-out);

\node (eva) [process, below of=solve-out, node distance=2.5cm] {Evaluate~$\Phi(\bm w^\prime)$ using\\ 
Algorithm~\ref{alg:bnb} from Section~\ref{subsubsec:b&b} or\\
Algorithm~\ref{alg:bnc} from Section~\ref{subsubsec:b&c} };

\draw [arrow] (solve-out) -- node[anchor=west, align=left] {optimal~$\bm w^\prime$, \\ optimal objective value $\theta^\prime$} (eva);

\node (check) [decision, below of=eva, node distance=2.1cm, aspect=2, align=center] {$|\Phi(\bm w^\prime)- \theta^\prime| \leq \varepsilon$?};

\draw [arrow] (eva) -- (check);

\node (report) [startstop, below of=check, node distance=2cm] 
{Report optimal solution \\ $\bm w^\star,\;\left(\bm \psi^\star,\;\bm x^\star,\;\bm \left\{\bm y^{k^\star}\right\}_{k\in\sets K}\right)$};
\draw [arrow] (check) -- node[anchor=west] {Yes} (report);

\node (cut) [process, left of=eva, node distance=6.7cm]
{Add cut \eqref{eq:integer-cut} or \eqref{eq:feasible-cut} from Section~\ref{subsec:l-shaped-outter}, \\ and cuts \eqref{eq:benders-cut}, \eqref{eq:strengthen-feas},  \\ \eqref{eq:information-cut} from Section~\ref{subsec:tighter-cuts} if applicable};

\coordinate[left of = check, node distance = 6.7cm] (aux0);
\draw [arrow] (check) -- node[anchor=south] {No} (aux0) -- (cut);

\draw [arrow] (cut) |- (solve-out);
\end{tikzpicture}
}
\caption{Overview of our proposed decomposition Algorithm \ref{alg:l-shaped}}
\label{fig:overview}
\end{figure}

\subsection{\texorpdfstring{$L$-shaped Algorithm to Optimize~$\bm w$
}{}}
\label{subsec:l-shaped-outter}

The objective function of~\eqref{eq:k-nested} can be written as a function of~$\bm w$. We present in this section an algorithm to solve the decomposed problem
\begin{equation}
    \begin{array}{cl}
         \displaystyle\min & \quad \Phi(\bm w) \\
         \st & \quad {\bm w} \in \sets W,
    \end{array}\label{eq:mins-w}
\tag{$\mathcal P$}
\end{equation}
where
\begin{equation}
    \begin{array}{ccl}
         \Phi(\bm w) = & \min & \displaystyle \bm c^\top \bm \psi +\max_{\overline{\bm \xi} \in \Xi} \displaystyle \left\{ 
         \begin{array}{cl}
            \displaystyle\min_{ k \in \sets K } 
             &\displaystyle\max_{ {\bm \xi} \in \Xi({\bm w},\overline{\bm \xi}) }   {\bm \xi}^\top {\bm C}  {\bm x} + {\bm \xi}^\top {\bm D}  {\bm w}  + {\bm \xi}^\top {\bm Q}  {\bm y^k} - \bm\psi^\top \bm g(\bm\xi) \\
             \st &  {\bm T} ({\bm \xi}){\bm x} + {\bm V} ({\bm \xi}){\bm w} + {\bm W} ({\bm \xi}){\bm y^k} \leq {\bm H}{\bm \xi} \;\;\; \forall {\bm \xi} \in \Xi({\bm w},\overline{\bm \xi})
         \end{array}
          \right\}  \\
         & \;\st & \;\; \bm\psi\in\mathbb{R}_+^{N_g},\;{\bm x} \in \sets X,\; \bm y^k \in\mathcal{Y}, \;\forall k\in\sets K.
    \end{array}\label{eq:Phi(w)}
\tag{$\mathcal {IP}_{\bm w}$}
\end{equation}
The function~$\Phi(\cdot)$ is, in general, nonconvex on~$[0, 1]^{N_{\xi}}$. However, since~$\bm w$ is binary,~\citet[Proposition 1]{Laporte1993TheRecourse} show that~$\Phi(\cdot)$ can be replaced by its convex lower envelope, a convex piecewise linear function on~$[0,1]^{N_{\xi}}$, that coincides with~$\Phi(\cdot)$ at all binary points in~$\{0,1\}^{N_{\xi}}$. We construct the convex lower envelope using two sets of constraints.

The first set of constraints contains the integer optimality inequalities first proposed by~\citet{Laporte1993TheRecourse} for a stochastic problem with integer first-stage variables, defined as
\begin{equation}
\Phi(\bm w) \geq\left(\Phi(\bm w^{r})-L\right)\left(\sum_{i \in \sets I_{r}} w_{i}-\sum_{i \notin \sets I_{r}} w_{i}\right)-\left(\Phi(\bm w^{r})-L\right)\left(\left|\sets I_{r}\right|-1\right)+L \quad \forall \bm w^r \in \sets W,
\tag{$\mathcal C_\text{opt}$}
\label{eq:integer-cut}
\end{equation}
where~$L$ is any lower bound on~$\Phi(\bm w)$.~$\sets I_r$ collects the indices of the components of vector~$\bm w^r$ that take a value \newqj{of} 1. We can verify that the integer optimality cuts are both valid and tight at each binary point~$\bm x$. When $\bm w = \bm w^r$, the term $\sum_{i \in \sets I_{r}} w_{i} - \sum_{i \notin \sets I_{r}} w_{i}$ is equal to $|\sets I_r|$, and the right-hand side of the inequality \eqref{eq:integer-cut} is equal to $\Phi(\bm w^r)$. In the case where $\bm w \neq \bm w^r$, the value of the term $\sum_{i \in \sets I_{r}} w_{i} - \sum_{i \notin \sets I_{r}} w_{i}$ is less than $|\sets I_r|$, thereby ensuring that the right-hand side of the inequality \eqref{eq:integer-cut} does not exceed $L$. The lower bound~$L$ can be determined by the problem structure or obtained heuristically, for example, by solving a deterministic problem which replaces the uncertainty sets in both the inner and outer maximization problems with a singleton~$\{\bar{\bm \xi}\}$, where~$\bar{\bm \xi} \in \Xi$.

For points~$\bm w^r$ such that~$\Phi(\bm w^r) = +\infty$, we can either set~$\Phi(\bm w^r)$ as a sufficiently large number in inequalities~\eqref{eq:integer-cut}, or exclude them from the feasible region~$\sets W$ by enforcing the following set of feasibility cuts (see~\citet{Balas1972CanonicalHypercube} and~\citet{Nannicini2012Rounding-basedMINLPs})
\begin{equation}
\sum_{i \in \sets I_r} \left(1-w_i\right) + \sum_{i\notin \sets I_r} w_i\geq 1, \quad \forall \bm w^r \in \sets W_\text{inf},
\tag{$\mathcal C_\text{feas}$}
\label{eq:feasible-cut}
\end{equation}
where~$\sets W_\text{inf}=\{\bm w: \bm w \in \sets W,\; \Phi(\bm w) = +\infty \}$. The cuts state that any feasible points~$\bm w$ should be different on at least one component from all infeasible points~$\bm w^r \in \sets W_\text{inf}$. Expressing the feasibility cut separately in~\eqref{eq:feasible-cut}, we can then exclude them from~\eqref{eq:integer-cut} and define the set of optimality inequalities on the set~$\sets W_\text{feas} = \sets W \setminus \sets W_\text{inf}.$ 

For any~$\widehat{\sets W}_\text{feas} \subseteq \sets W_\text{feas}$ and $\;\widehat{\sets W}_\text{inf} \subseteq \sets W_\text{inf},$ the following problem \eqref{eq:mins-w-master} is a relaxation of problem~\eqref{eq:mins-w} since its constraints enforce a subset of the feasibility cuts~\eqref{eq:feasible-cut}, and the objective takes an epigraph form on a subset of the optimality cuts~\eqref{eq:integer-cut}.
\begin{equation}
    \begin{array}{cl}
        \displaystyle \min & \;\; \theta \\
         \st & \;\; \theta \in \mathbb R,\;{\bm w} \in \sets W \\
         & \;\; 
              \theta \geq \displaystyle \left(\Phi(\bm w^{r})-L\right)\left(\sum_{i \in \sets I_{r}} w_{i}-\sum_{i \notin \sets I_{r}} w_{i}\right)
              \begin{array}[t]{r}-\left(\Phi(\bm w^{r})-L\right)\left(\left|\sets I_{r}\right|-1\right)+L \\
              \forall \bm w^r \in \widehat{\sets W}_\text{feas} 
         \end{array}\\
         & \;\;\displaystyle{\sum_{i \in \sets I_r} \left(1-w_i\right) + \sum_{i\notin \sets I_r} w_i\geq 1 \quad \forall \bm w^r \in \widehat{\sets W}_\text{inf}}.
    \end{array}\label{eq:mins-w-master}
\tag{$\mathcal P_\text{relax}$}
\end{equation}

We use the relaxed problem \eqref{eq:mins-w-master} as the main problem of the~$L$-shaped algorithm. Algorithm~\ref{alg:l-shaped} describes the~$L$-shaped procedure. It iteratively solves problem~\eqref{eq:mins-w-master}, evaluates the value~$\Phi(\bm w^\prime)$ at the current optimal solution~$\bm w^\prime$, adds new cuts, and updates the bounds. The proof of finite convergence in~\citet[Proposition 1]{Laporte1993TheRecourse} for a two-stage integer stochastic problem can be easily adapted. To keep the paper concise, we omit it here.

\begin{algorithm}[htb]
\setstretch{1.35}
\SetAlgoLined
\textbf{Inputs:} Any lower bound~$L$ of~$\theta$, optimality gap tolerance~$\varepsilon$.

\textbf{Output:} Optimal solution~$\bm w^\star$, and optimal objective value~$\theta^\star$.

\textbf{Initialization:} Set~$\theta^\star \leftarrow +\infty,\;\bm w^* \leftarrow \emptyset$,~$\widehat{\sets W}_\text{feas} \leftarrow \emptyset$,~$\widehat{\sets W}_\text{inf} \leftarrow \emptyset$. 

\begin{enumerate}[rightmargin=\dimexpr\linewidth-11.25cm-\leftmargin\relax, leftmargin=1.3em]
\item Solve problem~\eqref{eq:mins-w-master}, let~$(\bm w^\prime,\;\theta^\prime)$ be an optimal solution. If $|\theta^\star - \theta^\prime| \leq \varepsilon$, stop.

\item Evaluate~$\Phi(\bm w^\prime) $ by solving problem~\eqref{eq:Phi(w)}.

\item If $\Phi(\bm w^\prime) = +\infty$, set~$\widehat{\sets W}_\text{inf} \leftarrow \widehat{\sets W}_\text{inf} \bigcup \{\bm w^\prime\}$, go to step 1; else, go to step 4.

\item
If $\Phi(\bm w^\prime) < \theta^\star$, set~$\bm w^\star \leftarrow \bm w^\prime$ and~$\theta^\star \leftarrow \Phi(\bm w^\prime)$. Set~$\widehat{\sets W}_\text{feas} \leftarrow \widehat{\sets W}_\text{feas} \bigcup \{\bm w^\prime\}$, go to step~1.
\end{enumerate}

\textbf{Result:} Return~$\bm w^\star$ and~$\theta^\star$.
\caption{$L$-shaped algorithm for solving problem~\eqref{eq:mins-w}}
\label{alg:l-shaped}
\end{algorithm}
% \noteqj{I think the Algorithm \ref{alg:l-shaped} is correct. $\theta^\prime$ is the lower bound and~$\theta^\star$ is the current incumbent. In step 1, in each iteration, check if LB is greater than UB; if yes, then stop. In step 3, if $\Phi(\bm x^\prime)$ is better than the current incumbent, update the incumbent.}

\subsection{\texorpdfstring{Branch-and-Cut Algorithm to Evaluate~$\Phi(\bm w)$}{}}
\label{subsec:b&c}
We present in this section a branch-and-cut algorithm for evaluating the objective function~$\Phi(\bm w^\prime)$ at the current solution $\bm w^\prime$, which is a crucial component of Algorithm~\ref{alg:l-shaped}. We begin the section by discussing a basic branch-and-bound algorithm adapted from~\citet{Subramanyam2020} to find an optimal solution for a min-max-min problem with~$K$-adaptability. Then, we design a cutting plane algorithm to solve the inner min-max problem and merge it into the branch-and-bound tree to form the final branch-and-cut algorithm. Finally, we prove the convergence of the branch-and-cut algorithm.

Before proceeding with the algorithm, we introduce an equivalent formulation of problem~\eqref{eq:Phi(w)} as given in the following proposition.
\begin{proposition}
The following problem is equivalent to problem~\eqref{eq:Phi(w)}.
\begin{equation}
    \begin{array}{cl}
         \min & \;\; \theta \\
         \st & \;\; \theta\in\mathbb{R}, \; \bm\psi\in\mathbb{R}_+^{N_g}, \; {\bm x} \in \sets X, \; {\bm y}^k \in \sets Y \ \  \forall k \in \sets K, \; \Xi^k \subseteq \reals^{N_{\xi}} \;\; \forall k \in \sets K\\
         & \;\; \bigcup_{k \in \sets K} \Xi^k = \Xi \\
         & \;\; \left. \begin{array}{l} 
        \theta \geq \bm c^\top \bm \psi + \displaystyle\max_{ {\bm \xi} \in \Xi({\bm w ^ \prime},\overline{\bm \xi}) } \; \; 
        \begin{array}[t]{r}
        {\bm \xi}^\top {\bm C}  {\bm x} + {\bm \xi}^\top {\bm D}  {\bm w ^ \prime}  + {\bm \xi}^\top {\bm Q}  {\bm y^k} \\
        - \bm\psi^\top \bm g(\bm\xi)
        \end{array}\\
        {\bm T} ({\bm \xi}){\bm x} + {\bm V} ({\bm \xi}){\bm w ^ \prime} + {\bm W} ({\bm \xi}){\bm y^k} \leq {\bm H}{\bm \xi} \;\;\; \forall {\bm \xi} \in \Xi({\bm w ^ \prime},\overline{\bm \xi})
         \end{array}  \right \}
         \begin{array}{l}
              \forall \overline{\bm \xi} \in \Xi^k,\\
              \forall k \in \sets K,
         \end{array}
    \end{array}
\label{eq:disjuctive}
\tag{$\mathcal {IP}_{\bm w}({\Xi})$}
\end{equation}
where, for each~$k \in \sets K$, $\Xi^k$ denotes the set of scenarios for which the corresponding optimal policy is~$\bm y^k$.
\label{prop: disjunctive}
\end{proposition}

In~problem \eqref{eq:Phi(w)}, the best candidate policy~$\bm y^k$ is selected for each realization~$\overline{\bm \xi}$. Problem \eqref{eq:disjuctive} achieves the same goal by optimizing the decision variables~$\Xi^k,\; k \in \sets K$. This is equivalent to finding an optimal partition of the uncertainty set~$\Xi$. Problem \eqref{eq:disjuctive} is a robust optimization problem with a decision-dependent uncertainty set, since the uncertainty sets~$\{\Xi^k\}_{k\in \sets K}$ in the constraints are decision variables themselves.

\subsubsection{Constraint Generation Algorithm} 
\label{subsubsec:b&b}
To solve problem~\eqref{eq:disjuctive}, we use a constraint generation algorithm. The main problem~$\mathcal {IP}_{\bm w}(\widehat{\Xi})$ is a relaxation of~\eqref{eq:disjuctive} where the set~$\Xi$ in the partition constraint is substituted with a finite set~$\widehat{\Xi} \subseteq \Xi$. In each iteration, a scenario is generated and added to~$\widehat{\Xi}$.

We now describe our proposed procedure for solving the main problem. To this end, we present the following fixed partition problem, where~$\{\widehat{\Xi}^k\}_{k\in \sets K}$ are parameter inputs rather than decision variables.  
\begin{equation}
    \begin{array}{ccl}
        & \min & \theta \\
         & \st &  \theta\in\mathbb{R}, \; \bm\psi\in\mathbb{R}_+^{N_g}, \; {\bm x} \in \sets X, \; \bm y^k \in\mathcal{Y} \quad \forall k\in\sets K\\
         &&\left. \!\!\begin{array}{l} 
         \theta \geq \bm c^\top \bm \psi + \displaystyle\max_{ {\bm \xi} \in \Xi({\bm w^\prime},\overline{\bm \xi}) } \; \;
         \begin{array}[t]{r}
         {\bm \xi}^\top {\bm C} \; {\bm x} + {\bm \xi}^\top {\bm D} \; {\bm w^\prime}  + {\bm \xi}^\top {\bm Q} \; {\bm y^k} \\
         - \bm\psi^\top \bm g(\bm\xi) 
         \end{array}
         \\
         {\bm T} ({\bm \xi}){\bm x} + {\bm V} ({\bm \xi}){\bm w^\prime} + {\bm W} ({\bm \xi}){\bm y^k} \leq {\bm H}{\bm \xi} \quad \forall {\bm \xi} \in \Xi({\bm w^\prime},\overline{\bm \xi})
         \end{array} \; \right\} 
         \begin{array}{l}
              \overline{\bm \xi} \in \widehat{\Xi}^k,  \\
              \forall k \in \sets K. 
         \end{array} 
    \end{array}
\label{eq:scenario}
\tag{$\mathcal {IP}_{\bm w}^\prime(\{\widehat{\Xi}^k\}_{k\in \sets K})$}
\end{equation}

Any problem \eqref{eq:scenario} defined on partition~$\{\widehat{\Xi}^k\}_{k \in \sets K}$, where~$\bigcup_{k \in \sets K} \widehat{\Xi}^k = \widehat{\Xi}$ serves as a conservative approximation to~$\mathcal {IP}_{\bm w}({\widehat{\Xi}})$. Among these fixed partition problems, the one that achieves an optimal partition of~$\widehat{\Xi}$ will have the lowest objective value, which coincides with the value of~$\mathcal{IP}_{\bm w}({\widehat{\Xi}})$.

Once we have obtained an optimal solution~$(\overline{\theta}, \; \overline{\bm \psi}, \; \overline{\bm x},\; \{\overline{\bm y}^k\}_{k\in\mathcal{K}})$ to~$\mathcal {IP}_{\bm w}({\widehat{\Xi}})$, we solve the following separation problem to check whether the current solution is feasible and optimal across the entire uncertainty set~$\Xi$
\begin{equation}
    \begin{array}{ccl}
        z & = \displaystyle \max_{\overline{\bm \xi} \in \Xi} \;\;\min_{ { k} \in \sets K } \;\displaystyle\max_{ {\bm \xi} \in \Xi({\bm w^\prime},\overline{\bm \xi}) } \; & \max \left\{ \bm c^\top \overline{\bm \psi} + {\bm \xi}^\top {\bm C} \; \overline {\bm x} + {\bm \xi}^\top {\bm D} \; {\bm w^\prime}  + {\bm \xi}^\top {\bm Q} \; \overline{{\bm y}}^k - \overline{\bm \psi}^\top \bm g(\bm\xi) - \overline{\theta} \; ,\right. \\
        && \left. \max\limits_{l\in\left\{1,\ldots,L\right\}} \left\{ [\bm V(\bm \xi)]_l \overline{\bm x} + [\bm T(\bm \xi)]_l{\bm w^\prime} + [\bm W(\bm \xi)]_l{\overline{\bm y}^k} -{[\bm H]_l}{\bm \xi}\right\}
          \right\},
    \end{array}
\tag{$\mathcal {IP}_{\bm w}^{\text{out}}$}
\label{eq:separation}
\end{equation}
where~$l \in \{1,\ldots,L\}$ indices the uncertain constraints. If~$z > 0$, this indicates that either the solution is infeasible due to the maximum constraint violation taking a positive value or that there is an optimal~$\overline{\bm \xi}$ causing an objective value larger than~$\theta^\prime$. Theorem \ref{thm:sep-mip} shows that problem~\eqref{eq:separation} is amenable to an MIO reformulation. 

To obtain an optimal solution of~\eqref{eq:disjuctive}, we extend the branch-and-bound algorithm proposed in~\citet[Section 3.1]{Subramanyam2020} to a two-and-a-half\newqj{-}stage robust optimization problem. The algorithm starts by solving problem~\eqref{eq:scenario} with~$\widehat{\Xi} = \left\{\overline{\bm \xi}_0 \right\}$, where~$\overline{\bm \xi}_0$ is an arbitrary scenario from~$\Xi$. After getting an optimal solution, it solves the separation problem~\eqref{eq:separation} to check for any violations. Note that instead of optimizing over~$\{\widehat{\Xi}^k\}_{k \in \sets K}$, problem~\eqref{eq:scenario} fixes the partition of~$\widehat{\Xi}$. So to achieve an optimal solution to~\eqref{eq:disjuctive}, it separately adds the violated scenario into~$\widehat{\Xi}^1,\ldots,\widehat{\Xi}^k$ to create~$K$ new scenario-based problems. The procedure stops when there is no violation. We summarize the algorithm in Algorithm~\ref{alg:bnb}.

\begin{algorithm}[htbp!]
\setstretch{1.35}
\SetAlgoLined
\textbf{Inputs:} Initial scenario $\overline{\bm \xi}_0$ and number of policies $K$.

\textbf{Output:} Optimal solution~$(\bm \psi^\star ,\; \bm x^\star, \; \{{\bm {y}^k}^\star\}_{k\in\mathcal{K}})$, and optimal objective value~$\theta^*$ of problem \eqref{eq:Phi(w)}.

\textbf{Initialize:} Let~$\{\widehat{\Xi}^1,\ldots, \widehat{\Xi}^K\} \leftarrow \{ \{\overline{\bm \xi}_0\}, \{\}, \ldots, \, \{\} \}$, set~problem list~$\sets Q \leftarrow \{\mathcal P^\prime(\{\widehat{\Xi}^k\}_{k\in \sets K})\}$, and $\theta^\star \leftarrow +\infty$.

\While{$\sets Q \neq \emptyset$}
{
\begin{enumerate}[rightmargin=\dimexpr\linewidth-11.25cm-\leftmargin\relax,leftmargin=1.5em]
    \item Select, solve and remove a problem~$\mathcal P^\prime(\{\widehat{\Xi}^k\}_{k\in \sets K})$ from~$\sets Q$. Let~$(\overline{\theta}, \; \overline{\bm \psi}, \; \overline{\bm x},\; \{\overline{\bm y}^k\}_{k\in\mathcal{K}})$ be an optimal solution; if~$\theta^\prime > \theta^\star$, continue.
    \item Solve the separation problem~\eqref{eq:separation} let~$z$ be the optimal objective value and~$\overline {\bm \xi}$ be an optimal solution. If~$z > 0$, create~$K$ problems~$\mathcal P^\prime (\{\widehat{\Xi}^1,\ldots,\widehat{\Xi}^k \cup \{\overline{\bm \xi}\},\ldots,\widehat{\Xi}^K\}),$ for each~$k \in \sets K$, and add it to the problem list $\sets Q$; else, set~$(\theta^\star,\;\bm \psi^\star, \; \bm x^\star, \; \{{\bm {y}^k}^\star\}_{k\in\mathcal{K}}) \leftarrow (\overline{\theta}, \; \overline{\bm \psi}, \; \overline{\bm x},\; \{\overline{\bm y}^k\}_{k\in\mathcal{K}})$.
\end{enumerate}
}
\textbf{Result:} Report~$(\theta^\star, \; \bm \psi^\star, \bm x^\star, \; \{{\bm {y}^k}^\star\}_{k\in\mathcal{K}})$ as an optimal solution if~$\theta^\star < +\infty$, and problem is infeasible otherwise.
\caption{branch-and-bound algorithm for solving problem~\eqref{eq:Phi(w)}}
\label{alg:bnb}
\end{algorithm}

\subsubsection{A Branch-and-Cut Algorithm.}
\label{subsubsec:b&c}
In step 1 of Algorithm \ref{alg:bnb}, problem \eqref{eq:scenario} can be reformulated as a single\newqj{-}stage MILO by dualization. However, \newqj{T}heorem \ref{thm:det-scenario} shows that the size of the problem grows with~$|\widehat{\Xi}|$ with the number of dual variables being~$O(|\widehat{\Xi}|)$ and the number of constraints being~$O(|\widehat{\Xi}|)$, which is computationally challenging. Therefore, we propose to solve the robust optimization problem using a constraint generation algorithm\newqj{,} which will be embedded within the branch-and-cut algorithm that we now present. 

We start with the cut generation algorithm for solving problem~\eqref{eq:scenario}. Given collections of scenarios~$\{\widehat{\Xi}^k\}_{k\in \mathcal K}$, the following relaxation of~\eqref{eq:scenario} is solved
\begin{equation}
    \begin{array}{cl}
        \min & \quad \displaystyle \theta \\
         \st & \quad \theta\in\mathbb{R}, \; \bm\psi\in\mathbb{R}_+^{N_g}, \; {\bm x} \in \sets X, \; \bm y^k \in\mathcal{Y} \quad \forall k\in\sets K.\\
          & \quad \!\! \left. \begin{array}{l} 
         \theta \geq \bm c^\top \bm \psi +  \displaystyle \max_{\bm \xi \in \widehat\Xi^k_\text{in}} \; 
         \begin{array}[t]{r} {\bm \xi}^\top {\bm C} \; {\bm x} + {\bm \xi}^\top {\bm D} \; {\bm w^ \prime}  + {\bm \xi}^\top {\bm Q} \; {\bm y^k} \\
         - \bm\psi^\top \bm g(\bm\xi) 
         \end{array}\\
         {\bm T} ({\bm \xi}){\bm x} + {\bm V} ({\bm \xi}){\bm w ^\prime} + {\bm W} ({\bm \xi}){\bm y^k} \leq {\bm H}{\bm \xi} \quad \forall \bm \xi \in \widehat\Xi^k_\text{in}
         \end{array} \; \right\} \; \forall k \in \sets K,
    \end{array}
\label{eq:scenario-relaxed}
\tag{$\mathcal {IP}^{\prime \prime}_{\bm w}(\{\widehat{\Xi}_{\text{in}}^k\}_{k\in \sets K})$}
\end{equation}
where~$\widehat\Xi^k_\text{in}$ is a finite subset of $\bigcup_{\overline{\bm \xi} \in \widehat{\Xi}^k} \Xi({\bm w^\prime},\overline{\bm \xi})$. After getting a solution~$(\overline{\theta}, \; \overline{\bm \psi}, \; \overline{\bm x},\; \{\overline{\bm y}^k\}_{k\in\mathcal{K}})$ of~\eqref{eq:scenario-relaxed}, the following subproblems are solved to evaluate the violation of the constraints in~\eqref{eq:scenario}
\begin{equation}
    \begin{array}{cl}
        v^k = \displaystyle \max_{\overline{\bm \xi} \in \widehat{\Xi}^k} \;\displaystyle\max_{ {\bm \xi} \in \Xi({\bm w^\prime},\overline{\bm \xi}) } \; & \max \left\{ \bm c^\top \overline{\bm \psi} + {\bm \xi}^\top {\bm C} \; \overline{\bm x} + {\bm \xi}^\top {\bm D} \; {\bm w^\prime}  + {\bm \xi}^\top {\bm Q} \; \overline{{\bm y}}^k - \overline{\bm {\psi}}^\top \bm g(\bm\xi) - \overline{\theta} ,\right. \\
        & \left. \max\limits_{l\in\left\{1,\ldots,L\right\}} \left\{ \bm \xi^\top {\bm T_l}\overline{\bm x} + \bm \xi^\top{\bm V_l}{\bm w^\prime} + \bm \xi^\top{\bm W}_l{\overline{\bm y}^k} -{[\bm H]_l}{\bm \xi}\right\}
          \right\}\;\forall k \in \sets K.
    \end{array}
\tag{$\mathcal {IP}_{\bm w}^\text{in}$}
\label{eq:scenario-violated}
\end{equation}
We have~$v^k \leq 0$\ for all $k \in \sets K$, if and only if a robust optimal solution of~\eqref{eq:scenario} has been found.

The branch-and-cut algorithm for evaluating~$\Phi(\bm w^\prime)$ starts by solving problem~\eqref{eq:scenario-relaxed}. At incumbent nodes, it checks the inner robust feasibility by solving problem~\eqref{eq:scenario-violated} and adds a cut until a feasible solution to~\eqref{eq:scenario} is found. The feasible solution is then used to create branches or update the bounds as stated in  Algorithm~\ref{alg:bnb}. Algorithm~\ref{alg:bnc} summarizes the procedure, which we also illustrate in Figure 2.

% \linespread{0.75}
\begin{algorithm}[htbp!]
\setstretch{1.5}
\SetAlgoLined
\textbf{Inputs:} Initial scenario $\overline{\bm \xi}_0$, vector~$\bm w^\prime$, and  number of policies $K$.

\textbf{Output:} Optimal solution~$(\bm \psi^\star,\;\bm x^\star, \; \{{\bm {y}^k}^\star\}_{k\in\mathcal{K}})$ and optimal objective value~$\theta^\star$ of problem \eqref{eq:Phi(w)}.

\textbf{Initialize:} Set $\{\widehat{\Xi}^1,\ldots,\widehat{\Xi}^K\} \leftarrow \{ \{\overline{\bm \xi}_0\}, \{\}, \ldots, \, \{\}\}$,~$\widehat \Xi_\text{in}^k \leftarrow \widehat{\Xi}^k,\;\forall k\in \sets K$, $\theta^\star\leftarrow +\infty$. Start the branch-and-bound tree from problem~$\mathcal P^{\prime \prime}(\{\widehat{\Xi}_{\text in}^k\}_{k\in \sets K})$.

\textbf{1.} When an incumbent node is reached, get solution~$(\overline{\theta}, \; \overline{\bm \psi}, \; \overline{\bm x},\; \{\overline{\bm y}^k\}_{k\in\mathcal{K}})$, and solve problem~\eqref{eq:scenario-violated}. If~$\max_{k\in \sets K} v^k > 0$, go to step 2; else, go to step 3.

\textbf{2.} For all~$k \in \sets K$, let~$\bm \xi^k$ be an optimal solution to the~$k$th problem; if~$v^k > 0$, set~$\widehat{\Xi}_\text{in}^{k} \leftarrow \widehat{\Xi}_\text{in}^{k} \bigcup \{\bm \xi^k\}.$ Go to step 5.

\textbf{3.} Solve problem~\eqref{eq:separation}. If~$z > 0$, let~$\overline{\bm \xi}$ be an optimal solution, go to step 4; else, if~$\overline \theta \leq \theta^\star$, set~$\theta^\star \leftarrow \overline \theta$. Go to step 1.

\textbf{4.} Create~$K$ branches, each for problem~$\mathcal P^\prime (\{\widehat{\Xi}^1,\ldots,\widehat{\Xi}^k \cup \{\overline{\bm \xi}\},\ldots,\widehat{\Xi}^K\}),$~$k \in \sets K$. For each branch, start with solving the corresponding relaxation problem~$\mathcal P^{\prime \prime} (\{\widehat{\Xi}_{\text in}^1,\ldots,\widehat{\Xi}_{\text in}^k \cup \{\overline{\bm \xi}\},\ldots,\widehat{\Xi}_{\text in}^K\})$. Go to step 5. 

\textbf{5.} If the global lower bound~$lb_{\text{global}} \geq \theta^\star$ or no feasible node exist, terminate with status \texttt{Opt} or \texttt{Infeas}, and report the current solution; else if the node lower bound~$lb_{\text{node}} > \theta^\star$, fathom the node. Go to step 1.

\textbf{Result:} The branch-and-bound tree will end with either giving an optimal solution~$(\theta^\star,\;\bm \psi^\star,\;\bm x^\star, \; \{{\bm {y}^k}^\star\}_{k\in\mathcal{K}})$ or identifying that the problem as infeasible.

\caption{branch-and-cut algorithm for solving problem~\eqref{eq:Phi(w)}.}

\label{alg:bnc}
\end{algorithm}

\begin{figure}[htbp!]
\centering
\tikzstyle{startstop} = [rectangle, rounded corners, minimum width=3cm, minimum height=1cm,text centered, draw=black]
\tikzstyle{io} = [trapezium, trapezium left angle=70, trapezium right angle=110, minimum width=3cm, minimum height=1cm, text centered, draw=black]
\tikzstyle{process} = [rectangle, minimum width=3cm, minimum height=1cm, text centered, draw=black, execute at begin node=\setlength{\baselineskip}{3.5ex}]
\tikzstyle{decision} = [diamond, minimum width=3cm, minimum height=1cm, text centered, draw=black]

\tikzstyle{arrow} = [thick,->,>=stealth]

\scalebox{0.88}{
\begin{tikzpicture}[align=center,scale=0.75]
\small

\node (ini-in) [startstop] 
{Initialize and\\ start to solve~\eqref{eq:scenario-relaxed}};

\node (inc-node) [process, below of=ini-in, node distance=1.75cm] 
{Whenever at incumbent node,  \\get solution~$(\overline{\theta}, \; \overline{\bm \psi}, \; \overline{\bm x},\; \{\overline{\bm y}^k\}_{k\in\mathcal{K}})$};

\draw [arrow] (ini-in) -- (inc-node);

\node (in-vio) [process, below of=inc-node, node distance=1.75cm] 
{Check inner robust feasibility\\ by solving~\eqref{eq:scenario-violated} \\ for each~$k \in \sets K$};
\draw [arrow] (inc-node) -- (in-vio);

\node (in-feas) [decision, below of=in-vio, node distance=2cm, aspect=2, align=center] {$\max_{k\in\sets K} v^k > 0 ?$};
\draw [arrow] (in-vio) -- (in-feas);

\node (add-in) [process, right of=in-feas, node distance=4.2cm] 
{$\forall k\in \sets K$, if~$v^k > 0$ \\add~$\bm \xi^k$ to~$\widehat \Xi^k_\text{in}$};
\draw [arrow] (in-feas) -- node[anchor=south] {Yes} (add-in);

\node (out-vio) [process, below of=in-feas, node distance=2cm] 
{Check outer robust \\ feasibility by solving \\\eqref{eq:separation}};
\draw [arrow] (in-feas) -- node[anchor=west] {No} (out-vio);

\node (term) [decision, right of=out-vio, node distance=4.2cm, aspect=2, align=center]
{Terminate\\by Step 5?};
\draw [arrow] (add-in) -- (term);

\node (out-feas) [decision, below of=out-vio, node distance=2cm, aspect=2, align=center] {$z > 0 ?$};
\draw [arrow] (out-vio) -- (out-feas);

\node (add-out) [process, right of=out-feas, node distance=4.2cm] 
{Create~$K$ branches};
\draw [arrow] (out-feas) -- node[anchor=south] {Yes} (add-out);
\draw [arrow] (add-out) -- (term);

\node (output) [startstop, right of=add-out, node distance=3.7cm] {End with status\\ and solution};
\coordinate[right of = term, node distance = 3.7cm] (aux0);
\draw[arrow] (term) -- (aux0) -- node[anchor=west] {Yes} (output);

\node (fathom) [process, above of=output, node distance=5.5cm] 
{Check and fathom \\using Step 5};
\draw[arrow] (aux0) -- node[anchor=west] {No} (fathom);
\draw[arrow] (fathom) |- (inc-node);

\node (update) [process, left of=in-feas, node distance=3.7cm] 
{Update solution\\if~$\overline \theta < \theta^\star$ };
\coordinate[left of = out-feas, node distance = 3.7cm] (aux1);
\draw [arrow] (out-feas) -- node[anchor=south] {No} (aux1)  -- (update);

\draw [arrow] (update) |- (inc-node);
% \draw [arrow] (add-out) -- (term);

\end{tikzpicture}
}
\caption{Diagram of Algorithm \ref{alg:bnc} for solving problem~\eqref{eq:Phi(w)}}
\label{fig:alg-flowchart}
\end{figure}

\subsubsection{Convergence Analysis} We discuss the correctness as well as the asymptotic convergence of Algorithm~\ref{alg:bnc} in the following theorem.
\begin{theorem}
If the branch-and-cut Algorithm~\ref{alg:bnc} terminates, then it either returns an optimal solution to problem~\eqref{eq:Phi(w)} or correctly identifies it as infeasible. If there is no finite convergence, then in an infinite branch, every accumulation point of the sequence of the solutions obtained in step 1 of Algorithm~\ref{alg:bnc} corresponds an optimal solution~$(\bm \psi^\star,\;\bm x^\star,\;\{{\bm {y}^k}^\star\}_{k\in\mathcal{K}})$ of problem~\eqref{eq:Phi(w)} with objective value~$\theta^\star$.
\label{thm:bnc-correct}
\end{theorem}

\subsection{\texorpdfstring{Improved $L$-shaped Algorithm}{}}
\label{subsec:tighter-cuts}
In this section, we introduce three types of cuts to be added to Algorithm~\ref{alg:l-shaped}, which have the potential to improve its performance (see Table~\ref{table:r&d-per}).

\subsubsection{Benders' Cut}
\label{subsubsec:benders'}

A popular family of cuts in the integer programming literature is Benders' cuts, see~\citet{Benders1962PartitioningProblems}. In each iteration of Algorithm \eqref{alg:l-shaped}, for a fixed~$\bm w^r$, we relax the integrality of the decisions and solve an LP relaxation of problem~\eqref{eq:Phi(w)} and get the optimal dual solution to calculate the Benders' cut. Specifically, we consider the following LP problem
\begin{equation}
    \begin{array}{cl}
         \min & \quad \displaystyle \theta \\
         \st & \quad \theta\in\mathbb{R}, \; \bm\psi\in\mathbb{R}_+^{N_g}, \; {\bm x} \in \sets X^{\text {LP}}, \; \bm y(\bm \xi) \in \sets Y^{\text {LP}}\\
         & \quad \bm w = \bm w^r\\
         & \quad \!\! \left. \begin{array}{l} 
         \theta \geq \bm c^\top \bm \psi +  {\bm \xi}^\top {\bm C} \; {\bm x} + {\bm \xi}^\top {\bm D} \; {\bm w}  + {\bm \xi}^\top {\bm Q} \; {\bm y(\bm \xi)} - \bm\psi^\top \bm g(\bm\xi)  \\
         {\bm T} ({\bm \xi}){\bm x} + {\bm V} ({\bm \xi}){\bm w} + {\bm W} ({\bm \xi}){\bm y(\bm \xi)} \leq {\bm H}{\bm \xi} 
         \end{array} \quad \right\} \quad \forall \bm \xi \in \widehat \Xi,
    \end{array}
\label{eq:sub-relaxed}
\end{equation}
where~$\sets X^{\text {LP}}$ and $ \sets Y^{\text {LP}}$ are linear relaxations of~$\sets X$ and $ \sets Y$ by allowing all variables to be continuous, and~$\widehat \Xi$ can be any finite subset of~$\Xi$.

Problem~\eqref{eq:sub-relaxed} is a relaxation of~\eqref{eq:Phi(w)} in the sense that (i) it solves a fully adaptable problem on a finite set~$\widehat \Xi$, (ii) it ignores the fact that~$\bm y(\bm \xi)$ needs to be robust to the unobserved portion of~$\bm \xi$, and (iii) the mixed-integer sets~$\sets X$ and~$\sets Y$ are replaced by their linear relaxation. The Benders' cut \eqref{eq:benders-cut} added in step 3 of Algorithm~\ref{alg:l-shaped} is generated by obtaining  the optimal objective value $\overline \theta$ of problem~\eqref{eq:sub-relaxed}, and the values~$\bm \pi$ of the dual multipliers associated with the constraints $\bm w = \bm w^r$
\begin{equation}
\theta \geq \overline \theta + \bm \pi^\top (\bm w - \bm w^r).
\tag{$\mathcal C_\text{Benders}$}
\label{eq:benders-cut}
\end{equation}
In the implementation, we can either randomly generate or collect the scenario generated in step 3 of Algorithm~\ref{alg:bnc} to construct set~$\widehat{\Xi}$.

\subsubsection{Strengthened Feasibility Cut}
\label{subsubsec:strenthen-feas-cut}
If problem \eqref{eq:endo_whole} contains a subset of constraints in the form
\begin{equation}
\bm V(\bm \xi) \bm w \leq \bm H \bm \xi, \;\; \forall \bm \xi \in \Xi,
\label{eq: strengthen-feas-cstr}
\end{equation}
and the matrix $\bm V(\bm \xi)$ is component-wise non-negative for all $\bm \xi \in \Xi$, then for any infeasible solution $\bm w^r$, a solution $\bm w$ with $\bm w \geq \bm w^r$ is also infeasible for problem \eqref{eq:mins-w}. In this scenario, we incorporate the following strengthened feasibility cut in step 2 of Algorithm \ref{alg:l-shaped}
\begin{equation}
\sum_{i \in \sets I_r} 1 - w_i \geq 1.
\tag{$\mathcal C_\text{s-feas}$}
\label{eq:strengthen-feas}
\end{equation}
The cut \eqref{eq:strengthen-feas} is stronger than \eqref{eq:feasible-cut} in the sense that the former excludes all solutions~$\bm w \geq \bm w^r$ from the feasible region, while the latter excludes only~$\bm w^r$ itself.
\begin{theorem}
If problem \eqref{eq:endo_whole} contains a subset of constraints in the form of \eqref{eq: strengthen-feas-cstr}, and the matrix $\bm V(\bm \xi)$ is component-wise non-negative for all $\bm \xi \in \Xi$, then for any solution~$\bm w^r$ that is infeasible to problem \eqref{eq:endo_whole}, any solutions that violate constraint \eqref{eq:strengthen-feas} are also infeasible.
\label{thm:strengthen-feas}
\end{theorem}

\subsubsection{Strengthened Optimality Cut}
When the variable~$\bm w$ only appears in the objective function, and the cost vector of it is deterministic, problem~\eqref{eq:Phi(w)} can be written as
\begin{equation}
    \begin{array}{cl}
         \displaystyle\min & \; \bm d^\top \bm w + \Phi(\bm w) \\
         \st & \; {\bm w} \in \sets W,
    \end{array}\label{eq:mins-w-det-w}
\end{equation}
where
\begin{equation}
    \begin{array}{ccl}
         \Phi(\bm w) \;=\; &\displaystyle\min & \; \displaystyle \bm c^\top \bm \psi +\max_{\overline{\bm \xi} \in \Xi} \displaystyle \left\{ 
         \begin{array}{cl}
            \displaystyle\min_{ k \in \sets K } 
             &\displaystyle\max_{ {\bm \xi} \in \Xi({\bm w},\overline{\bm \xi}) } \;  {\bm \xi}^\top {\bm C} \; {\bm x} + {\bm \xi}^\top {\bm Q} \; {\bm y^k} - \bm\psi^\top \bm g(\bm\xi) \\
             \st & \;\; {\bm T} ({\bm \xi}){\bm x} + {\bm W} ({\bm \xi}){\bm y^k} \leq {\bm H}{\bm \xi} \;\;\; \forall {\bm \xi} \in \Xi({\bm w},\overline{\bm \xi})
         \end{array}
          \right\}  \\
         &\st & \; \bm\psi\in\mathbb{R}_+^{N_g},\;{\bm x} \in \sets X,\; \bm y^k \in\mathcal{Y},\; \forall k\in\sets K.
    \end{array}\label{eq:Phi(w)-det-w}
\end{equation}

For the function $\Phi(\bm w)$ defined in \eqref{eq:Phi(w)-det-w}, consider two vectors $\bm w$ and $\bm w^r$. If $\bm w$ is component-wise smaller than $\bm w^r$, i.e., $\bm w \leq \bm w^r$, then $\Phi(\bm w) \geq \Phi(\bm w^r)$. The intuition behind this is that if more entries of $\bm w$ take the value 1, the decision-maker has more information when making the recourse decision. This leads to the function $\Phi(\bm w)$ having a non-decreasing value, as discussed in \citet[Proposition 2]{paradiso2022exact}.
When applicable, we use the following strengthened optimality cut in step 3 of Algorithm~\ref{alg:l-shaped}
\begin{equation}
\theta \geq \Phi(\bm w^\prime) - \left(\Phi(\bm w^\prime)-L\right)\sum_{i \notin \sets I_r} w_{i}.
\tag{$\mathcal C_\text{s-opt}$}
\label{eq:information-cut}
\end{equation}

\subsubsection{Upper Cutoff}
\label{subsubsec:upper-cutoff}

To avoid spending time evaluating the function value $\Phi(\bm w^\prime)$ with suboptimal $\bm w^\prime$, in the initialization step of Algorithm~\ref{alg:bnc}, we set the cutoff value as the current upper bound $\theta^\star$ found by Algorithm~\ref{alg:l-shaped}. The termination condition for Algorithm~\ref{alg:bnc} is met when the current lower bound exceeds the cutoff value. In such cases, $\bm w^\prime$ is proven to be suboptimal. In the case of early termination, the integer cut \eqref{eq:integer-cut} is generated using the cutoff value.

The alternative to step 2 of Algorithm~\ref{alg:l-shaped} is as follows:

\textbf{$2^\prime$}: Evaluate the value~$\Phi(\bm w^\prime)$ with cutoff value~$\theta^\star$.

The additional initialization and add-on to step 5 of Algorithm~\ref{alg:bnc} are given as follows:

\textbf{Initialize}: Set cutoff value~$\theta_{\text {ub}} \leftarrow \theta^\star$. 

\textbf{$5^+$}: If the global lower bound~$\text {lb}_{\text {global}} > \theta_{\text {ub}}$, terminate and return the current incumbent solution.

%%%%%%%%%%%%%%%%%%%%%%%%%%%%%%%%%%%%%%%%%%%%%%%%%%%%%%%%%%%%%%%%%%%%%%%%%%%%%%%%%%%%%%%%%%%%%%%%
%%%%%%%%%%%%%%%%%%%%%%%%%%%%%%%%%%%%%%%%%%%%%%%%%%%%%%%%%%%%%%%%%%%%%%%%%%%%%%%%%%%%%%%%%%%%%%%%
%%%%%%%%%%%%%%%%%%%%%%%%%%%%%%%% NUMERICAL EXPERIMENT %%%%%%%%%%%%%%%%%%%%%%%%%%%%%%%%%%%
%%%%%%%%%%%%%%%%%%%%%%%%%%%%%%%%%%%%%%%%%%%%%%%%%%%%%%%%%%%%%%%%%%%%%%%%%%%%%%%%%%%%%%%%%%%%%%%%
%%%%%%%%%%%%%%%%%%%%%%%%%%%%%%%%%%%%%%%%%%%%%%%%%%%%%%%%%%%%%%%%%%%%%%%%%%%%%%%%%%%%%%%%%%%%%%%%

\section{Numerical Results}
\label{sec:numerical}
In this section, we investigate the performance of our approach on the best box and R\&D project portfolio optimization problems. We show the trade-off between optimality and computational effort as the number of policies~$K$ is varied and the suboptimality of the robust optimization approach proposed in \citet{Vayanos2020} under worst-case distribution, thereby demonstrating the value of our method relative to existing literature. In both problems, we report the results of solving the MINLO problem by Gurobi 9.0.1. In the R\&D project portfolio optimization problem, we showcase the attractive scalability properties of our proposed decomposition algorithm. To solve the robust optimization counterpart, where the ambiguity set \eqref{eq:ambi_set} only includes the support, we employ the \ROCPP\, package~\citet{Vayanos2020ROC++:C++} for coding and solving the problem. Furthermore, in Appendix \ref{ecsec:r&d-loan}, we present the numerical results of an extension to the R\&D problem. The results show that the proposed solution method works equally effectively as in the R\&D problem.

All of our experiments are performed on the High-Performance Computing Cluster of our university. Each job is allocated 2GB \newqj{of} RAM, 4 cores, and a 2.1 GHz Xeon processor. For each instance, an optimality gap of~$0.1\%$ and a time limit of 7,200 seconds are allowed to solve the problem. 
% All implementation codes are available in \citet{droddid}.

\subsection{Distributionally Robust Best Box problem}
\subsubsection{Problem Description}
We perform experiments on several instances of the two-stage best box problem that subsumes selection problems as special cases. It is a distributionally robust variant of the robust problem in~\citet{Vayanos2020, Vayanos2020ROC++:C++}.
In this problem, we have~$N$ boxes indexed in the set~$\sets N = \{ 1, \ldots, N\}$ to choose from. Each box~$n$ contains an unknown prize with value~$\xi^\mathrm{r}_n$, and opening it incurs an unknown cost~$\xi^\mathrm{c}_n$. The return and cost are uncertain and will only be revealed if the box is opened. In the first stage, the decision-maker decides whether to open each box~$n \in \sets N$, which we indicate with the decision variable~$w^r_n \in \{0, \, 1\}$. Thus,~$\xi_n^r$ and~$\xi_n^c$ are revealed between the first and second time stage if only~$w^r_n = 1$. The total budget available to open boxes is~$B$. In the second stage, the decision-maker keeps one of the opened boxes and earns its value, indicated with the decision variable~$y_n \in \{0, \, 1\}$. The goal of the decision-maker is to select the boxes to open and keep to maximize the worst-case expected value of the box kept. We assume that the values and costs are expressible as~$ \xi^\mathrm{r}_n=(1+\bm \Phi_n^\top\bm \zeta /2)r_n$ and ~$\xi^\mathrm{c}_n=(1+\bm \Psi_n^\top\bm \zeta/2) c_n$, where~$r_n$ and $c_n$ correspond to the nominal value and cost for box~$n$, respectively,~$\bm \zeta \in \left[-1,\; 1\right]^M$ are~$M$ risk factors, and the vectors~$\bm{\Phi}_n,\; \bm{\Psi}_n \in \mathbb{R}^M$ collect the loading factors for  box~$n$. The support of the distribution of the uncertain parameters is described as follows
\begin{equation*}
\begin{array}[t]{r}
    \Xi=\left\{\left(\bm \xi^{\mathrm{r}}, \bm{\xi}^{\mathrm{c}}\right) \in \mathbb{R}^{2 N}: \exists \bm{\zeta} \in[-1,1]^{M}: \xi_{n}^{\mathrm{r}}=\left(1+\mathbf{\Phi}_{n}^{\top} \bm{\zeta} / 2\right) r_{n}, \; \xi_{n}^{\mathrm{c}}=\left(1+\bm{\Psi}_{n}^{\top} \bm{\zeta} / 2\right) c_{n}, \right .\\ \left. n \in \sets N\right\}.
\end{array}
\end{equation*}
% We use the following ambiguity set for the uncertain values and costs:
% \begin{equation}
% \begin{array}[t]{cl}
%      \mathcal{P}_\text{BB} = \Big\{\mathbb{P}\in\mathcal{M}_+(\mathbb{R}^{2N})\; : & \; \mathbb{P}((\bm \xi^r, \bm \xi^c)\in\Xi) = 1, \\
%      & \; \mathbb{E}_{\mathbb{P}}\left[\left\{|\mathbf e^\top \left((\bm \xi^r, \bm \xi^c) - (\bm r, \bm c) \right )\right| \right ]\leq 0.25 N^{-1/2} \mathbf e^\top (\bm r, \bm c))\Big\}.
% \end{array}
% \label{eq: bb-ambi_set}
% \end{equation}
We define $\bm \xi = \left(\bm \xi^\mathrm{r},\bm\xi^\mathrm{c}\right)$,~$\bm \mu = \left(\bm r,\bm c\right)$, and use the following ambiguity set for the uncertain values and costs
\begin{equation}
     \mathcal{P}_\text{BB} = \Big\{\mathbb{P}\in\mathcal{M}_+(\mathbb{R}^{N_\xi})\; :\; \mathbb{P}(\bm \xi\in\Xi) = 1,\;\; \mathbb{E}_{\mathbb{P}}\left[\left\{|\mathbf e^\top \left(\bm \xi - \bm \mu \right )\right| \right ]\leq 0.25 N^{-1/2} \mathbf e^\top \bm \mu)\Big\}.
\label{eq: bb-ambi_set}
\end{equation}
The last constraint in the ambiguity set follows the intuition of the central limit theorem, and it imposes an upper bound on the cumulative deviation of the boxes' returns and costs from their expected values.

The distributionally robust best box problem can be written as
\begin{equation}
    \begin{array}{cl}
        \min & \quad \displaystyle \sup_{\mathbb{P} \in \mathcal{P}} \;\;\;\mathbb{E}_{\mathbb{P}}\left[  -{\bm \xi^\mathrm{r}}^\top  {\bm y}({\bm \xi}) \right] \\
         \st & \quad {\bm w} = (\bm w^r,\, \bm w^r),\, \bm w^r \in \left\{0,\;1\right\}^N \\
         %& \quad {\bm R}{\bm x} + {\bm S}{\bm w} \; \leq \; {\bm t} \\
         & \quad \!\! \left. \begin{array}{l} 
         {\bm y}({\bm \xi}) \in \left\{0,\;1\right\}^N  \\
         {\bm \xi^\mathrm{c}}^\top \bm w^r  \leq B \\
         \bm y(\bm \xi) \leq \bm w^r \\
         \mathbf e^\top \bm y(\bm \xi) = 1
         \end{array} \quad \right\} \quad \forall {\bm \xi} \in \Xi, \\
         & \quad {\bm y}({\bm \xi}) = {\bm y}({\bm \xi}') \quad \forall {\bm \xi}, \; {\bm \xi}' \in \Xi \; : \; {\bm w} \circ {\bm \xi} = {\bm w} \circ {\bm \xi}'.
    \end{array}
\label{eq:bb}
\end{equation}
Besides the support, the ambiguity set~\eqref{eq: bb-ambi_set} has only one constraint, leading to a scalar dual variable in the MINLO reformulation of~\eqref{eq:bb}. Therefore, the reformulation has only one bilinear term with continuous variables, and the reformulated MINLO problem can be handled by Gurobi efficiently.

\subsubsection{Computational Results}
\label{subsubsec: b&b-res}

We evaluate the performance of our approach on 80 randomly generated instances of problem~\eqref{eq:bb}:~$20$ instances for each~$N \in \{10,20,30,40\}$ with risk factors~$M \in \{5, 10, 15, 20\}$. In each instance,~$\bm c$ is drawn uniformly at random from the hypercube~$[0,10]^N$, and~$\bm r = \bm c/5,\;B = \mathbf{e}^\top \bm c/2$. The risk loading factors~$\bm \Phi_n$ and~$\bm \Psi_n$ are sampled uniformly at random from the hypercube~$[-1,1]^M$. We solve the~$K$-adaptability counterpart of~\eqref{eq:bb} for $K \in \{1,2,3,4\}$. We summarize our results in the following.

In the best box problem, as the number of bilinear terms between continuous variables is small, solving the monolithic MINLO is more efficient than using the decomposition algorithm in all cases. Table~\ref{table:bb-per} summarizes the computational results across the generated instances. The improvement over the static solution increases rapidly with the number of policies. The average improvements over the static solution for all sizes are~$87.1\%,\, 111.3\%,\,121.3\%$ for~$K = 2,\;3,\,4$, respectively. The improvement in solution quality comes at a computational cost. 

\newqj{To study the quality of the solutions obtained by the~$K$-adaptability approximation, we solve a small-sized problem with~$N = 10$ to optimality for~$K$ varying from~$1$ to~$11$. We report the average improvements over the static solution as~$K$ increases in Section~\ref{ecsubsub:bb-improve-k}. The results show that a small number of policies (3 out of 11) is sufficient for the~$K$-adaptability approximation \eqref{eq:k-nested} to be optimal in the original problem~\eqref{eq:nested}.}

\begin{table}[!ht]
\renewcommand{\arraystretch}{1.5}
% \center
\begin{tabular}{ ccccc } 
\hline
$K$ & $N=10$ & $N=20$ & $N=30$ & $N=40$ \\ 
 \hline
$ 1 $ &$ 0.0\% / 0.1s / 0.0\% $&$ 0.0\% / 0.2s / 0.0\% $&$ 0.0\% / 0.4s / 0.0\% $&$ 0.0\% / 0.5s / 0.0\% $\\
$ 2 $ &$42.6\% / 0.2s / 0.0\% $&$ 86.5\% / 6.0s / 0.0\% $&$ 118.2\% / 1090.5s / 0.0\% $&$ 101.1\% / 2530.6s / 0.0\% $\\
$ 3 $ &$43.9\% / 0.4s / 0.0\% $&$117.2\% / 8.8s / 0.0\% $&$148.9\% / 1168.9s / 0.0\% $&$ 135.2\% / 2297.8s / 0.0\% $\\
$ 4 $ &$ 43.9\% / 0.8s / 0.0\% $&$ 125.4\% / 16.8s / 0.0\% $&$ 162.2\% / 3031.6s / 0.0\% $&$ 153.5\% / 4819.0s / 0.0\% $\\
\hline
\end{tabular}
\caption{Summary of computational results for the best box problem. Each entry corresponds to: the average objective value improvement of the best~$K$-adaptability solution found in the time limit over the best static solution found in the time limit; the improvement is calculated as the ratio of the difference between the objective values of the $K$-adaptability solution and the static solution to the objective value of the static solution / the average solution time of the solved instances/ the average optimality gap across unsolved cases.}
\label{table:bb-per}

\end{table}

We now assess the suboptimality of an optimal solution to the RO problem in the DRO problem. To compute such suboptimality, we first solve the robust optimization counterpart of problem~\eqref{eq:bb}, where the ambiguity set~\eqref{eq: bb-ambi_set} includes only support information. Let $(\bm x_\text{RO},\; \bm w_\text{RO}, \;\{\bm y_\text{RO}^k\}_{k \in \sets K})$ denote an optimal solution to the RO problem, we fix these values in \eqref{eq:bb} and solve the inner maximization problem to obtain the objective value of the RO solution under the worst case distribution, denoted by $z^\star_\text{RO}$. Accordingly, we let $z^\star_\text{DRO}$ denote the optimal objective value obtained by optimizing problem \eqref{eq:bb}. Table~\ref{table:bb-subopt} summarizes the relative difference between $z^\star_\text{RO}$ and $z^\star_\text{DRO}.$ From the table, we see that the average suboptimality of the RO solution is as high as~$40.59\%$.

\begin{table}[!ht]
\renewcommand{\arraystretch}{1.5}
% \begin{center}
\begin{tabular}{ ccccc } 
 \hline
$K$ & $N=10$ & $N=20$ & $N=30$ & $N=40$ \\ 
 \hline
$ 1 $ & $ 40.59\% $ & $ 13.19\% $ & $ 3.94\% $ & $ 6.54\% $ \\
$ 2 $ & $ 25.40\% $ & $ 1.60\% $ & $ 1.18\% $ & $ 1.43\% $ \\
$ 3 $ & $ 20.09\% $ & $ 2.10\% $ & $ 1.52\% $ & $ 0.77\% $ \\
$ 4 $ & $ 19.65\% $ & $ 1.98\% $ & $ 1.05\% $ & $ 1.75\% $ \\
 \hline
\end{tabular}
% \end{center}
\caption{Suboptimality of the RO solutions to the best box problem. Each cell represents an average of 20 instances. In each instance, the relative differences are calculated by $\displaystyle \frac{(z^\star_\text{DRO} - z^\star_\text{RO})}{(z^\star_\text{DRO} + z^\star_\text{RO})/2}$.}
\label{table:bb-subopt}
\end{table}

To further assess the performance of the DRO and RO solutions, in Appendix \ref{ecsubsub:bb-out-sample}, we assess the out-of-sample performance of these solutions. The results show that the DRO solutions, in general, outperform the RO solutions also under non-worst-case distributions.

\subsection{Distributionally Robust R\&D Project Portfolio Optimization}
\label{subsec:r&d}
\subsubsection{Problem Description}
We next perform experiments on several instances of the distributionally robust R\&D project portfolio optimization problem -- this is a DRO variant of the robust version proposed in \citet{Vayanos2020, Vayanos2020ROC++:C++}. In this problem, there are~$N$ projects indexed in the set~$\sets N = \{ 1, \ldots, N\}$, and the decision-maker decides whether and when to invest in a project. Each project~$n$ has an uncertain return~$\xi_n^\mathrm{r}$, and undertaking it incurs an uncertain cost~$\xi_n^\mathrm{c}$. The return and cost of \newqj{the} project~$n$ can only be observed if we undertake this project. The decision to invest in project~$n$ in the first (resp. second) stage is denoted by~$w^r_n$ (resp.~$y_n$),~$n \in \sets N$. Thus,~$\xi_n^r$ and~$\xi_n^c$ are revealed between the first and second time stage if only~$w^r_n = 1$. The total budget available across two years is~$B$. A project can be invested in at most once, and only a fraction~$\theta \in [0, \,1)$ of the return is obtained for investments made in the second stage. The goal of the decision-maker is to maximize the worst-case expected return. We assume the returns and costs are expressible as~$ \xi^\mathrm{r}_n=(1+\bm \Phi_n^\top\bm \zeta /2) r_n$ and ~$\xi^\mathrm{c}_n=(1+\bm \Psi_n^\top\bm \zeta/2) c_n$, where~$r_n,\;c_n$ correspond to the nominal return and cost of project~$n$, respectively,~$\bm \zeta \in \left[-1,\; 1\right]^M$ are~$M$ risk factors, and the vectors~$\bm{\Phi}_n$ and~$\bm{\Psi}_n \in \mathbb{R}^M$ collect the loading factors for project~$n$. The support of the distribution of the uncertain parameters is described as follows
\begin{equation*}
\begin{array}[t]{r}
    \Xi=\left\{\left(\bm \xi^{\mathrm{r}}, \bm{\xi}^{\mathrm{c}}\right) \in \mathbb{R}^{2 N}: \exists \bm{\zeta} \in[-1,1]^{M}: {\xi}_{n}^{\mathrm{r}}=\left(1+\mathbf{\Phi}_{n}^{\top} \bm{\zeta} / 2\right){r}_{n}, \;{\xi}_{n}^{\mathrm{c}}=\left(1+\bm{\Psi}_{n}^{\top} \bm{\zeta} / 2\right) {c}_{n}, \right .\\ \left . n\in \sets N\right\}.
\end{array}
\end{equation*}
We define $\bm \xi = \left(\bm \xi^\mathrm{r},\bm\xi^\mathrm{c}\right)$,~$\bm \mu = \left(\bm r,\bm c\right)$, and use the following ambiguity set for the uncertain returns and costs
\begin{equation}
\begin{array}{cl}
    \mathcal{P}_\text{RD} = \Big\{\mathbb{P}\in\mathcal{M}_+(\mathbb{R}^{N_\xi})\; :\; \mathbb{P}(\bm \xi\in\Xi) = 1,\;\;&\mathbb{E}_{\mathbb{P}}\left[\left|\bm \xi - \bm \mu \right| \right ]\leq 0.15 \bm \mu,\\ &\mathbb{E}_{\mathbb{P}}\left[\left|\mathbf e^\top \left(\bm \xi - \bm \mu \right )\right| \right ]\leq 0.15 N^{-1/2} \mathbf e^\top \bm \mu)\Big\}.
\end{array}
\label{eq: r&d-ambi_set}
\end{equation} 
The second set of constraints in the ambiguity set enforces the mean absolute deviation of the individual returns and costs to be bounded. The last constraint imposes an upper bound on the cumulative deviation of the projects' returns and costs from their expected value.

The distributionally robust R\&D project portfolio optimization problem can be written as
\begin{equation}
    \begin{array}{cl}
         \min & \quad \displaystyle \sup_{\mathbb{P} \in \mathcal{P}} \;\;\;\mathbb{E}_{\mathbb{P}}\left[  -{\bm \xi^r}^\top \left(\bm w + \theta {\bm y}({\bm \xi}) \right)\right] \\
         \st & \quad {\bm w} = (\bm w^r,\; \bm w^r),\; \bm w^r \in \left\{0,\;1\right\}^N \\
         %& \quad {\bm R}{\bm x} + {\bm S}{\bm w} \; \leq \; {\bm t} \\
         & \quad \!\! \left. \begin{array}{l} 
         {\bm y}({\bm \xi}) \in \left\{0,\;1\right\}^N  \\
         {\bm \xi^c}^\top \left(\bm w^r + {\bm y}({\bm \xi}) \right) \leq B \\
         \bm w^r + \bm y(\bm \xi) \leq \mathbf e
         \end{array} \quad \right\} \quad \forall {\bm \xi} \in \Xi \\
         & \quad {\bm y}({\bm \xi}) = {\bm y}({\bm \xi}') \quad \forall {\bm \xi}, {\bm \xi}' \in \Xi \; : \; {\bm w} \circ {\bm \xi} = {\bm w} \circ {\bm \xi}'.
    \end{array}
\label{eq:r&d}
\end{equation}

\subsubsection{Computational Results}
\label{subsubsec: r&d-result}

We run experiments on~$80$ randomly generated instances of the~$K$-adaptability counterpart of problem~\eqref{eq:r&d}:~$20$ instances for each~$N\in\{5,10,15,20\}$~with risk factors~$M \in \{4,5,8,10\}$, respectively. We draw~$\bm c$ uniformly at random from the box~$[0,10]^N$ and let~$\bm r = \bm c/5,\;B = \mathbf e^\top \bm c/2$. The risk load factor~$\bm \Phi_n$ and~$\bm \Psi_n$ are sampled uniformly at random from the standard simplex. 

We report the results of solving problem \eqref{eq:r&d} by utilizing the monolithic MINLO reformulation and applying the decomposition Algorithm~\ref{alg:l-shaped} with and without the cuts introduced in Section~\ref{subsec:tighter-cuts}. In the decomposition Algorithm~\ref{alg:l-shaped}, we set the time limit of evaluating~$\Phi(\bm w)$ to~$[300s, 600s, 1200s, 2400s]$ for problems with size~$N\in\{5,\;10,\;15,\;20\}$. The lower bound value~$L$ in Algorithm \ref{alg:l-shaped} is found by solving a deterministic problem where the uncertain parameters are fixed to~$\bm \xi = (\bm r, \bm c)$. If the subproblem evaluating~$\Phi(\bm w)$ is not solved to optimality within the time limit, we use the current lower bound of the subproblem to generate a valid integer cut of the form~\eqref{eq:integer-cut} and use the upper bound to update the objective value~$\theta^\star$ in step 4. If the main problem is not solved to optimality within the time limit, we 
report~$\theta^\star$ as the objective value and the current optimal objective value of \newqj{the} main problem~\eqref{eq:mins-w-master} as the lower bound.

\begin{landscape}
\setlength{\tabcolsep}{2pt}
\begin{table}[!ht]
\renewcommand{\arraystretch}{1.5}
\begin{footnotesize}
\begin{tabular}{cccccccccccccccccccc } 
 \hline
 & & & \multicolumn{5}{c}{MINLO} & & \multicolumn{5}{c}{ Decomposition} & & \multicolumn{5}{c}{ Decomposition with Added Cuts}\\
 $N$ & $K$ & $\,$ & Opt($\#$) & Time(s) & Gap & Better($\#$) & Improvement & $\,$ &  Opt($\#$) & Time(s) & Gap & Better($\#$) & Improvement & $\,$ &  Opt($\#$) & Time(s) & Gap & Better($\#$) & Improvement\\ 
 \hline
 \multirow{3}{*}{5} & 2 && $ \textbf{20} $ & $ \textbf{1.6} $ & $ \textbf{0.0\%} $ & $ \textbf{0} $ & $ \textbf{17.5\%}$  && $ \textbf{20} $ & $ 54.3 $ & $ \textbf{0.0\%} $ & $ \textbf{0} $ & $ \textbf{17.5\%} $ &
& $ \textbf{20} $ & $ 16.0 $ & $ \textbf{0.0\%} $ & $ \textbf{0} $ & $ \textbf{17.5\%} $ \\
 ~ & 3 && $ 15 $ & $ 912.3 $ & $ 0.03 \% $ & $ \textbf{0} $ & $ \textbf{21.3\%} $ && $ \textbf{20} $ & $ 435.4 $ & $ \textbf{0.0\%} $ & $ \textbf{0} $ & $ \textbf{21.3\%} $ && $ \textbf{20} $ & $ \textbf{239.0} $ & $ \textbf{0.0\%} $ & $ \textbf{0} $ & $ \textbf{21.3\%} $ 
 \\
 ~ & 4 && $ 0 $ & $ 7200 $ & $ 0.2\% $ & $ \textbf{1} $ & $ \textbf{22.3\%} $ && $ \textbf{20} $ & $ 854.8 $ & $ \textbf{0.0\%} $ & $ 0 $ & $ \textbf{22.3\%} $ && $ \textbf{20} $ & $ \textbf{475.0} $ & $ \textbf{0.0\%} $ & $ \textbf{1} $ & $ \textbf{22.3\%} $  \\
\hline
 \multirow{3}{*}{10} & 2 && $ 19 $ & $ \textbf{730.5} $ &  -  & $ \textbf{0} $ & $ 16.9 \%$ &
& $ 0 $ & $ \textbf{7200} $ & $ 27.2\% $ & $ \textbf{0} $ & $ 20.3\%$ &
& $ \textbf{20} $ & $ 1826.8 $ & $ \textbf{0.0\%} $ & $ \textbf{0} $ & $ \textbf{23.1\%} $ \\
 ~ & 3 && $ \textbf{0} $ & $ \textbf{7200} $ & - & $ 2 $ & $ 11.8 \%$ &
& $ \textbf{0} $ & $ \textbf{7200} $ & $ 30.2\% $ & $ 1 $ & $ 20.7\%$ &
& $ \textbf{0} $ & $ \textbf{7200} $ & $ \textbf{16.7\%} $ & $ \textbf{16} $ & $ \textbf{26.7\%} $  \\
 ~ & 4 && $ \textbf{0} $ & $ \textbf{7200} $ &  - & $ 1 $ & $ 23.8 \%$ &
& $ \textbf{0} $ & $ \textbf{7200} $ & $ 29.4\% $ & $ 2 $ & $ 21.6\%$ &
& $ \textbf{0} $ & $ \textbf{7200} $ & $ \textbf{16.5\%} $ & $ \textbf{16} $ & $ \textbf{27.4\%} $  \\
 \hline
 \multirow{3}{*}{15} & 2 && $ \textbf{0} $ & $ \textbf{7200} $ &  -  & $ \textbf{18} $ & $ 16.1 \%$ &
& $ \textbf{0} $ & $ \textbf{7200} $ & $ 32.9\% $ & $ 0 $ & $ 8.4\%$ &
& $ \textbf{0} $ & $ \textbf{7200} $ & $ \textbf{23.9\%} $ & $ 2 $ & $ \textbf{20.1\%} $  \\
 ~ & 3 && $ \textbf{0} $ & $ \textbf{7200} $ & - & $ 7 $ & $ 21.5 \%$ &
& $ \textbf{0} $ & $ \textbf{7200} $ & $ 30.5\% $ & $ 0 $ & $ 12.4\%$ &
& $ \textbf{0} $ & $ \textbf{7200} $ & $ \textbf{22.2\%} $ & $ \textbf{13} $ & $ \textbf{23.0\%} $ \\
 ~ & 4 && $ \textbf{0} $ & $ \textbf{7200} $ &  -  & $ 7 $ & $ 20.9 \%$ &
& $ \textbf{0} $ & $ \textbf{7200} $ & $ 28.8\% $ & $ 0 $ & $ 15.1\%$ &
& $ \textbf{0} $ & $ \textbf{7200} $ & $ \textbf{21.9\%} $ & $ \textbf{13} $ & $ \textbf{23.2\%} $ \\
 \hline
 \multirow{3}{*}{20} & 2 && $ \textbf{0} $ & $ \textbf{7200} $ & - & $ \textbf{20} $ & $ \textbf{21.9\%} $ &
& $ \textbf{0} $ & $ \textbf{7200} $ & $ 36.9\% $ & $ 0 $ & $ 0.9\%$  &
& $ \textbf{0} $ & $ \textbf{7200} $ & $ \textbf{24.5\%} $ & $ 0 $ & $ 18.7 \%$ \\
 ~ & 3 && $ \textbf{0} $ & $ \textbf{7200} $ & - & $ 7 $ & $ 20.7 \%$ &
& $ \textbf{0} $ & $ \textbf{7200} $ & $ 34.8\% $ & $ 0 $ & $ 4.3\%$ &
& $ \textbf{0} $ & $ \textbf{7200} $ & $ \textbf{22.5\%} $ & $ \textbf{13} $ & $ \textbf{21.6\%} $  \\
 ~ & 4 && $ \textbf{0} $ & $ \textbf{7200} $ & - & $ 1 $ & $ 18.9 \%$ &
& $ \textbf{0} $ & $ \textbf{7200} $ & $ 33.9\% $ & $ 0 $ & $ 5.7\%$ &
& $ \textbf{0} $ & $ \textbf{7200} $ & $ \textbf{21.5\%} $ & $ \textbf{19} $ & $ \textbf{23.6\%} $  \\
 \hline
\end{tabular}
\caption{Summary of computational results of the R$\&$D problem. $\mathbf{Decomposition}$ and $\mathbf{ Decomposition\;\;with\;\;Added\;\;Cuts}$ 
refer to Algorithm~\ref{alg:l-shaped} without and with the cuts introduced in Section~\ref{subsec:tighter-cuts}, respectively. $\mathbf{Opt(\#)}$ corresponds to the number of instances solved to optimality, $\mathbf{Time(s)}$ to the average computational time (in seconds) for instances solved to optimality, and $\mathbf{Gap}$ to the average optimality gap for the instances not solved within the time limit. We write '-' when no valid lower bound was found. $\mathbf{Better(\#)}$ denotes the number of instances in each method that achieved a better objective value than the other method. $\mathbf{Improvement}$ denotes the average improvement in the objective value of the $K$-adaptability solution found in the time limit over the static solution found in the time limit. The improvement is calculated as the ratio of the difference between the objective values of the $K$-adaptability solution and the static solution to the objective value of the static solution.}
\label{table:r&d-per}
\end{footnotesize}
\end{table}
\end{landscape}

Table~\ref{table:r&d-per} summarizes the computational results across those instances. From the table, we observe that the improved decomposition algorithm achieves optimality or finds better solutions for more instances than the basic decomposition algorithm and solves the monolithic MINLO problem when $K$ is large. Across all instances, the proposed algorithm has a significantly smaller optimality gap,\newqj{,} which decreases with~$K$. We notice that solving the MINLO problem often results in the inability to obtain a satisfactory or valid lower bound, especially when $N \geq 10$. As the number of policies~$K$ increases, the number of instances for which the improved decomposition algorithm achieves a better objective value increases. For example, when~$N=20$, the number of better solutions found by our approach is~$\{0, \, 13, \, 19\}$ with~$K \in \{2, \,3, \,4\}$, respectively. The improvements in the objective value over the static solutions of the improved decomposition algorithm are $\{18.7\%, \, 21.6\%, \, 23.6\%\}$, whereas the improvements obtained by solving the MINLO problem and using the basic decomposition algorithm are $\{21.9\%, \, 20.7\%, \,18.9\%\}$ and $\{0.9\%, \, 4.3\%, \, 5.7\%\}$, respectively. Note that in some settings, MINLO solutions achieve a better objective value in more instances, but the average improvement in the objective value is worse. This discrepancy arises when the MINLO solutions are sometimes significantly suboptimal, thereby reducing the overall average improvement. Finally, in all cases, the improvement in the objective value over the static solutions of the decomposition algorithm increases with the number of policies~$K$. This aligns with the rationale of the~$K$-adaptability approach, where a larger~$K$ offers more flexibility for the recourse decisions, leading to a better objective value. However, the trend is not clear in MINLO solutions when~$N \geq 10$ due to computational difficulties. 

\newqj{To study the quality of the solutions obtained by the~$K$-adaptability approximation, in Section~\ref{ecsubsub:r&d-improve-k}, we solve a small-sized problem with~$N = 4$ and~$M = 2$ to optimality for all~$K$ varying from~$1$ to~$16$. We report the average improvements over the static solution as~$K$ increases. The results show that a small number of policies (5 out of 16) is sufficient for the~$K$-adaptability approximation \eqref{eq:k-nested} to be optimal in the original problem~\eqref{eq:nested}.}

Next, we examine the suboptimality of the RO solution. We define~$z^\star_\text{RO}$ and~$z^\star_\text{DRO}$ as in Section \ref{subsubsec: b&b-res}. The relative differences between~$z^\star_\text{RO}$ and~$z^\star_\text{DRO}$ are summarized in Table~\ref{table:r&d-subopt}.
\begin{table}[!ht]
\renewcommand{\arraystretch}{1.5}
\begin{tabular}{ ccccc } 
 \hline
$K$ & $N=5$ & $N=10$ & $N=15$ & $N=20$ \\ 
 \hline
$1$ & $ 0.37\% $ & $ 1.60\% $ & $ 3.24\% $ & $ 2.08\% $  \\
$2$ & $ 17.33\% $ & $ 23.03\% $ & $ 20.37\% $ & $ 19.61\% $ \\
$3$ & $ 14.62\% $ & $ 19.15\% $ & $ 17.65\% $ & $ 18.33\% $ \\
$4$ & $ 15.05\% $ & $ 11.25\% $ & $ 7.81\% $ & $ 8.73\% $ \\
\hline
\end{tabular}
\caption{Suboptimality of RO solution to the R$\&$D problem. Each cell represents an average of 20 instances. In each instance, the relative differences are calculated by $\displaystyle \frac{(z^*_\text{DRO} - z^*_\text{RO})}{(z^*_\text{DRO} + z^*_\text{RO})/2}$.}
\label{table:r&d-subopt}
\end{table}
From the table, we see that the average suboptimality of the RO solution among all instances is up to 23.03\%, which shows the importance of incorporating distributional information in the uncertainty set.

To further assess the performance of the DRO and RO solutions, in Appendix \ref{ecsubsub:r&d-out-sample}, we assess the out-of-sample performance of these solutions. The results show that the DRO approach, in general, performs better than RO on non-worst-case distributions.

\newqj{
\section{Conclusion}
\label{sec:conclusion}
In this paper, we studied two-stage distributionally robust optimization problems with decision-dependent information discovery. We proposed a framework for modeling these problems, a~$K$-adaptability scheme to approximately solve them, and a decomposition algorithm to exactly solve the resulting~$K$-adaptable problem. We demonstrated the effectiveness of our approaches in terms of computational efficiency, the benefits of adaptability, and out-of-sample performance, compared to the robust optimization framework, which assumes only support information is known.

There are several promising directions for future research. First, in this paper, we adopted a moment-based ambiguity set, motivated by the fact that data are not available during the design phase and only become available after implementing the first- and second-stage decisions. A natural extension would be to explore data-driven ambiguity sets in settings where data are available. This requires the development of alternative ambiguity sets, possibly in the spirit of \cite{Chen2020RobustRSOME} and \cite{zhang2024optimal}, which in turn would lead to different reformulations and corresponding solution methods. Second, it would be valuable to examine how solution quality varies with the parameter $K$. While there are some theoretical results in the literature regarding robust optimization, see \cite{kurtz2024many}, it is of practical importance to assess the performance of the $K$-adaptability method in the DRO problem with DDID. In our paper, we take a step in this direction through numerical experiments that evaluate the empirical improvement over the static solution as the number of policies increases. It would be interesting to further derive theoretical results, particularly regarding how the optimality gap varies with $K$ and how many policies are required to achieve optimality. Finally, it would be interesting to extend the proposed algorithm to problems with multiple time stages, as this would broaden its applicability and enable it to address a wider range of practical problems.

}

\section*{Acknowledgements}
Phebe Vayanos and Qing Jin are funded in part by the National Science Foundation under grant 1763108. Grani A. Hanasusanto is funded in part by the National Science Foundation under grants 2342505, 2343869, and 2404413. They are grateful for the support.

\clearpage
% \bibliography{references}

%%%%%%%%%%%%%%%%%%%%%%%%%%%%%%%%%%%%%%%%%%%%%%%%%%%%%%%%%%%%%%%%%%%%%%%%%%%%%%%%%%%%%
%%%%%%%%%%%%%%%%%%%%%%%%%%%%%%%%%%%%%%%%%%%%%%%%%%%%%%%%%%%%%%%%%%%%%%%%%%%%%%%%%%%%%
%%%%%%%%%%%%%%%%%%%%%%%%%%%%%%%%%%%%%%%% E-COMPANION %%%%%%%%%%%%%%%%%%%%%%%%%%%%%%%%
%%%%%%%%%%%%%%%%%%%%%%%%%%%%%%%%%%%%%%%%%%%%%%%%%%%%%%%%%%%%%%%%%%%%%%%%%%%%%%%%%%%%%
%%%%%%%%%%%%%%%%%%%%%%%%%%%%%%%%%%%%%%%%%%%%%%%%%%%%%%%%%%%%%%%%%%%%%%%%%%%%%%%%%%%%%

% Here starts the e-companion (EC)
%%%%%%%%%%%%%%%%%%%%%%%%%%%%%%%%%%%%%%%%%%%%%%%%%%%%%%%%%
\newpage

%\ECDisclaimer
%%%%%%%%%%%%%%%%%%%%%%%%%%%%%%%%%%%%%%%%%%%%%%%%%%%%%%%%%

%% Main head for the e-companion
\begin{appendix}
\setcounter{equation}{0}
\renewcommand{\theequation}{\thesection.\arabic{equation}}
\setcounter{theorem}{0}
\renewcommand{\thetheorem}{\thesection.\arabic{theorem}}
\setcounter{lemma}{0}
\renewcommand{\thelemma}{\thesection.\arabic{lemma}}
%%%%%%%%%%%%%%%%%%%%%%%%%%%%%%%%%%%%%%%%%%%%%%%%%%%%%%%%%
%%%%%%%%%%%%%%%%%%% OMITTED PROOFS %%%%%%%%%%%%%%%%%%%%%%
%%%%%%%%%%%%%%%%%%%%%%%%%%%%%%%%%%%%%%%%%%%%%%%%%%%%%%%%%

\section{Omitted proofs}
\label{ec:proofs}

The following lemma is used to prove the strong duality result in Theorem~\ref{thm:min-max}.
\begin{lemma}
The point~$(1,\;\bm c)$ resides in the interior of the convex cone
\begin{equation}
\mathcal{V}=\left\{(a, \bm b) \in \mathbb{R} \times \mathbb{R}^{N_g}: \exists \mu \in \mathcal{M}_{+}(\Xi) \text { such that } \begin{array}{ll} 
&\int \mu(\mathrm{d} \bm \xi)=a, \\
& \int \bm g(\bm \xi) \mu(\mathrm{d} \bm \xi) \leq \bm b
\end{array}\right\}.
\end{equation}
\label{lemma:interior}
\end{lemma}

{\proof
Let~$\mathbb B_r(\bm o)$ be the closed Euclidean ball of radius~$r \geq 0$ centered at~$\bm o$. Choose any point~$(s,\;\bm s)\in \mathbb B_{\kappa}(1) \times \mathbb B_{\kappa}(\bm c)$, where~$\kappa$ is a sufficiently small constant. Consider a point~$\tilde{\bm \xi} \in \interior \Xi \text{ such that } \bm g(\tilde{\bm \xi}) < \bm c$. The existence of such a point is guaranteed by the Slater condition. Because the function~$\bm g$ is continuous, we can always find an interior point~$\tilde{\bm \xi}$ of $\Xi$. Next, consider a scaled Dirac measure~$s \cdot \delta_{\tilde{\bm \xi}/s}$that places mass~$s$ at~$\tilde{\bm \xi}/s$. By construction, this measure satisfies~$\int s \cdot \delta_{\tilde{\bm \xi}/s} (\mathrm{d} \bm \xi) = s$. Moreover, for  sufficiently small~$\kappa$ the measure is supported on~$\Xi$ (since~$\tilde{\bm \xi} \in \interior \Xi$) and satisfies~$\int \bm g(\bm \xi) s \cdot \delta_{\tilde{\bm \xi}/s} (\mathrm{d} \bm \xi) \leq \bm s$ (since~$\bm g(\tilde{\bm \xi}) < \bm c$ and the function~$\bm g$ is continuous).} \endproof

\proof[Proof of Theorem~\ref{thm:min-max}] For fixed~$\bm{x}, \bm{w}, \bm{y}(\cdot)$, we can express the objective function of~\eqref{eq:endo_whole} as the optimal value of the moment problem
\begin{equation}
\begin{array}{ll}
\max &\displaystyle \int_{\Xi}  {\bm \xi}^\top {\bm C} \; {\bm x} + {\bm \xi}^\top {\bm D} \; {\bm w}  + {\bm \xi}^\top {\bm Q} \; {\bm y}({\bm \xi})\;\; \mu(\text{d} \bm\xi)\\

\st &  \mu\in\mathcal{M}_+(\mathbb{R}^{N_\xi})\\
&\displaystyle\int_{\Xi} \mu(\text{d} \bm\xi) = 1\\
& \displaystyle\int_{\Xi} \bm g(\bm \xi) \mu(\text{d} \bm\xi) \leq \bm c.\\
\end{array}
\label{primal}
\end{equation} 
The dual problem is given by
\begin{equation}
\begin{array}{ll}
\min & \phi + \bm c^\top \bm \psi \\
\st &  \phi\in\mathbb{R},\; \bm\psi\in\mathbb{R}_+^{N_g}\\
&\phi+ \bm\psi^\top \bm g(\bm\xi)\geq  {\bm \xi}^\top {\bm C} \; {\bm x} + {\bm \xi}^\top {\bm D} \; {\bm w}  + {\bm \xi}^\top {\bm Q} \; {\bm y}({\bm \xi})\qquad \forall\bm\xi\in\Xi.
\end{array}
\label{eq:dual}
\end{equation}
Strong duality holds by Lemma \ref{lemma:interior} and \citet[Proposition 3.4]{Shapiro2001OnProblems}.
Eliminating the variable $\phi$ and combining  with the outer minimization for~$(\bm x,\;\bm w,\; \bm y(\cdot))$, we get
\begin{equation}
    \begin{array}{cl}
    \displaystyle \min & \quad \displaystyle \max_{\bm \xi \in \Xi} \;\;\bm c^\top \bm \psi + {\bm \xi}^\top {\bm C} \; {\bm x} + {\bm \xi}^\top {\bm D} \; {\bm w}  + {\bm \xi}^\top {\bm Q} \; {\bm y}({\bm \xi}) - \bm\psi^\top \bm g(\bm \xi)\\

    \st & \quad \bm\psi\in\mathbb{R}_+^{N_g},\;{\bm x} \in \sets X, \; {\bm w} \in \sets W\\
    
    & \quad \!\! \left. \begin{array}{l} 
    {\bm y}({\bm \xi}) \in \sets Y  \\
    {\bm T}(\bm \xi) {\bm x} + {\bm V}(\bm \xi){\bm w} + {\bm W}(\bm \xi){\bm y}({\bm \xi}) \leq {\bm H}{\bm \xi} 
         \end{array} \quad \right\} \quad \forall {\bm \xi} \in \Xi \\
    & \quad {\bm y}({\bm \xi}) = {\bm y}({\bm \xi}') \quad \forall {\bm \xi}, {\bm \xi}' \in \Xi \; : \; {\bm w} \circ {\bm \xi} = {\bm w} \circ {\bm \xi}'.
    \end{array}
\label{eq:dual-full}
\end{equation}
The problem can be interpreted as a two-stage robust optimization with DDID, where the dual multipliers $\bm \psi$ constitute first-stage decisions.

We now show the equivalence between problem~\eqref{eq:nested} and~\eqref{eq:dual-full}. We first show that problem~\eqref{eq:nested} lower bounds~\eqref{eq:dual-full}. Let~$(\bm x, \; \bm w, \; \bm \psi,\;\bm y(\cdot))$ be an optimal solution of problem~\eqref{eq:dual-full}. By construction,~$(\bm x, \; \bm w, \; \bm \psi)$ is also feasible in~\eqref{eq:nested}. For any fixed~$\overline{\bm \xi}$, we set
$$
    \begin{array}{cl}
          \bm y^{\prime}(\overline{\bm \xi})  \in \displaystyle \argmin_{\bm y \in \mathcal Y} \left\{ 
    \begin{array}{cl}
             \displaystyle\max_{ {\bm \xi} \in \Xi({\bm w},\overline{\bm \xi}) } & \; \bm c^\top \bm \psi + {\bm \xi}^\top {\bm C} \; {\bm x} + {\bm \xi}^\top {\bm D} \; {\bm w}  + {\bm \xi}^\top {\bm Q} \; {\bm y} - \bm\psi^\top \bm g(\bm\xi) \\
             \st & \;\; {\bm T} ({\bm \xi}){\bm x} + {\bm V} ({\bm \xi}){\bm w} + {\bm W} ({\bm \xi}){\bm y} \leq {\bm H}{\bm \xi} \;\;\; \forall {\bm \xi} \in \Xi({\bm w},\overline{\bm \xi})
         \end{array}
    \right\}.
\end{array}
$$
It holds that
\begin{equation}
    \begin{array}{cl}
& \displaystyle \max_{\overline{\bm \xi} \in \Xi} \displaystyle \max _{\bm \xi \in \Xi(\bm w, \overline{\bm \xi})} \bm c^\top \bm \psi + \bm \xi^{\top} \bm{C} \bm{x}+\bm{\xi}^{\top} \bm{D} \bm{w}+\bm{\xi}^{\top} \bm{Q} \bm{y}^{\prime}(\overline{\bm{\xi}}) - \bm \psi^\top \bm g(\bm \xi)\\
= & \displaystyle \max _{\bm{\xi} \in \Xi} \bm c^\top \bm \psi + \bm{\xi}^{\top} \bm{C} \bm{x}+\bm{\xi}^{\top} \bm{D} \bm{w}+\bm{\xi}^{\top} \bm{Q} \bm{y}^{\prime}(\bm{\xi}) - \bm \psi^\top \bm g(\bm \xi) \\
\leq & \displaystyle \max _{\bm{\xi} \in \Xi} \bm c^\top \bm \psi + \bm{\xi}^{\top} \bm{C} \bm{x}+\bm{\xi}^{\top} \bm{D} \bm{w}+\bm{\xi}^{\top} \bm{Q} \bm{y}(\bm{\xi}) - \bm \psi^\top \bm g(\bm \xi)
\end{array}
\end{equation}
Thus, we have shown that the optimal solution of~\eqref{eq:dual-full} is feasible in~\eqref{eq:nested} with an objective value no greater than the one attained in~\eqref{eq:dual-full}. 

Next, we show conversely that  problem~\eqref{eq:dual-full} lower bounds~\eqref{eq:nested}. To this end, let~$(\bm x, \bm w, \bm \psi)$ be an optimal solution in~\eqref{eq:nested} and define~$\bm y(\bm \xi) = \bm y^\prime(\bm w \circ \bm \xi)$. By construction, the quadruple~$(\bm x, \; \bm w, \; \bm \psi, \; \bm y(\bm \xi))$ is feasible in~\eqref{eq:dual-full} and attains the objective value
\begin{equation}
    \begin{array}{cl}
& \displaystyle \max _{\bm{\xi} \in \Xi} \bm c^\top \bm \psi + \bm{\xi}^{\top} \bm{C} \bm{x}+\bm{\xi}^{\top} \bm{D} \bm{w}+\bm{\xi}^{\top} \bm{Q} \bm{y}(\bm{\xi}) - \bm \psi^\top \bm g(\bm \xi) \\
= & \displaystyle \max_{\overline{\bm \xi} \in \Xi} \displaystyle \max _{\bm \xi \in \Xi(\bm w, \overline{\bm \xi})} \bm c^\top \bm \psi + \bm \xi^{\top} \bm{C} \bm{x}+\bm{\xi}^{\top} \bm{D} \bm{w}+\bm{\xi}^{\top} \bm{Q} \bm{y}^{\prime}(\overline{\bm{\xi}}) - \bm \psi^\top \bm g(\bm \xi)
\end{array}
\end{equation}
which is the same value attained in~\eqref{eq:nested}.

Thus, we have shown that the optimal objective values attained by the two problems coincide. And we can construct an optimal solution in problem~\eqref{eq:dual-full} from an optimal solution in~\eqref{eq:nested}. This completes the proof.\endproof

\proof[Proof of Theorem~\ref{thm:objunc-bilinear}]
We linearize the vector~$\bm g(\bm \xi^k),\;\forall k \in \mathcal{K}$ of piecewise linear functions by introducing an epigraph uncertain parameter~$\zeta^k_s$ for~$g_s(\bm \xi^k),\;\forall s \in \mathcal{S}, \forall k \in \mathcal{K}\;$\newqj{,} and write the inner maximization problem in the following epigraph form
\begin{equation}
        \begin{array}{cl}
         \max & \quad \tau   \\
         \st & \quad \tau \in \reals, \; \overline{\bm \xi} \in \reals^{N_\xi}, \; {\bm \xi}^k \in \reals^{N_\xi}, k \in \sets K \\
         & \quad  \tau \; \leq \; ({\bm \xi}^k)^\top {\bm C} \; {\bm x} + ({\bm \xi}^k)^\top {\bm D} \; {\bm w}  + ({\bm \xi}^k)^\top {\bm Q} \; {\bm y^k}  - \bm\psi^\top \bm \zeta^k \quad \forall k \in \sets K \\
         & \quad {\bm A}\overline{\bm \xi} \leq {\bm b} \\
         & \quad {\bm A}{\bm \xi}^k \leq {\bm b}  \quad \forall k \in \sets K \\
         & \quad  {\bm w} \circ {\bm \xi}^k = {\bm w} \circ \overline{\bm \xi}  \quad \forall k \in \sets K \\
         & \quad \zeta_s^k \leq -\bm g_{s,t}^\top \bm \xi^k \qquad \forall s \in \sets S, \, \forall t\in\sets T,\,\forall k \in \sets K.
    \end{array}
    \label{subproblem1}
\end{equation} 
Problem~\eqref{subproblem1} is a feasible LP problem, so strong duality implies that it has the same optimal value as its dual problem. Dualizing and grouping with the outer minimization yields problem~\eqref{eq:objunc-bilinear} in Theorem~\ref{thm:objunc-bilinear}.\endproof

\proof[Proof of Theorem~\ref{thm:cstrunc-bilinear}]

Applying the uncertainty set lifting procedure in~\citet[Lemma 1]{Vayanos2020} to problem~\eqref{eq:k-nested}, we can lift the space of the uncertainty set in the inner maximization, exchange it with the inner minimization, and finally get the equivalent min-max-min problem
\begin{equation}
    \begin{array}{cl}
         \min & \bm c^\top \bm \psi +\displaystyle \max_{ \{ {\bm \xi}^k \}_{k \in \sets K} \in \Xi^K({\bm w}) }  \left\{ 
         \begin{array}{cl}
             \displaystyle\min_{ k\in\sets K} & ({\bm \xi}^k)^\top {\bm C}  {\bm x} + ({\bm \xi}^k)^\top {\bm D}  {\bm w}  + ({\bm \xi}^k)^\top {\bm Q}  {\bm y^k} - \bm\psi^\top \bm g(\bm \xi^k) \\
            \st & {\bm T}(\bm \xi) {\bm x} + {\bm V}(\bm \xi) {\bm w} + {\bm W}(\bm \xi){\bm y}^k \leq {\bm H}{\bm \xi}
         \end{array}
          \right\}  \\
         \st & \; \bm\psi\in\mathbb{R}_+^{N_g},\;{\bm x} \in \sets X, \; {\bm w} \in \sets W,\;\bm y^k \in\mathcal{Y}\quad \forall k\in\sets K,
    \end{array}
    \label{eq:k-objunc-cstr}
\end{equation}
where
\begin{equation}
\Xi^K({\bm w}) := 
\left\{
\{ {\bm \xi}^k\}_{k \in \sets K} \in \Xi : \exists \;\overline {\bm \xi} \in \Xi \;\; \textup{such\;that}\;{\bm \xi}^k \in \Xi( {\bm w},\overline {\bm \xi})\;\;\forall k \in \sets K
\right\}.
\label{eq:lift-unc}
\end{equation}

We then shift the second stage constraints into the uncertainty set and formulate the~$K$-adaptability problem with constraint uncertainty~\eqref{eq:k-nested} equivalently as
\begin{equation}
    \begin{array}{cl}
         \min & \;\; \bm c^\top \bm \psi +\displaystyle \max _{\bm \ell \in \mathcal{L}} \max _{\{\bm \xi^{k}\}_{k \in \mathcal{K}}\in \Xi^{K}({\bm w}, \bm \ell)} \min_{\genfrac{}{}{0pt}{2}{k \in K:}{\bm \ell_k=0}} \; 
         \begin{array}[t]{r}
         \left\{ ({\bm \xi}^k)^\top {\bm C} \; {\bm x} + ({\bm \xi}^k)^\top {\bm D} \; {\bm w}  + ({\bm \xi}^k)^\top {\bm Q} \; {\bm y} \right. \\
         - \left. \bm\psi^\top \bm g(\bm \xi^k)
          \right\}
          \end{array} \\
         \st & \;\;  \bm\psi\in\mathbb{R}_+^{N_g},\;{\bm x} \in \sets X, \; {\bm w} \in \sets W,\;\bm y^k \in\mathcal{Y}\quad \forall k\in\sets K,
    \end{array}
    \label{eq:k-cstr-l}
\end{equation}
where~$\sets L:=\{1, \ldots, L\}^{K}$,~$L$ is the number of second-stage constraints with uncertainty in problem~\eqref{eq:k-nested}. The uncertainty sets~$\Xi^K(\bm w, \bm \ell), \; \bm \ell \in \sets L$, are defined as
\begin{equation}
\Xi^{K}(\boldsymbol{w}, \boldsymbol{\ell}):=\left\{\begin{array}{ll}
&\boldsymbol{w} \circ \boldsymbol{\xi}^{k}=\boldsymbol{w} \circ \overline{\boldsymbol{\xi}}, \;\;
\forall k \in \mathcal{K} \text { for some } \overline{\boldsymbol{\xi}} \in \Xi 
\\
\left\{\boldsymbol{\xi}^{k}\right\}_{k \in \mathcal{K}} \in \Xi^{K}:
 &
 \begin{array}{r}
 \boldsymbol{T}(\bm \xi^k) \boldsymbol{x} +\boldsymbol{V}(\bm \xi^k) \boldsymbol{w}  +\boldsymbol{W}(\bm \xi^k) \boldsymbol{y}^{k} \leq \boldsymbol{H} \boldsymbol{\xi}^{k} \\   \forall k \in \mathcal{K}: \bm \ell_k=0 
 \end{array}
\\
&
\begin{array}{r}
\left[\boldsymbol{T}(\bm \xi^k) \boldsymbol{x}  + \boldsymbol{V}(\bm \xi^k) \boldsymbol{w}  +\boldsymbol{W}(\bm \xi^k) \boldsymbol{y}^{k}\right]_{\bm \ell_{k}}> \left[\boldsymbol{H} 
\boldsymbol{\xi}^{k}\right]_{\bm \ell_{k}} \\ \forall k \in \mathcal{K}: \bm \ell_k \neq 0 
\end{array}
\end{array}\right\}.
\label{eq:uncset-l}
\end{equation}
\label{pro:uncset-l}
The components of vector~$\bm \ell \in \sets L$ encode which policies are robust feasible for the
parameter realizations~$\left\{\boldsymbol{\xi}^{k}\right\}_{k \in \mathcal{K}}$. Policy~$\bm y^k$ is robust feasible in problem~\eqref{eq:k-nested} if~$\ell_k = 0$. If~$\ell_k \neq 0$, then the~$\ell_k$-th constraint in problem~\eqref{eq:k-nested} is violated for some~$\bm \xi \in \Xi(\bm w, \overline{\bm \xi})$.

Problem~\eqref{eq:k-cstr-l} involves many open uncertainty sets. Thus, we employ inner approximations~$\Xi_\epsilon^{K}(\boldsymbol{w}, \boldsymbol{\ell})$ of the uncertainty sets~\eqref{eq:uncset-l} that are parameterized by scalar~$\epsilon$. The approximate problem can be written as
\begin{equation}
    \begin{array}{cl}
         \min & \;\; \bm c^\top \bm \psi +\displaystyle \max _{\bm \ell \in \mathcal{L}} \max _{\left\{\bm \xi^{k}\right\}_{k \in \mathcal{K}}\in \Xi_\epsilon^{K}({\bm w}, \bm \ell)} \min_{\genfrac{}{}{0pt}{2}{k \in K:}{\bm \ell_k=0}} \; 
         \begin{array}[t]{r}
         \left\{ 
         ({\bm \xi}^k)^\top {\bm C} \; {\bm x} + ({\bm \xi}^k)^\top {\bm D} \; {\bm w}  + ({\bm \xi}^k)^\top {\bm Q} \; {\bm y} \right. \\ 
         \left.- \bm\psi^\top \bm g(\bm\xi^k)
          \right\} 
          \end{array} \\
         \st & \;\;  \bm\psi\in\mathbb{R}_+^{N_g},\;{\bm x} \in \sets X, \; {\bm w} \in \sets W,\;\bm y^k \in\mathcal{Y}\quad \forall k\in\sets K,
    \end{array}
    \label{eq:k-cstr-eps}
\end{equation}
where
\begin{equation*}
\Xi_\epsilon^{K}(\boldsymbol{w}, \boldsymbol{\ell}):=\left\{\begin{array}{ll}
&\boldsymbol{w} \circ \boldsymbol{\xi}^{k}=\boldsymbol{w} \circ \overline{\boldsymbol{\xi}}, \;\; \forall k \in \mathcal{K} : \exists\, \overline{\boldsymbol{\xi}} \in \Xi \\
\left\{\boldsymbol{\xi}^{k}\right\}_{k \in \mathcal{K}} \in \Xi^{K}: &
\begin{array}{r}
\boldsymbol{T}(\bm \xi^k) \boldsymbol{x}+\boldsymbol{V}(\bm \xi^k) \boldsymbol{w}+\boldsymbol{W}(\bm \xi^k) \boldsymbol{y}^{k} \leq \boldsymbol{H} \boldsymbol{\xi}^{k} \\ \forall k \in \mathcal{K}: \bm \ell_k=0 
\end{array}
\\
&
\begin{array}{r}
\left[\boldsymbol{T}(\bm \xi^k) \boldsymbol{x}+\boldsymbol{V}(\bm \xi^k) \boldsymbol{w}+\boldsymbol{W}(\bm \xi^k) \boldsymbol{y}^{k}\right]_{\bm \ell_{k}} \\ \geq \left[\boldsymbol{H} \boldsymbol{\xi}^{k}\right]_{\bm \ell_{k}} +\epsilon  \\ \forall k \in \mathcal{K}: \bm \ell_k \neq 0
\end{array}
\end{array}\right\}
\label{eq:uncset-eps}
\end{equation*}
We next reformulate the approximate problem~\eqref{eq:k-cstr-eps} as a mixed binary bilinear program.
The outer maximization of the approximation problem~\eqref{eq:k-cstr-eps} is identical to
\begin{equation}
    \max _{\bm \ell \in \mathcal{L}} \max _{\left\{\bm \xi^{k}\right\}_{k \in \mathcal{K}} \in \Xi_{\epsilon}^{K}(\boldsymbol{w}, \boldsymbol{\bm \ell})} \min _{\boldsymbol{\delta} \in \Delta_{K}(\ell)}
    \begin{array}[t]{r}
    \left\{\displaystyle \sum_{k \in \mathcal{K}} \boldsymbol{\delta}_{k}\left[\left(\boldsymbol{\xi}^{k}\right)^{\top} \boldsymbol{C} \boldsymbol{x}+\left(\boldsymbol{\xi}^{k}\right)^{\top} \boldsymbol{D} \boldsymbol{w}+\left(\boldsymbol{\xi}^{k}\right)^{\top} \boldsymbol{Q} \boldsymbol{y}^{k} \right. \right. \\
    \left. \left. - \bm \psi^\top \bm g(\bm \xi)\right]\right\},
    \end{array}
\end{equation}
 where $\Delta_{K}(\bm \ell):=\left\{\boldsymbol{\delta} \in \mathbb{R}_{+}^{K}: \mathbf{e}^{\top} \boldsymbol{\delta}=1, \boldsymbol{\delta}_{k}=0 \quad \forall k \in \mathcal{K}: \ell_{k} \neq 0\right\}$. If~$\Xi_\epsilon^K(\bm w, \bm \ell) = \emptyset$ for all~$\bm \ell \in \sets L_+$, apply the min-max theorem, the problem is equivalent to
\begin{equation}
    \min _{\bm \delta(\ell) \in \Delta_{K}(\ell),} \max _{\ell \in \partial \mathcal{L}} \max _{\left\{\xi^{k}\right\}_{k \in \mathcal{K}} \in \Xi_{\epsilon}^{K}(w, \ell)}    \begin{array}[t]{r}
    \left\{\displaystyle \sum_{k \in \mathcal{K}} \boldsymbol{\delta}_{k}\left[\left(\boldsymbol{\xi}^{k}\right)^{\top} \boldsymbol{C} \boldsymbol{x}+\left(\boldsymbol{\xi}^{k}\right)^{\top} \boldsymbol{D} \boldsymbol{w}+\left(\boldsymbol{\xi}^{k}\right)^{\top} \boldsymbol{Q} \boldsymbol{y}^{k} \right. \right. \\
    \left. \left. - \bm \psi^\top \bm g(\bm \xi)\right]\right\}.
    \end{array}
\end{equation}
Use an epigraph formulation, we conclude that~\eqref{eq:k-cstr-eps} is equivalent to the following formulation
\begin{equation}
        \begin{array}{cl}
         \max & \quad \tau   \\
         \st & \quad x \in \mathcal{X}, \; \boldsymbol{w} \in \mathcal{W}, \; \boldsymbol{y}^{k} \in \mathcal{Y}, \; k \in \mathcal{K}, \; \bm \psi \in \reals^{N_g}_+ \\
         & \quad \tau \in \mathbb{R}, \; \boldsymbol{\bm \delta}(\bm \ell) \in \bm \Delta_{K}(\bm \ell), \; \bm \ell \in \partial \mathcal{\bm L}\\
         & \quad 
         \begin{array}{lr}
         \tau \geq \sum_{k \in \mathcal{K}} \boldsymbol{\bm \delta}_{k}(\boldsymbol{\bm \ell})\left[\left(\boldsymbol{\xi}^{k}\right)^{\top} \boldsymbol{C} \boldsymbol{x} \right.& \left.  +\left(\boldsymbol{\xi}^{k}\right)^{\top} \boldsymbol{D} \boldsymbol{w}+\left(\boldsymbol{\xi}^{k}\right)^{\top} \boldsymbol{Q} \boldsymbol{y}^{k}- \bm \psi^\top \bm g(\bm \xi) \right] \\ 
         & \forall \bm \ell \in \partial \mathcal{L},\left\{\boldsymbol{\xi}^{k}\right\}_{k \in \mathcal{K}} \in \Xi_{\epsilon}^{K}(\boldsymbol{w}, \boldsymbol{\ell}) 
         \end{array} \\
        & \quad \Xi_{\epsilon}^{K}(\boldsymbol{w}, \boldsymbol{\ell})=\emptyset \quad \forall \bm \ell \in \mathcal{L}_{+}.
    \end{array}
    \label{eq:cstr-epi}
\end{equation}
The semi-infinite constraint associated with~$\bm \ell \in \partial \mathcal{L}$ is satisfied if and only if the optimal value of
\begin{equation}
    \begin{array}{cl}
    \max & \sum_{k \in \mathcal{K}} \boldsymbol{\bm \delta}_{k}(\boldsymbol{\ell})\left[\left(\boldsymbol{\xi}^{k}\right)^{\top} \boldsymbol{C} \boldsymbol{x}+\left(\boldsymbol{\xi}^{k}\right)^{\top} \boldsymbol{D} \boldsymbol{w}+\left(\boldsymbol{\xi}^{k}\right)^{\top} \boldsymbol{Q} \boldsymbol{y}^{k} + \bm\psi^\top \bm \zeta^k \right] \\
    \st & \overline{\boldsymbol{\xi}} \in \mathbb{R}^{N_{\xi}}, \boldsymbol{\xi}^{k} \in \mathbb{R}^{N_{\xi}}, k \in \mathcal{K} \\
    & \boldsymbol{A} \overline{\boldsymbol{\xi}} \leq \boldsymbol{b} \\
    & \boldsymbol{A} \boldsymbol{\xi}^{k} \leq \boldsymbol{b} \quad \forall k \in \mathcal{K} \\
    & \boldsymbol{T}(\bm \xi^k) \boldsymbol{x}+\boldsymbol{V}(\bm \xi^k) \boldsymbol{w}+\boldsymbol{W}(\bm \xi^k) \boldsymbol{y}^{k} \leq \boldsymbol{H} \boldsymbol{\xi}^{k} \quad \forall k \in \mathcal{K}: \ell_k=0 \\
    & {\left[\boldsymbol{T}(\bm \xi^k) \boldsymbol{x}+\boldsymbol{V}(\bm \xi^k) \boldsymbol{w}+\boldsymbol{W}(\bm \xi^k) \boldsymbol{y}^{k}\right]_{\ell_{k}} \geq\left[\boldsymbol{H} \boldsymbol{\xi}^{k}\right]_{\ell_{k}}+\epsilon} \quad \forall k \in \mathcal{K}: \ell_k \neq 0 \\
    & \boldsymbol{w} \circ \boldsymbol{\xi}^{k}=\boldsymbol{w} \circ \overline{\boldsymbol{\xi}} \quad \forall k \in \mathcal{K} \\
    & \bm \zeta_s^k \leq -\bm g_{s,t}^\top \bm \xi^k \quad \forall s \in \sets S, \, \forall t\in\sets T,\;\forall k\in\sets K
    \end{array}
    \label{eq:partial-L}
\end{equation}
does not exceed~$\tau$.

The last constraint in~\eqref{eq:cstr-epi} is satisfied for~$\bm \ell \in \sets L$ whenever the linear program problem
\begin{equation}
        \begin{array}{cll}
         \max & \quad 0  & \\
         \st & \quad \overline{\bm \xi} \in \reals^{N_\xi}, \; {\bm \xi}^k \in \reals^{N_\xi}, k \in \sets K \\
         & \quad {\bm A}\overline{\bm \xi} \leq {\bm b} & \\
         & \quad {\bm A}{\bm \xi}^k \leq {\bm b}  & \forall k \in \sets K \\
         & \quad  {\bm w} \circ {\bm \xi}^k = {\bm w} \circ \overline{\bm \xi}  & \forall k \in \sets K\\
         & \quad {\left[\boldsymbol{T}(\bm \xi^k) \boldsymbol{x}+\boldsymbol{V}(\bm \xi^k) \boldsymbol{w}+\boldsymbol{W}(\bm \xi^k) \boldsymbol{y}^{k}\right]_{\ell_{k}}\geq\left[\boldsymbol{H}\boldsymbol{\xi}^{k}\right]_{\ell_{k}}+\epsilon} & \forall k \in \sets K\\
    \end{array}
    \label{eq:l+}
\end{equation}
is infeasible.
Strong duality holds for both~\eqref{eq:partial-L} and~\eqref{eq:l+}. After taking dualization of~\eqref{eq:partial-L},~\eqref{eq:l+},\newqj{,} and combining the dual problem with the outer minimization problem, we can get the following problem
 \begin{equation}
\begin{array}{cl}
\min & \bm c^{\top}\bm \psi + \mathcal{\tau}\\
\st & x \in \mathcal{X}, \; \bm w \in \mathcal{W}, \; \bm y^k \in \sets Y, \; k \in \sets K, \; \tau \in \reals, \; \bm \psi \in \reals^{N_g}_+\\

& \left.\begin{array}{l} 
\bm \alpha(\bm \ell) \in \reals_{+}^{R}, \; \bm \alpha^{k}(\bm \ell) \in \reals_{+}^{R}, \; k \in \sets K,  \bm \eta^{k}(\bm \ell) \in \reals^{N_{\xi}}, \; k \in \sets K \\

\bm \Lambda^k(l) \in \reals^{N_g \times T}_+ \; k \in \sets K,\; \bm \beta^{k}(\bm \ell) \in \mathbb{R}_{+}^{L}, \; k \in \mathcal{K}\\

\bm \delta(\bm \ell) \in \Delta_{K}(\bm \ell),\; \bm \gamma(\bm \ell) \in \mathbb{R}_{+}^{K} \\

\bm A^{\top} \bm \alpha(\bm \ell)=\sum\limits_{k \in \mathcal{K}} \bm w \circ \bm \eta^{k}(\bm \ell) \\
\bm \delta_k(\bm \ell)\bm \psi = \bm \Lambda^k(\bm \ell)\bm e\\

\begin{array}{r}
\bm \delta_{k}(\bm \ell) \left[\bm {C x + D w + Q}\bm y^{k}\right] - \sum^L_{l = 1}\left(\bm{T}_l \bm x+ \bm V_l \bm w+ \bm W_l \bm y^{k}\right)\bm \beta_l^{k}(\bm \ell) \\ 
 = \bm A^{\top} \bm \alpha^{k}(\bm \ell)- \bm H^{\top} \bm \beta^{k}(\bm \ell) + \bm w \circ \bm \eta^{k}(\bm \ell) + \sum\limits_{s \in \sets S}\sum\limits_{t \in \sets T}\bm \Lambda_{st}(\bm \ell)\bm g_{st} \\
\;\; \forall k \in \mathcal{K}: \bm \ell_{k}=0 
\end{array}\\

\begin{array}{r}
\bm \delta_{k}(\bm \ell) \left[\bm {C x + D w + Q}\bm y^{k}\right] + \left(\bm{T}_{\bm \ell_{k}} \bm x+ \bm V_{\bm \ell_k} \bm w+ \bm W_{\bm \ell_k} \bm y^{k}\right)\bm \gamma_{k}(\bm \ell) \\
 = \bm A^{\top} \bm \alpha^{k}(\bm \ell) +[\bm H]_{\bm \ell_{k}} \bm \gamma_{k}(\bm \ell) + \bm w \circ \bm \eta^{k}(\bm \ell) + \sum\limits_{s\in\sets S}\sum\limits_{t\in\sets T}\bm \Lambda_{st}(\bm \ell)\bm g_{st} \\  \forall k \in \mathcal{K}: \bm \ell_{k} \neq 0 
\end{array}\\

% \begin{array}{ll}
\tau \geq \bm b^{\top}\left(\bm \alpha(\bm \ell)+\sum_{k \in \mathcal{K}} \bm \alpha^{k}(\bm \ell)\right) -\epsilon \sum\limits_{\genfrac{}{}{0pt}{2}{k \in K:}{\ell_k \neq 0}} \bm \gamma_{k}(\bm \ell)
% \end{array}

\end{array}\right\} \; \forall \bm \ell \in \partial \mathcal{L}\\

& \left.\begin{array}{l} 

\bm \phi(\bm \ell) \in \reals_{+}^{R}, \; \bm \phi^{k}(\bm \ell) \in \reals_{+}^{R}, \; k \in \sets K,  \; \bm \rho(\bm \ell) \in \reals_+^k, \; \bm \chi^{k}(\bm \ell) \in \reals^{N_{\xi}}, \; k \in \sets K\\

\bm A^{\top} \bm \phi(\bm \ell)=\sum\limits_{k \in \mathcal{K}} \bm w \circ \bm \chi^{k}(\bm \ell) \\

\begin{array}{r}
\bm A^{\top} \bm \phi^{k}(\bm \ell)+[\bm H]_{\bm \ell_{k}} \bm \rho_{k}(\bm \ell) + \bm w \circ \bm \chi^{k}(\bm \ell) \\ =\left(\bm{T}_{\bm \ell_{k}} \bm x+ \bm V_{\bm \ell_k} \bm w+ \bm W_{\bm \ell_k} \bm y^{k}\right) \bm \rho_{k}(\bm \ell) \\ \forall k \in \mathcal{K} 
\end{array}
\\

\bm b^{\top}\left(\bm \phi(\bm \ell)+\sum\limits_{k \in \mathcal{K}} \bm \phi^{k}(\bm \ell)\right)-\epsilon \sum\limits_{k \in \mathcal{K}} \bm \rho_{k}(\bm \ell) \leq-1

\end{array}\right\} \; \forall \bm \ell \in \mathcal{L}_+,

\end{array}
\label{eq:cstrunc-bilinear}
\end{equation}
where~$\sets L:=\{1, \ldots, L\}^{K}$,~$\Delta_{K}(\bm \ell):=\left\{\boldsymbol{\delta} \in \mathbb{R}_{+}^{K}: \mathbf{e}^{\top} \boldsymbol{\delta}=1, \boldsymbol{\delta}_{k}=0 \quad \forall k \in \mathcal{K}: \bm \ell_k \neq 0\right\},$ $ \partial \mathcal{L}:=\{\bm \ell \in \mathcal{L}: \bm \ell \not>\mathbf{0}\}$ and  $\mathcal{L}^{+}:=\{\bm \ell \in \mathcal{L}: \bm \ell>\mathbf{0}\}.$ Matrices~$\bm T_l,\;\bm V_l,\;\bm W_l$ are the coefficient matrices in the~$l$th constraint, in other words,~$[\bm T(\bm \xi)]_l = \bm \xi^\top \bm T_l,\; [\bm V(\bm \xi)]_l = \bm \xi^\top \bm V_l,\; [\bm W(\bm \xi)]_l = \bm \xi^\top \bm W_l$. \endproof

\proof[Proof of Proposition~\ref{prop: disjunctive}]
For any feasible solution~$(\bm \phi, \, \bm x,\,\{\bm y^k\}_{k \in \sets K})$ in~\eqref{eq:disjuctive}, we show that it is also feasible in~\eqref{eq:Phi(w)} with the same or lower objective value. Let  $\{\Xi^k\}_{k \in \sets K}$  be the corresponding optimal partition. Since $\{\Xi^k\}_{k \in \sets K}$ exhausts all points in $\Xi$, we see that $(\bm \phi, \, \bm x,\,\{\bm y^k\}_{k \in \sets K})$ is indeed feasible in~\eqref{eq:Phi(w)} because for any $\overline{\bm\xi}\in\Xi$, there must exist $k\in\mathcal K$ such that $\overline{\bm\xi}\in\Xi^k$ and, therefore, $\bm y^k$ is feasible in the inner minimization problem. Next, we show that the optimal objective value of \eqref{eq:Phi(w)} is at most $\theta$ by bounding it from above, written as 
\begin{equation}
    \begin{array}{cl}
       &  \displaystyle \max_{\overline{\bm \xi} \in \Xi} \displaystyle \left\{ 
         \begin{array}{cl}
            \displaystyle\min_{ k \in \sets K } \;\;
             &\displaystyle\max_{ {\bm \xi} \in \Xi({\bm w},\overline{\bm \xi}) } \;  \bm c^\top \bm \psi + {\bm \xi}^\top {\bm C} \; {\bm x} + {\bm \xi}^\top {\bm D} \; {\bm w}  + {\bm \xi}^\top {\bm Q} \; {\bm y^k} - \bm\psi^\top \bm g(\bm\xi) \\
             \st & \;\; {\bm T} ({\bm \xi}){\bm x} + {\bm V} ({\bm \xi}){\bm w} + {\bm W} ({\bm \xi}){\bm y^k} \leq {\bm H}{\bm \xi} \;\;\; \forall {\bm \xi} \in \Xi({\bm w},\overline{\bm \xi})
         \end{array}
          \right\} \\
  = &  \displaystyle \max_{j\in\mathcal K}\max_{\overline{\bm \xi} \in \Xi^j} \displaystyle \left\{ 
         \begin{array}{cl}
            \displaystyle\min_{ k \in \sets K } \;\;
             &\displaystyle\max_{ {\bm \xi} \in \Xi({\bm w},\overline{\bm \xi}) } \; \bm c^\top \bm \psi +  {\bm \xi}^\top {\bm C} \; {\bm x} + {\bm \xi}^\top {\bm D} \; {\bm w}  + {\bm \xi}^\top {\bm Q} \; {\bm y^k} - \bm\psi^\top \bm g(\bm\xi) \\
             \st & \;\; {\bm T} ({\bm \xi}){\bm x} + {\bm V} ({\bm \xi}){\bm w} + {\bm W} ({\bm \xi}){\bm y^k} \leq {\bm H}{\bm \xi} \;\;\; \forall {\bm \xi} \in \Xi({\bm w},\overline{\bm \xi})
         \end{array}
          \right\} \\
            \leq &  \displaystyle \max_{j\in\mathcal K}\max_{\overline{\bm \xi} \in \Xi^j} \displaystyle 
             \displaystyle\max_{ {\bm \xi} \in \Xi({\bm w},\overline{\bm \xi}) } \; \bm c^\top \bm \psi +  {\bm \xi}^\top {\bm C} \; {\bm x} + {\bm \xi}^\top {\bm D} \; {\bm w}  + {\bm \xi}^\top {\bm Q} \; {\bm y^j} - \bm\psi^\top \bm g(\bm\xi) \\\leq & \theta.
          \end{array}
\end{equation}
Here, the first inequality holds because we do not optimize for the best index $k$ but instead simply set it to the index $j$ from the outer maximization. 

Conversely, we show that any solution~$(\bm \phi, \, \bm x,\,\{\bm y^k\}_{k \in \sets K})$ feasible in~\eqref{eq:Phi(w)} is also feasible in~\eqref{eq:disjuctive} with the same or lower objective value. To this end, we construct a  partition  feasible for   \eqref{eq:disjuctive}, written as
\begin{equation*}
\begin{array}{r}
\Xi^j=\left\{\bm\xi\in\Xi:               
\begin{array}{cl}
            \displaystyle j\in\argmin_{ {k} \in \sets K } \;\;
             &\displaystyle\max_{ {\bm \xi} \in \Xi({\bm w},\overline{\bm \xi}) } \;  {\bm \xi}^\top {\bm C} \; {\bm x} + {\bm \xi}^\top {\bm D} \; {\bm w}  + {\bm \xi}^\top {\bm Q} \; {\bm y^k} - \bm\psi^\top \bm g(\bm\xi) \\
             \quad\;\;\st & \;\; {\bm T} ({\bm \xi}){\bm x} + {\bm V} ({\bm \xi}){\bm w} + {\bm W} ({\bm \xi}){\bm y^k} \leq {\bm H}{\bm \xi} \;\;\; \forall {\bm \xi} \in \Xi({\bm w},\overline{\bm \xi})
         \end{array}\right\}\\\forall j\in\mathcal K.
\end{array}
\end{equation*}
The set $\Xi^j$ contains all points $\bm\xi\in\Xi$ for which the policy $\bm y^j$ is optimal in the second-stage subproblem. 
Under this partition, one can verify that the solution $(\bm \phi, \, \bm x,\,\{\bm y^k\}_{k \in \sets K})$ is feasible in~\eqref{eq:disjuctive} with the same objective value. Since the partition can still be optimized, we conclude that the optimal value of problem~\eqref{eq:disjuctive} is at most equal to that of problem~\eqref{eq:Phi(w)}, which completes the proof. \endproof

\begin{theorem}
Problem~\eqref{eq:separation} is equivalent to the following MIO problem
\begin{equation}
    \begin{array}{cl}
         \max & \quad z \\
         \st & \quad  z \in \mathbb{R}, \; \overline{\bm \xi} \in \reals^{N_\xi}, \; {\bm \xi}^k \in \reals^{N_\xi}, \forall k \in \sets K, \; {z}_{kl} \in\{0,1\},(k, l) \in \mathcal{K} \times\{0,1, \ldots, L\} \\
         & \quad {\bm A}\overline{\bm \xi} \leq {\bm b} \\
         & \quad \!\! \left. \begin{array}{l} 
         {\bm A}{\bm \xi}^k \leq {\bm b}\\
         {\bm w^\prime} \circ {\bm \xi}^k = {\bm w^\prime} \circ \overline{\bm \xi}\\
        {z}_{k0}=1 \; \Rightarrow \; 
        \begin{array}[t]{r}
        z \leq \bm c^\top \bm \psi^\prime + { {\bm \xi}^k}^\top {\bm C} \; {\bm x^\prime} + {{\bm \xi}^k}^\top {\bm D} \; {\bm w^\prime}  + {{\bm \xi}^k}^\top {\bm Q} \; {\bm y^k}^\prime \\ 
        - {{\bm\psi}^\prime}^\top \bm g({\bm \xi}^k) - \theta 
        \end{array}
        \\
        z_{k l}=1 \; \Rightarrow \; 
        \begin{array}[t]{r}
        z \leq {{\bm \xi}^k}^\top {\bm T_l}{\bm x^\prime} + {{\bm \xi}^k}^\top{\bm V_l}{\bm w^\prime} + {{\bm \xi}^k}^\top{\bm W}_l{\bm y^k}^\prime -{\bm H_l}{\bm \xi} \\ \forall l \in\{1, \ldots, L\}
        \end{array}
        \\
        \sum_{l=0}^{L} {z}_{kl}=1
         \end{array} \quad \right\} \quad \forall k \in\sets K.
    \end{array}
\label{eq:separation-milp}
\end{equation}
\label{thm:sep-mip}
\end{theorem}

\proof
Following the same uncertainty lifting procedure that we used in the proof of Theorem~\ref{thm:cstrunc-bilinear}, we can write~\eqref{eq:separation} equivalently as
\begin{equation}
    \begin{array}{r}
        \displaystyle \max_{ \{ {\bm \xi}^k \}_{k \in \sets K} \in \Xi^K({\bm w^\prime})} \;\; \min_{ k \in \sets K } \;\; \max \left\{ \bm c^\top \bm \psi^\prime + ({\bm \xi}^k)^\top {\bm C} \; {\bm x^\prime} + ({\bm \xi}^k)^\top {\bm D} \; {\bm w^\prime}  + ({\bm \xi}^k)^\top {\bm Q} \; {{\bm y}^k}^\prime \right. \\ - {{\bm\psi}^\prime}^\top \bm g(\bm\xi^k) - \theta^\prime \; ,\\
         \;\; \left. \max\limits_{l\in\left\{1,\ldots,L\right\}} \left\{ \bm \xi^\top {\bm T_l}{\bm x^\prime} + \bm \xi^\top{\bm V_l}{\bm w^\prime} + \bm \xi^\top{\bm W}_l{\bm y^k}^\prime -{\bm H_l}{\bm \xi}\right\}
          \right\},
    \end{array}
\label{eq:separation2}
\end{equation}
where~$\Xi^K(\bm w^\prime)$ is defined as in~\eqref{eq:lift-unc}.

Then, using an epigraph to represent the minimization problem over~$k$ and assigning indicator variables for each objective of the inner maximization, we can write problem~\eqref{eq:separation2} as~\eqref{eq:separation-milp}.
\endproof

\begin{theorem}
Problem~\eqref{eq:scenario} admits an equivalent deterministic counterpart whose size grows linearly with the cardinality of~$\widehat{\Xi}$.
\label{thm:det-scenario}
\end{theorem} 

\proof
For each~$k \in \sets K$, consider any scenario~$\overline{\bm \xi} \in \widehat{\Xi}^k$. The maximization problem on the right-hand side of the first constraint of problem~\eqref{eq:scenario} can be written as
\begin{equation}
    \begin{array}{cl}
         \displaystyle\max & \quad {\bm \xi}^\top {\bm C} \; {\bm x} + {\bm \xi}^\top {\bm D} \; {\bm w}  + {\bm \xi}^\top {\bm Q} \; {\bm y^k} + \bm\psi^\top \bm \zeta \\
         \st & \quad {\bm \xi} \in \mathbb R^{N_{\xi}}\\
         & \quad {\bm A}{\bm \xi} \leq {\bm b}\\
         & \quad  {\bm w} \circ {\bm \xi}^k = {\bm w} \circ \overline{\bm \xi}\\
         & \quad \bm \zeta_s \leq -\bm g_{s,t}^\top \bm \xi^k \qquad \forall s \in \sets S, \, \forall t\in\sets T.
    \end{array}
    \label{eq:scenario-robust-obj}
\end{equation}
This linear program admits a strong dual given by the minimization problem
\begin{equation}
    \begin{array}{cl}
         \min & \quad {\bm b}^\top {\bm \beta(\overline{\bm \xi})} + \left({\bm w} \circ \overline{\bm \xi} \right)^\top {\bm \gamma}(\overline{\bm \xi})\\
         \st & \quad {\bm \beta}(\overline{\bm \xi}) \in \mathbb R^R_+, \; {\bm \gamma}(\overline{\bm \xi}) \in \mathbb R^{N_\xi}, \; { \delta}(\overline{\bm \xi})_{s,t} \in \mathbb R_+, \; \forall s \in \sets S, \; \forall t\in\sets T\\
        & \displaystyle \quad {\bm A}^\top{\bm \beta}(\overline{\bm \xi}) + \sum_{s\in \sets S}\sum_{t\in \sets T}{ \delta}(\overline{\bm \xi})_{s,t} \bm g_{s,t} + {\bm w} \circ \bm \gamma(\overline{\bm \xi}) = {\bm C} {\bm x} + {\bm D} {\bm w} + {\bm Q} {\bm y}^k \\
        & \displaystyle\quad \sum_{t\in\sets T}\bm \delta(\overline{\bm \xi})_t = \bm  \psi. \qquad
    \end{array}
    \label{eq:scenario-obj-dual}
\end{equation}
%By strong duality, problem~\eqref{eq:scenario-robust-obj} and~\eqref{eq:scenario-obj-dual} will have the same optimal objective value.

Similarly, in the second set of constraints, for any scenario~$\overline{\bm \xi} \in \Xi^k$, the maximization problem on the left-hand side of the~$l$th constraint can be written as
\begin{equation}
    \begin{array}{cl}
         \displaystyle\max & \quad \bm {\xi}^\top {\bm T_l}{\bm x} + \bm \xi^\top{\bm V_l}{\bm w} + \bm \xi^\top{\bm W}_l{\bm y^k} -{\bm H_l}{\bm \xi} \\
         \st & \quad {\bm \xi} \in \mathbb R^{N_{\xi}}\\
         & \quad {\bm A}{\bm \xi} \leq {\bm b}\\
         & \quad  {\bm w} \circ {\bm \xi} = {\bm w} \circ \overline{\bm \xi}.
    \end{array}
    \label{eq:scenario-robust-cstr}
\end{equation}
By taking dual of problem~\eqref{eq:scenario-robust-cstr}, we obtain a minimization problem whose optimal value coincides with that of the primal being expressible as
\begin{equation}
    \begin{array}{cl}
         \min & \quad {\bm b}^\top {\bm \alpha(\overline{\bm \xi})_l} + \left({\bm w} \circ \overline{\bm \xi} \right)^\top {\bm \eta}(\overline{\bm \xi})_l\\
         \st & \quad {\bm \alpha}(\overline{\bm \xi})_l \in \mathbb R^R_+, \; {\bm \eta}(\overline{\bm \xi})_l \in \mathbb R^{N_\xi} \\
        & \displaystyle \quad {\bm A}^\top{\bm \alpha}(\overline{\bm \xi})_l + {\bm w} \circ \bm \eta(\overline{\bm \xi})_l = {\bm T_l}{\bm x} + {\bm V_l}{\bm w} + {\bm W}_l{\bm y^k} -{\bm H}_l^{\top}.
    \end{array}
    \label{eq:scenario-cstr-dual}
\end{equation}
Therefore, replacing the constraints in \eqref{eq:scenario} with the respective  dual reformulations~\eqref{eq:scenario-obj-dual} and~\eqref{eq:scenario-cstr-dual}, we can write~\eqref{eq:scenario} as the deterministic problem
\begin{equation}
    \begin{array}{cl}
         \min & \quad \displaystyle \theta \\
         \st & \quad \theta\in\mathbb{R}, \; \bm\psi\in\mathbb{R}_+^{N_g}, \; {\bm x} \in \sets X, \; \bm y^k \in\mathcal{Y} \quad \forall k\in\sets K\\
         & \quad \!\! \left. \begin{array}{l} 
         
         {\bm \beta}(\overline{\bm \xi}) \in \mathbb R^R_+, \; {\bm \gamma}(\overline{\bm \xi}) \in \mathbb R^{N_\xi}, \; { \delta}(\overline{\bm \xi})_{s,t} \in \mathbb R_+, \; \forall s \in \sets S, \, \forall t\in\sets T \\
         
         {\bm \alpha}(\overline{\bm \xi})_l \in \mathbb R^R_+, \; {\bm \eta}(\overline{\bm \xi})_l \in \mathbb R^{N_\xi}, \; \forall l \in \sets L\\
         
         \theta \geq \bm c^\top \bm \psi + {\bm b}^\top {\bm \beta(\overline{\bm \xi})} + \left({\bm w} \circ \overline{\bm \xi} \right)^\top {\bm \gamma}(\overline{\bm \xi})  \\

         \begin{array}{r}
         {\bm A}^\top{\bm \beta}(\overline{\bm \xi}) + \sum_{s\in \sets S}\sum_{t\in \sets T}{ \delta}(\overline{\bm \xi})_{s,t} \bm g_{s,t} + {\bm w} \circ \bm \gamma(\overline{\bm \xi}) \\
         = {\bm C} {\bm x} + {\bm D} {\bm w} + {\bm Q} {\bm y}^k 
         \end{array}
         \\
         
        \displaystyle \sum_{t\in\sets T}\bm \delta(\overline{\bm \xi})_t = \bm  \psi \\
        
        \left. \begin{array}{l}
        
        {\bm b}^\top {\bm \alpha(\overline{\bm \xi})_l} + \left({\bm w} \circ \overline{\bm \xi} \right)^\top {\bm \eta}(\overline{\bm \xi})_l \leq 0\\
        
        \displaystyle{\bm A}^\top{\bm \alpha}(\overline{\bm \xi})_l + {\bm w} \circ \bm \eta(\overline{\bm \xi})_l = {\bm T_l}{\bm x} + {\bm V_l}{\bm w} + {\bm W}_l{\bm y^k} -{\bm H}_l^{\top}
        
        \end{array} \; \right\} \; \forall l \in \sets L
        
        \end{array} \; \right\} \; 
        
        \begin{array}{l}
             \forall \overline{\bm \xi} \in \widehat{\Xi}^k, \\
             \forall k \in \sets K,
        \end{array}
    \end{array}
\label{eq:scenario-det}
\end{equation}
where~$\bigcup_{k \in \sets K} \widehat{\Xi}^k = \widehat{\Xi}$. Note that whenever a new scenario~$\overline{\bm \xi}$ is added to the set~$\widehat{\Xi}$, we need to add a set of dual variables~$( {\bm \beta}(\overline{\bm \xi}), {\bm \gamma}(\overline{\bm \xi}), { \bm \delta}(\overline{\bm \xi}), \{{\bm \alpha}(\overline{\bm \xi})_l, {\bm \eta}(\overline{\bm \xi})_l\}_{l \in \sets L})$ into the deterministic problem~\eqref{eq:scenario-det}. \endproof

\begin{lemma}
If problem~\eqref{eq:Phi(w)} is feasible, then~$\bm \psi$ admits a finite optimal solution. 
\label{lemma: bounded-psi}
\end{lemma}

\proof
Consider the equivalent reformulation \eqref{eq:dual-full} of problem \eqref{eq:endo_whole} in the proof of Theorem~\ref{thm:min-max}. By Lemma \ref{lemma:interior} and \citet[Proposition 3.4]{Shapiro2001OnProblems}, for any feasible~$\bm x, \bm w, \bm y(\bm \xi)$, an optimal solution~$\bm \psi$ is attained by the dual problem. 
This claim also holds when we restrict $\bm y(\bm \xi)$ in  \eqref{eq:dual-full}  to feasible $K$-adaptable policies. The claim then follows since problem \eqref{eq:Phi(w)} constitutes an equivalent $K$-adaptable reformulation of~\eqref{eq:dual-full}. \endproof

\begin{lemma}
The cutting-plane algorithm for solving the scenario-based problem~\eqref{eq:scenario} converges in finitely many iterations. 
\label{lemma: cg-finit}
\end{lemma}

\proof
Fix a given quadruplet~$(\overline{\theta}, \overline{\bm \psi}, \overline{\bm x}, \{\overline{\bm y}^k\}_{k\in\mathcal{K}})$ in~\eqref{eq:scenario-violated}, for any~$\overline{\bm \xi} \in \widehat{\Xi}^k,\,\forall k \in \sets K$, the inner maximization problem of~$\bm \xi$ is a maximization problem of a convex function of~$\bm \xi$ over a  polyhedron~$\Xi(\bm w, \overline{\bm \xi})$. Thus,~\citet[Corollary 2.1]{Blankenship1976InfinitelyProblems} guarantees the finite convergence of the cutting-plane algorithm. \endproof

\proof[Proof of Theorem~\ref{thm:bnc-correct}]
We first establish the correctness of the algorithm after termination. If problem~\eqref{eq:Phi(w)} is infeasible, then in the branch-and-cut algorithm, there is no solution that can give a negative objective value in problem~\eqref{eq:separation} and be identified as feasible. Assume now that problem \eqref{eq:Phi(w)} is feasible and let~$(\theta^\star,\bm \psi^*,\bm x^*,\{{\bm {y}^k}^\star\}_{k\in\mathcal{K}})$ be an optimal solution with objective value~$\theta^\star$. The branch-and-cut tree must have a leaf node since it terminates. Consider an arbitrary leaf node in the branch-and-cut tree for which the optimal solution~$(\theta^\star,\bm \psi^*,\bm x^*,\{{\bm {y}^k}^\star\}_{k\in\mathcal{K}})$ is feasible. Such a leaf node exists because an optimal solution is feasible in the root node, and our branching mechanism ensures that the feasibility is maintained in at least \newqj{one} child node every time the algorithm branches. By definition, the leaf node was not branched. The node could have been fathomed in  step~1, where the current incumbent solution must also be optimal since it weakly dominates $(\theta^\star,\bm \psi^*,\bm x^*,\{{\bm {y}^k}^\star\}_{k\in\mathcal{K}})$. Otherwise, the node could have been fathomed in step 3, where the incumbent solution was updated to  $(\theta^\star,\bm \psi^*,\bm x^*,\{{\bm {y}^k}^\star\}_{k\in\mathcal{K}})$. 

For the asymptotic result, within any infinite branch, the finite convergence result derived in Lemma~\ref{lemma: cg-finit} ensures the existence of an infinite sequence of solutions in step 1 of Algorithm~\ref{alg:bnc}. For each solution in this sequence, the value $v_{k^\star}$ obtained by solving problem~\eqref{eq:scenario-violated} is less than or equal to 0.
We denote the sequence as~$\{(\theta^n,\bm \psi^n,\bm x^n,\{{\bm {y}^k}\}^n_{k\in\mathcal{K}})\}_1^{+\infty}$ and the associated solutions of the separation problem as~$\{\bm \xi^n\}_1^{+\infty}$,~where~$n$ is the order of node appear in the branch. Since~$\sets X$ and $\sets Y$ are compact sets, and an optimal solution~$\bm \psi$ is attained by Lemma~\ref{lemma: bounded-psi}, the Bolzano–Weierstrass theorem implies that the sequence~$\{(\theta^n,\bm \psi^n,\bm x^n,\{{\bm {y}}^k\}^n_{k\in\mathcal{K}}),\bm \xi^n\}_1^{+\infty}$ has at least one accumulation point, denoted as~$\{(\theta^\star,\bm \psi^*,\bm x^*,\{{\bm {y}^k}^\star\}_{k\in\mathcal{K}}),\bm \xi^\star\}$.

We first prove the feasibility of the accumulation point. Denote the value of problem~\eqref{eq:separation} under the solution~$(\theta,\bm \psi,\bm x,\{{\bm {y}}^k\}_{k\in\mathcal{K}})$ and scenario~$\bm \xi$ as~$z((\theta,\bm \psi,\bm x,\{{\bm {y}}^k\}_{k\in\mathcal{K}}), \bm \xi)$. Suppose that, given the accumulation solution~$(\theta^\star,\bm \psi^*,\bm x^*,\{{\bm {y}^k}^\star\}_{k\in\mathcal{K}})$, there exists a scenario~$\bm \xi^\prime$ such that~$z((\theta^\star,\bm \psi^*,\bm x^*,\{{\bm {y}^k}^\star\}_{k\in\mathcal{K}}),\bm \xi^\prime) > 0$. Since problem~\eqref{eq:separation} finds the worst case scenario  for the solution~$(\theta^\star,\bm \psi^*,\bm x^*,\{{\bm {y}^k}^\star\}_{k\in\mathcal{K}})$, the largest violation~$z((\theta^\star,\bm \psi^*,\bm x^*,\{{\bm {y}^k}^\star\}_{k\in\mathcal{K}}),\bm \xi^\star)\geq z((\theta^\star,\bm \psi^*,\bm x^*,\{{\bm {y}^k}^\star\}_{k\in\mathcal{K}}),\bm \xi^\prime) > 0$. The violation of the sequence~$\{(\theta^{n+1},\bm{\psi}^{n+1},\bm{x}^{n+1}, \{\bm y^k\}^{n+1}_{k\in\mathcal{K}}),\bm \xi^n)\}$ is always less than zero, and the sequence also converges to the same accumulation point. Thus, we have~$z((\theta^\star,\bm \psi^*,\bm x^*,\{{\bm {y}^k}^\star\}_{k\in\mathcal{K}}),\bm \xi^\star) \leq 0$ and raise the desired contradiction.

We now show the optimality of the accumulation point of the solution sequence corresponding to an infinite branch. If there is a solution of the problem~\eqref{eq:Phi(w)}, denoted as~$(\overline{\theta}, \overline{\bm \psi}, \overline{\bm x},\{\overline{\bm y}^k\}_{k\in\mathcal{K}})$, where~${\theta}^\prime \leq \theta^\star$. Then, using the same logic as the correctness proof, it is either a node solution of the branch-and-cut tree or an accumulation point of a sequence. In each case, the current objective value will be updated to~$\theta^\prime+\delta$, where~$\delta \geq 0$ and is an arbitrarily small number. Consequently, the branch corresponding to the sequence~$\{(\theta^n,\bm \psi^n,\bm x^n,\{{\bm {y}^k}\}^{n}_{k\in\mathcal{K}})\}_1^{+\infty}$ must be finite since it will be cut by a better solution. Therefore, we reach a contradiction with the infinite branch and conclude our proof. \endproof

\proof[Proof of Theorem \ref{thm:strengthen-feas}]

For any given $\bm w^r$, if a solution~${\bm w}$ violates the cut \eqref{eq:strengthen-feas}, then we have 
\begin{equation}
    \sum_{i \in \sets I_r} 1 - w_i \leq 0,
\end{equation}
and thus~$\bm w \geq \bm w^r$. If instead the solution $\bm w^r$ violates the constraint \eqref{eq: strengthen-feas-cstr} under the scenario~$\bar{\bm \xi}$, then for any solution ${\bm w}$ such that~$\bm w \geq \bm w^r$, we have
\begin{equation}
    \bm V(\bar{\bm \xi}) \bm w \geq \bm V(\bar{\bm \xi}) \bm w^r \geq \bm H \bar{\bm \xi},
\end{equation}
where the first inequality holds because the matrix $\bm V(\bm \xi)$ is component-wise non-negative, and the second one holds because of the infeasibility of~$\bar{\bm \xi}$. Thus~$\bm w$ is infeasible in problem \eqref{eq:endo_whole}. This concludes the proof.\endproof

\section{Supplementary Numerical Results}
\label{ec:results}
\setcounter{table}{0}
\renewcommand{\thetable}{B.\arabic{table}}
\setcounter{figure}{0} 
\renewcommand{\thefigure}{\Alph{section}.\arabic{figure}}

\subsection{Supplementary Numerical Results to Section \ref{subsubsec: b&b-res}}
\label{ecsub:bb}

\subsubsection{Supplementary Runtime Results}
\label{ecsubsub:bb-time}

In Section \ref{subsubsec: b&b-res}, we state that since the number of bilinear terms between continuous variables is small, solving the monolithic MINLO is more efficient than using the decomposition algorithm in all cases. From Table \ref{table:bb-per}, we know that Gurobi solves all instances to optimality within the given time limit. In this section, for computational completeness, we compare and report the performance of the decomposition Algorithm~\ref{alg:l-shaped}, which uses the Branch-and-Cut Algorithm \ref{alg:bnc} to evaluate~$\Phi(\bm w)$, both with and without the cuts introduced in Section \ref{subsec:tighter-cuts}. We also evaluate its performance when the Branch-and-Bound Algorithm \ref{alg:bnb} is used with the same cuts in Section \ref{subsec:tighter-cuts}. In Algorithm \ref{alg:l-shaped}, the time limit for evaluating~$\Phi(\bm w)$ is set to~$300$ seconds, and the lower bound~$L$ is obtained by solving a deterministic problem with~$\bm \xi = (\bm r, \bm c)$.

Table~\ref{table:bb-sup} summarizes the computational results for the same randomly generated instances. From the table, we observe that the improved decomposition algorithm, which uses Algorithm \ref{alg:bnb} to evaluate~$\Phi(\bm w)$, outperforms the basic decomposition algorithm. Across all instances, the proposed algorithm solves more instances to optimality and achieves a smaller optimality gap, which decreases as~$K$ increases. As the number of policies~$K$ increases, the improved decomposition algorithm outperforms the other two approaches in an increasing number of instances. For example, when~$N=40$, the number of times a better solution was found by the improved decomposition approach is~$\{12, \, 19, \, 20\}$ for~$K \in \{2, \,3, \,4\}$, respectively. Additionally, the improvements in objective value over the static solutions achieved by the improved decomposition algorithm are greater than those of the other two methods across all cases. Notably, the decomposition algorithm using Algorithm \ref{alg:bnb} to evaluate~$\Phi(\bm w)$ rarely finds high-quality solutions when the problem size is large (e.g.,\,~$N > 10$), and its solution time is longer than the other two methods for small-sized problems. This result highlights the efficiency of the proposed Branch-and-Cut Algorithm \ref{alg:bnc} and its advantages.

\begin{landscape}
\setlength{\tabcolsep}{2pt}

\begin{table}[!ht]
\renewcommand{\arraystretch}{1.5}
\begin{footnotesize}
\begin{tabular}{cccccccccccccccccccc } 
 \hline
 & & & \multicolumn{5}{c}{ Decomposition} & & \multicolumn{5}{c}{ Decomposition with Added Cuts} & & \multicolumn{5}{c}{Decomposition with Algorithm \ref{alg:bnb} and Cuts}\\
 $N$ & $K$ & $\,$ & Opt($\#$) & Time(s) & Gap & Better($\#$) & Improvement & $\,$ &  Opt($\#$) & Time(s) & Gap & Better($\#$) & Improvement & $\,$ &  Opt($\#$) & Time(s) & Gap & Better($\#$) & Improvement\\ 
 \hline
\multirow{3}{*}{10} & 2 && $ \textbf{20} $ & $ 20.7 $ & $ \textbf{0.0\%} $ & $ \textbf{0} $ & $ \textbf{42.6\%} $ && $ \textbf{20} $ & $ \textbf{4.2} $ & $ \textbf{0.0\%} $ & $ \textbf{0} $ & $ \textbf{42.6\%} $ && $ \textbf{20} $ & $ 14.2 $ & $ \textbf{0.0\%} $ & $ \textbf{0} $ & $ \textbf{42.6\%} $ \\
~ & 3 && $ \textbf{20} $ & $ 31.1 $ & $ \textbf{0.0\%} $ & $ \textbf{0} $ & $ \textbf{43.9\%} $ && $ \textbf{20} $ & $ \textbf{7.2} $ & $ \textbf{0.0\%} $ & $ \textbf{0} $ & $ \textbf{43.9\%} $ && $ \textbf{20} $ & $ 50.9 $ & $ \textbf{0.0\%} $ & $ \textbf{0} $ & $ \textbf{43.9\%} $ \\
~ & 4 && $ \textbf{20} $ & $ 52.8 $ & $ \textbf{0.0\%} $ & $ \textbf{0} $ & $ \textbf{43.9\%} $ && $ \textbf{20} $ & $ \textbf{26.6} $ & $ \textbf{0.0\%} $ & $ \textbf{0} $ & $ \textbf{43.9\%} $ && $ \textbf{20} $ & $ 101.8 $ & $ \textbf{0.0\%} $ & $ \textbf{0} $ & $ \textbf{43.9\%} $ \\
\hline
\multirow{3}{*}{20} & 2 && $ 0 $ & $ 7200.0 $ & $ 40.7\% $ & $ 1 $ & $ 84.6\% $ && $ \textbf{13} $ & $ \textbf{2769.5} $ & $ \textbf{23.6\%} $ & $ 1 $ & $ \textbf{85.9\%} $ && $ 0 $ & $ 7200.0 $ & $ 100.0\% $ & $ 0 $ & $ -100.0\% $ \\
~ & 3 && $ 0 $ & $ 7200.0 $ & $ 32.4\% $ & $ 4 $ & $ 110.6\% $ && $ \textbf{5} $ & $ \textbf{1071.4} $ & $ \textbf{19.2\%} $ & $ \textbf{8} $ & $ \textbf{108.1\%} $ && $ 0 $ & $ 7200.0 $ & $ 100.0\% $ & $ 0 $ & $ -100.0\% $ \\
~ & 4 && $ 0 $ & $ 7200.0 $ & $ 30.5\% $ & $ 3 $ & $ 116.6\% $ && $ \textbf{8} $ & $ \textbf{3494.4} $ & $ \textbf{15.7\%} $ & $ \textbf{13} $ & $ \textbf{120.3\%} $ && $ 0 $ & $ 7200.0 $ & $ 100.0\% $ & $ 0 $ & $ -100.0\% $ \\
\hline
\multirow{3}{*}{30} & 2 && $ 0 $ & $ 7200.0 $ & $ 41.4\% $ & $ 0 $ & $ 111.9\% $ && $ \textbf{2} $ & $ \textbf{1799.4} $ & $ \textbf{34.1\%} $ & $ \textbf{2} $ & $ \textbf{118.2\%} $ && $ 0 $ & $ 7200.0 $ & $ 100.0\% $ & $ 0 $ & $ -100.0\% $ \\
~ & 3 && $ \textbf{0} $ & $ \textbf{7200.0} $ & $ 36.2\% $ & $ 2 $ & $ 131.2\% $ && $ \textbf{0} $ & $ \textbf{7200.0} $ & $ \textbf{28.9\%} $ & $ \textbf{13} $ & $ \textbf{147.7\%} $ && $ \textbf{0} $ & $ \textbf{7200.0} $ & $ 100.0\% $ & $ \textbf{0} $ & $ -100.0\% $ \\
~ & 4 && $ 0 $ & $ 7200.0 $ & $ 37.0\% $ & $ 1 $ & $ 127.9\% $ && $ 1 $ & $ \textbf{1182.8} $ & $ \textbf{25.5\%} $ & $ \textbf{19} $ & $ \textbf{158.3\%} $ && $ 0 $ & $ 7200.0 $ & $ 100.0\% $ & $ 0 $ & $ -100.0\% $ \\
\hline
\multirow{3}{*}{40} & 2 && $ 0 $ & $ 7200.0 $ & $ 48.3\% $ & $ 2 $ & $ 86.4\% $ && $ \textbf{1} $ & $ \textbf{3356.4} $ & $ \textbf{44.9\%} $ & $ \textbf{12} $ & $ \textbf{98.0\%} $ && $ 0 $ & $ 7200.0 $ & $ 100.0\% $ & $ 0 $ & $ -100.0\% $ \\
~ & 3 && $ \textbf{0} $ & $ \textbf{7200.0} $ & $ 44.2\% $ & $ 0 $ & $ 101.0\% $ && $ \textbf{0} $ & $ \textbf{7200.0} $ & $ \textbf{35.8\%} $ & $ \textbf{19} $ & $ \textbf{129.6\%} $ && $ \textbf{0} $ & $ \textbf{7200.0} $ & $ 100.0\% $ & $ 0 $ & $ -100.0\% $ \\
~ & 4 && $ \textbf{0} $ & $ \textbf{7200.0} $ & $ 42.6\% $ & $ 0 $ & $ 106.8\% $ && $ \textbf{0} $ & $ \textbf{7200.0} $ & $ \textbf{30.1\%} $ & $ \textbf{20} $ & $ \textbf{149.3\%} $ && $ \textbf{0} $ & $ \textbf{7200.0} $ & $ 100.0\% $ & $ 0 $ & $ -100.0\% $ \\
\hline
\end{tabular}

\caption{Supplementary computational results of the best box problem. $\mathbf{Decomposition}$ and $\mathbf{ Decomposition\;\;with\;\;Added\;\;Cuts}$ 
refer to Algorithm~\ref{alg:l-shaped} using the Branch-and-Cut Algorithm \ref{alg:bnc} to evaluate~$\Phi(\bm w)$, without and with the cuts introduced in Section~\ref{subsec:tighter-cuts}, respectively. $\mathbf{Decomposition \;\; with \;\; Algorithm \;\; \ref{alg:bnb} \;\; and \;\; Cuts}$ refers to Algorithm~\ref{alg:l-shaped} using the Branch-and-Bound Algorithm \ref{alg:bnb} with the cuts introduced in Section~\ref{subsec:tighter-cuts}. $\mathbf{Opt(\#)}$ corresponds to the number of instances solved to optimality, $\mathbf{Time(s)}$ to the average computational time (in seconds) for instances solved to optimality, and $\mathbf{Gap}$ to the average optimality gap for the instances not solved within the time limit. $\mathbf{Better(\#)}$ denotes the number of instances in each method that achieved a better objective value than the other method. $\mathbf{Improvement}$ denotes the average improvement in the objective value of the $K$-adaptability solution found in the time limit over the static solution found in the time limit. The improvement is calculated as the ratio of the difference between the objective values of the $K$-adaptability solution and the static solution to the objective value of the static solution.}
\label{table:bb-sup}
\end{footnotesize}
\end{table}
\end{landscape}

\newqj{
\subsubsection{Improvement over~$K$}
\label{ecsubsub:bb-improve-k}

In this section, to study the quality of the solutions obtained by the~$K$-adaptability approximation, we randomly generated 20 instances of a small-sized best box problem with~$N = 10$,  which we can solve to optimality. Based on the structure of the problem, the original problem \eqref{eq:nested} with full adaptability needs at most 11 candidate policies (one for selecting each box at the second stage, and another one if no box is opened). We now summarize the average improvements over the static solution as~$K$ increases in Figure~\ref{fig:bb-k}. 

\begin{figure}[htbp!]
\centering
\begin{tikzpicture}
\begin{axis}[
    width=12cm,
    height=7cm,
    xmin=2, xmax=11,
    ymin=42, ymax=45,
    xtick={2,...,11},
    ytick={42, 43, 44, 45},
    xlabel={Number of $K$},
    ylabel={Improvement (\%)},
    grid=major,
    grid style={line width=.2pt, draw=gray!50},
    axis lines=box,                    % Adds frame around plot
    tick label style={font=\small}
]

% Data points
\addplot[
    mark=*,
    thick,
    color=black!80
]
coordinates {
    (2, 42.6)
    (3, 43.9)
    (4, 43.9)
    (5, 43.9)
    (6, 43.9)
    (7, 43.9)
    (8, 43.9)
    (9, 43.9)
    (10, 43.9)
    (11, 43.9)
};

% Dashed vertical lines at K=2 and K=3
% \addplot[dashed, thick, gray] coordinates {(3,42) (3,43.9)};

% Dashed horizontal lines
% \addplot[dashed, thick, gray] coordinates {(2,43.9) (3,43.9)};

\end{axis}
\end{tikzpicture}
\caption{\newqj{Average improvement over the static solution in the best box problem as~$K$ increases}}

\label{fig:bb-k}
\end{figure}

The average improvements over the static solution for~$K = 2$ and~$K \geq 3$ are~$42.6\%$ and~$43.9\%$, respectively. This result shows that, in the small-sized best box problem, $3$ (out of a total of~$11$) candidate policies are sufficient for the~$K$-adaptability approximation \eqref{eq:k-nested} to be optimal in the original problem~\eqref{eq:nested}.}

\subsubsection{Out of Sample Performance}

\label{ecsubsub:bb-out-sample}

In this section, we compare the out-of-sample performance of the RO solution and the DRO solution in the best box problem. For all the problems, we sample 200 distributions and calculate the expected costs under both the RO and DRO solutions. For each distribution, we draw 1000 samples of the random parameter~$\bar{\bm \zeta}$ from the uniform distribution~$U[-1, 1]$, and generate the scenarios~$\bar{\bm \xi}$ using the relationship between them defined in the support set:
\begin{equation*}
    \bar{\bm \xi}=\left(\bm \xi^{\mathrm{r}}, \bm{\xi}^{\mathrm{c}}\right)\text{, where }\bar{\xi}_{n}^{\mathrm{r}}=\left(1+\mathbf{\Phi}_{n}^{\top} \bar{\bm{\zeta}} / 2\right) r_{n}, \; \bar{\xi}_{n}^{\mathrm{c}}=\left(1+\bm{\Psi}_{n}^{\top} \bar{\bm{\zeta}} / 2\right) c_{n},\; n \in \sets N.
\end{equation*}

If the constraints in the ambiguity set are not satisfied, we decrease the interval length of the uniform distribution by~$2/3$ and resample the scenarios of~$\bm \zeta$ until the constraints are satisfied. For each sample~$\bar{\bm \xi}$ and DRO solution, denoted as~$\{\bm w, \bm x, \bm \psi, \{\bm y^k\}_{k\in\sets K}\}$, we solve the following inner maximization problem to find the best candidate policy~$\bm y^{k^\star}$: 
\begin{equation}
 \begin{array}{cl}
            \displaystyle\min_{ k \in \sets K } \;\;
             &\displaystyle\max_{ {\bm \xi} \in \Xi({\bm w},\overline{\bm \xi}) } \; \bm c^\top \bm \psi + {\bm \xi}^\top {\bm C} \; {\bm x} + {\bm \xi}^\top {\bm D} \; {\bm w}  + {\bm \xi}^\top {\bm Q} \; {\bm y^k} - \bm\psi^\top \bm g(\bm\xi) \\
             \st & \;\; {\bm T} ({\bm \xi}){\bm x} + {\bm V} ({\bm \xi}){\bm w} + {\bm W} ({\bm \xi}){\bm y^k} \leq {\bm H}{\bm \xi} \;\;\; \forall {\bm \xi} \in \Xi({\bm w},\overline{\bm \xi}).
         \end{array}
\label{eq:BB-DRO-inner}
\end{equation}
The objective value associated with~$\bar{\bm \xi}$ is calculated as~$\bar{\bm \xi}^\top {\bm C} {\bm x} + \bar{\bm \xi}^\top {\bm D} {\bm w} + \bar{\bm \xi}^\top {\bm Q}{\bm y^{k^\star}}$. We then use the objective value of each sample to calculate the expectation of the discrete distribution.

For the RO solution~$\{\bm w, \bm x, \{\bm y^k\}_{k\in\sets K}\}$, we solve the following inner maximization problem to find the best candidate policy:
\begin{equation}
 \begin{array}{cl}
            \displaystyle\min_{ k \in \sets K } \;\;
             &\displaystyle\max_{ {\bm \xi} \in \Xi({\bm w},\overline{\bm \xi}) } \; {\bm \xi}^\top {\bm C} \; {\bm x} + {\bm \xi}^\top {\bm D} \; {\bm w}  + {\bm \xi}^\top {\bm Q} \; {\bm y^k} \\
             \st & \;\; {\bm T} ({\bm \xi}){\bm x} + {\bm V} ({\bm \xi}){\bm w} + {\bm W} ({\bm \xi}){\bm y^k} \leq {\bm H}{\bm \xi} \;\;\; \forall {\bm \xi} \in \Xi({\bm w},\overline{\bm \xi}),
         \end{array}
\label{eq:BB-RO-inner}
\end{equation}
and calculate the expected value in the same way described above.
 
Table \ref{table:bb-out-sample} reports the relative difference of the median, the~$95\%$ and~$99\%$ quantiles between the two methods. From the simulation results, we can see that the DRO approach, in general, produces less conservative solutions that perform better than those of the RO approach even on non-worst-case distributions. 

\begin{table}[!ht]
\renewcommand{\arraystretch}{1.5}

\centering
\begin{tabular}{ ccccc } 
 \hline
$N$ & $K$ & median & $95\%$ quantile & $99\%$ quantile \\ 
 \hline
\multirow{4}{*}{10} & $ 1 $ & $ 7.10\% $ & $ 7.20\% $ & $ 7.32\% $ \\
~ & $ 2 $ & $ 2.44\% $ & $ 2.34\% $ & $ 2.31\% $ \\
~ & $ 3 $ & $ 4.04\% $ & $ 3.91\% $ & $ 3.90\% $ \\
~ & $ 4 $ & $ 4.03\% $ & $ 3.90\% $ & $ 3.89\% $ \\
\hline
\multirow{4}{*}{20}  & $ 1 $ & $ 6.84\% $ & $ 6.80\% $ & $ 6.83\% $ \\
~ & $ 2 $ & $ 0.23\% $ & $ 0.22\% $ & $ 0.23\% $ \\
~ & $ 3 $ & $ 0.96\% $ & $ 0.96\% $ & $ 0.97\% $ \\
~ & $ 4 $ & $ 1.06\% $ & $ 1.05\% $ & $ 1.07\% $ \\
\hline
\multirow{4}{*}{30} & $ 1 $ & $ 2.77\% $ & $ 2.75\% $ & $ 2.73\% $ \\
~ & $ 2 $ & $ 1.69\% $ & $ 1.69\% $ & $ 1.68\% $ \\
~ & $ 3 $ & $ 1.64\% $ & $ 1.65\% $ & $ 1.66\% $ \\
~ & $ 4 $ & $ 1.87\% $ & $ 1.87\% $ & $ 1.88\% $ \\
\hline
\multirow{4}{*}{40} & $ 1 $ & $ 3.68\% $ & $ 3.65\% $ & $ 3.67\% $ \\
~ & $ 2 $ & $ 0.58\% $ & $ 0.57\% $ & $ 0.57\% $ \\
~ & $ 3 $ & $ -0.33\% $ & $ -0.31\% $ & $ -0.31\% $ \\
~ & $ 4 $ & $ 0.43\% $ & $ 0.44\% $ & $ 0.44\% $ \\
 \hline
\end{tabular}

\caption{The relative difference between the median, the~$95\%$ and~$99\%$ quantiles of the expected objective values in the 200 simulated distributions given by the RO and DRO solutions to the Best Box problem. Each cell represents an average of 20 instances.}

\label{table:bb-out-sample}
\end{table}

\subsection{Supplementary Numerical Results for Section \ref{subsubsec: r&d-result}}
\label{ecsub:r&d}

\subsubsection{Supplementary Runtime Results}
\label{ecsubsub:r&d-time}

In this section, for computational completeness, we present the performance of the decomposition Algorithm \ref{alg:l-shaped} using the Branch-and-Bound Algorithm \ref{alg:bnb} with the cuts introduced in Section \ref{subsec:tighter-cuts}. In Algorithm~\ref{alg:l-shaped}, the time limits for evaluating~$\Phi(\bm w)$ are set to~$[300s, 600s, 1200s, 2400s]$ for problem sizes~$N \in \{5, \; 10, \; 15, \; 20\}$, consistent with Section \ref{subsubsec: r&d-result}. The lower bound~$L$ is again obtained by solving a deterministic problem with~$\bm \xi = (\bm r, \bm c)$.

Table~\ref{table:r&d-sup} summarizes the computational results for the same randomly generated instances and should be read in tandem with Table~\ref{table:r&d-per}. From the table, we observe that, across all cases, using Algorithm \ref{alg:bnc} to evaluate~$\Phi(\bm w)$ solves more instances to optimality and, if the time limit is reached, it returns solutions with better objective values. Algorithm \ref{alg:bnc} also requires shorter solving times, has smaller gaps, and shows more improvement. Particularly, when the problem size is large (e.g.,\,$N > 10$), the decomposition algorithm using Algorithm \ref{alg:bnc} finds solutions having considerable improvement over the static solutions, while Algorithm \ref{alg:bnb} has difficulty finding high-quality solutions. This again demonstrates the efficiency and advantages of the proposed Branch-and-Cut Algorithm \ref{alg:bnc}.

\setlength{\tabcolsep}{2pt}

\begin{table}[!htp]
\renewcommand{\arraystretch}{1.5}
\begin{footnotesize}
\centering
\begin{tabular}{cccccccccc } 
 \hline
 $N$ & $K$ & $\,$ & Opt($\#$) & Time(s) & Gap & Better($\#$) & Improvement \\ 
 \hline
\multirow{3}{*}{5} & 2 && $ 6 $ & $ 253.0 $ & $ 17.9\% $ & $ 0 $ & $ 10.1\% $\\
~ & 3 && $ 0 $ & $ 7200.0 $ & $ 17.8\% $ & $ 0 $ & $ 12.4\% $\\
~ & 4 && $ 0 $ & $ 7200.0 $ & $ 18.4\% $ & $ 0 $ & $ 12.4\% $\\
\hline
\multirow{3}{*}{10} & 2 && $ 0 $ & $ 7200.0 $ & $ 34.1\% $ & $ 0 $ & $ 4.4\% $\\
~ & 3 && $ 0 $ & $ 7200.0 $ & $ 28.1\% $ & $ 1 $ & $ 14.2\% $\\
~ & 4 && $ 0 $ & $ 7200.0 $ & $ 31.9\% $ & $ 0 $ & $ 5.6\% $\\
\hline
\multirow{3}{*}{15} & 2 && $ 0 $ & $ 7200.0 $ & $ 88.9\% $ & $ 0 $ & $ -73.0\% $\\
~ & 3 && $ 0 $ & $ 7200.0 $ & $ 91.1\% $ & $ 0 $ & $ -73.0\% $\\
~ & 4 && $ 0 $ & $ 7200.0 $ & $ 88.8\% $ & $ 0 $ & $ -73.0\% $\\
\hline
\multirow{3}{*}{20} & 2 && $ 0 $ & $ 7200.0 $ & $ 95.9\% $ & $ 0 $ & $ -79.8\% $\\
~ & 3 && $ 0 $ & $ 7200.0 $ & $ 94.4\% $ & $ 0 $ & $ -79.8\% $\\
~ & 4 && $ 0 $ & $ 7200.0 $ & $ 94.9\% $ & $ 0 $ & $ -79.8\% $\\
\hline
\end{tabular}

\caption{Computational results of Algorithm~\ref{alg:l-shaped} using the Branch-and-Bound Algorithm \ref{alg:bnb} with the cuts introduced in Section~\ref{subsec:tighter-cuts} for R$\&$D problem. $\mathbf{Opt(\#)}$ corresponds to the number of instances solved to optimality, $\mathbf{Time(s)}$ to the average computational time (in seconds) for instances solved to optimality, and $\mathbf{Gap}$ to the average optimality gap for the instances not solved within the time limit. $\mathbf{Better(\#)}$ denotes the number of instances in each method that achieved a better objective value than the other method. $\mathbf{Improvement}$ denotes the average improvement in the objective value of the $K$-adaptability solution found in the time limit over the static solution found in the time limit. The improvement is calculated as the ratio of the difference between the objective values of the $K$-adaptability solution and the static solution to the objective value of the static solution.}
\label{table:r&d-sup}
\end{footnotesize}
\end{table}

\newqj{
\subsubsection{Improvement over~$K$}
\label{ecsubsub:r&d-improve-k}

In this section, to study the quality of the solutions obtained by the~$K$-adaptability approximation, we consider a small-sized problem with~$N = 4$ and~$M = 2$, where we are able to solve the problem with the set of candidate policies containing all possible second-stage decisions. We randomly generate 20 instances, and summarize the average improvements over the static solution as~$K$ increases in Figure~\ref{fig:r&d-k}. 

\begin{figure}[htbp!]
\centering
\begin{tikzpicture}
\begin{axis}[
    width=12cm,
    height=7cm,
    xmin=2, xmax=16,
    ymin=34, ymax=43,
    xtick={2,4,6,8,10,12,14,16},
    ytick={35, 37.5, 40, 42.5},
    xlabel={Number of $K$},
    ylabel={Improvement (\%)},
    grid=major,
    grid style={line width=.2pt, draw=gray!50},
    axis lines=box,                    % Adds frame around plot
    tick label style={font=\small}
]

% Data points
\addplot[
    mark=*,
    thick,
    color=black!80
]
coordinates {
    (2, 35.0)
    (3, 41.82)
    (4, 41.92)
    (5, 42.06)
    (6, 42.06)
    (7, 42.06)
    (8, 42.06)
    (9, 42.06)
    (10, 42.06)
    (11, 42.06)
    (12, 42.06)
    (13, 42.06)
    (14, 42.06)
    (15, 42.06)
    (16, 42.06)
};
% Dashed vertical lines at K=2 and K=3
% \addplot[dashed, thick, gray] coordinates {(3,42) (3,43.9)};

% Dashed horizontal lines
% \addplot[dashed, thick, gray] coordinates {(2,42.1) (16,42.1)};

\end{axis}
\end{tikzpicture}
\caption{\newqj{Average improvement over static solution in the R\&D project portfolio
optimization problem as~$K$ increases}}

\label{fig:r&d-k}
\end{figure}

The average improvement over the static solutions for~$K=2,\,3,\ldots,\,6$ are~$35.0\%,\,41.8\%,\,41.9\%,\,42.1\%,\,42.1\%$, respectively. This result shows that, in the small-sized~R$\&$D problem, $5$ (out of a total of~$16$) candidate policies are sufficient for the~$K$-adaptability approximation \eqref{eq:k-nested} to be optimal compared to the original problem \eqref{eq:nested}.}

\subsubsection{Out of Sample Performance}
\label{ecsubsub:r&d-out-sample}

In this section, we compare the out-of-sample performance of the RO solution and the DRO solution in the R\&D problem. We sample 200 distributions that satisfy the moment constraints in~$\mathcal{P}_\text{RD}$. For each instance and distribution, we evaluate the expected return~$\mathbb E_{\mathbb P}[-{\bm \xi^r}^\top (\bm w + \theta \bm y(\bm \xi))]$ using~$\bm y(\bm \xi)$ as the solution computed by RO and DRO approaches. The detail\newqj{s} of the sampling procedure are similar to Section \ref{ecsubsub:bb-out-sample}. Table \ref{table:r&d-out-sample} reports the relative difference of the median, the~$95\%$ and~$99\%$ quantiles between the two methods. Similar to the best box problem, Table \ref{table:r&d-out-sample} shows that, in the R$\&$D problem, the DRO solutions, in general, also outperform the RO solutions on non-worst-case distributions. 

\begin{table}[!htp]
\renewcommand{\arraystretch}{1.5}
% \begin{center}
\begin{tabular}{ cccccc } 
 \hline
$N$ & $K$ & median & $95\%$ quantile & $99\%$ quantile \\ 
 \hline
\multirow{4}{*}{5} & $ 1 $ & $ -0.03\% $ & $ -0.03\% $ & $ -0.03\% $ \\
~ & $ 2 $ & $ 9.55\% $ & $ 9.78\% $ & $ 9.87\% $ \\
~ & $ 3 $ & $ 4.53\% $ & $ 4.68\% $ & $ 4.71\% $ \\
~ & $ 4 $ & $ 2.90\% $ & $ 3.03\% $ & $ 3.06\% $ \\
\hline
\multirow{4}{*}{10} & $ 1 $ & $ -0.27\% $ & $ -0.26\% $ & $ -0.25\% $ \\
~ & $ 2 $ & $ 22.92\% $ & $ 22.95\% $ & $ 22.96\% $ \\
~ & $ 3 $ & $ 4.50\% $ & $ 4.80\% $ & $ 4.92\% $ \\
~ & $ 4 $ & $ 4.46\% $ & $ 4.54\% $ & $ 4.56\% $ \\
\hline
\multirow{4}{*}{15} & $ 1 $ & $ -0.60\% $ & $ -0.59\% $ & $ -0.60\% $ \\
~ & $ 2 $ & $ 21.66\% $ & $ 21.67\% $ & $ 21.66\% $ \\
~ & $ 3 $ & $ 3.07\% $ & $ 3.51\% $ & $ 3.60\% $ \\
~ & $ 4 $ & $ 3.93\% $ & $ 3.95\% $ & $ 4.03\% $ \\
\hline
\multirow{4}{*}{20} & $ 1 $ & $ -0.26\% $ & $ -0.25\% $ & $ -0.25\% $ \\
~ & $ 2 $ & $ 21.44\% $ & $ 21.44\% $ & $ 21.43\% $ \\
~ & $ 3 $ & $ 4.85\% $ & $ 5.32\% $ & $ 5.47\% $ \\
~ & $ 4 $ & $ 4.20\% $ & $ 4.21\% $ & $ 4.25\% $ \\
 \hline
\end{tabular}
% \end{center}
\caption{The relative difference between the median, the~$95\%$ and~$99\%$ quantiles of the expected objective values in the 200 simulated distributions given by the RO and DRO solutions to the R$\&$D problem. Each cell represents an average of 20 instances.}

\label{table:r&d-out-sample}
\end{table}

\subsection{Distributionally Robust R$\&$D Project Portfolio Optimization with Loans}
\label{ecsec:r&d-loan}
\subsubsection{Problem Description}
\label{ecsubsec: r&d-loans-problem}
We perform experiments on several instances of a variant of the distributionally robust R\&D project portfolio optimization problem studied in Section \ref{subsec:r&d}. In this problem, in the first stage, the decision\newqj{-}maker is allowed to borrow an amount~$x$ of money which she must pay back, with interest~$\alpha$ in the final stage. We use the same ambiguity set as the one described in \ref{subsec:r&d}. The support of the distribution of the uncertain parameters is given by:
\begin{equation*}
\begin{array}[t]{r}
    \Xi=\left\{\left(\bm \xi^{\mathrm{r}}, \bm{\xi}^{\mathrm{c}}\right) \in \mathbb{R}^{2 N}: \exists \bm{\zeta} \in[-1,1]^{M}: {\xi}_{n}^{\mathrm{r}}=\left(1+\mathbf{\Phi}_{n}^{\top} \bm{\zeta} / 2\right){r}_{n}, \;{\xi}_{n}^{\mathrm{c}}=\left(1+\bm{\Psi}_{n}^{\top} \bm{\zeta} / 2\right) {c}_{n}, \right .\\ \left . n\in \sets N\right\}.
\end{array}
\end{equation*}
We define $\bm \xi = \left(\bm \xi^\mathrm{r},\bm\xi^\mathrm{c}\right)$,~$\bm \mu = \left(\bm r,\bm c\right)$, and use the following ambiguity set for the uncertain returns and costs:
\begin{equation}
\begin{array}{cl}
    \mathcal{P}_\text{RD} = \Big\{\mathbb{P}\in\mathcal{M}_+(\mathbb{R}^{N_\xi})\; :\; \mathbb{P}(\bm \xi\in\Xi) = 1,\;\;&\mathbb{E}_{\mathbb{P}}\left[\left|\bm \xi - \bm \mu \right| \right ]\leq 0.15 \bm \mu,\\ &\mathbb{E}_{\mathbb{P}}\left[\left|\mathbf e^\top \left(\bm \xi - \bm \mu \right )\right| \right ]\leq 0.15 N^{-1/2} \mathbf e^\top \bm \mu)\Big\}.
\end{array}
\label{eq: r&d-borrow-ambi_set}
\end{equation} 
The second set of constraints in the ambiguity set enforces the mean absolute deviation of the individual returns and costs to be bounded. The last constraint imposes an upper bound on the cumulative deviation of the projects' returns and costs from their expected value.

The distributionally robust R\&D project portfolio optimization problem with loans can be written as
\begin{equation}
    \begin{array}{cl}
         \min & \quad \displaystyle \sup_{\mathbb{P} \in \mathcal{P}} \;\;\;\mathbb{E}_{\mathbb{P}}\left[  -{\bm \xi^r}^\top \left(\bm w + \theta {\bm y}({\bm \xi}) \right) + (1 + \alpha) x \right] \\
         \st & \quad {\bm w} = (\bm w^r,\; \bm w^r),\; \bm w^r \in \left\{0,\;1\right\}^N,\; x \in \reals_+ \\
         %& \quad {\bm R}{\bm x} + {\bm S}{\bm w} \; \leq \; {\bm t} \\
         & \quad \!\! \left. \begin{array}{l} 
         {\bm y}({\bm \xi}) \in \left\{0,\;1\right\}^N  \\
         {\bm \xi^c}^\top \left(\bm w^r + {\bm y}({\bm \xi}) \right) \leq B + x\\
         \bm w^r + \bm y(\bm \xi) \leq \mathbf e
         \end{array} \quad \right\} \quad \forall {\bm \xi} \in \Xi \\
         & \quad {\bm y}({\bm \xi}) = {\bm y}({\bm \xi}') \quad \forall {\bm \xi}, {\bm \xi}' \in \Xi \; : \; {\bm w} \circ {\bm \xi} = {\bm w} \circ {\bm \xi}'.
    \end{array}
\label{eq:r&d-loans}
\end{equation}

\subsubsection{Computational Results}
\label{ecsubsec: r&d-loans-result}

We run experiments on~$80$ randomly generated instances of the~$K$-adaptability counterpart of problem~\eqref{eq:r&d}:~$20$ instances for each~$N\in\{5,10,15,20\}$~with risk factors~$M \in \{4,5,8,10\}$, respectively. We draw~$\bm c$ uniformly at random from the box~$[0,10]^N$ and let~$\bm r = 5 \bm c,\;B = \mathbf e^\top \bm c/2$ and~$\alpha = 0.2$. The risk load factors~$\bm \Phi_n$ and~$\bm \Psi_n$ are sampled uniformly at random from the standard simplex. 

We report the results of solving problem \eqref{eq:r&d} by utilizing the monolithic MINLO reformulation and applying the decomposition Algorithm~\ref{alg:l-shaped} with and without the cuts introduced in Section~\ref{subsec:tighter-cuts}. In the decomposition Algorithm~\ref{alg:l-shaped}, we set the time limit of evaluating~$\Phi(\bm w)$ to~$[300s, 600s, 1200s, 2400s]$ for problems with size~$N\in\{5,\;10,\;15,\;20\}$. The lower bound value~$L$ in Algorithm \ref{alg:l-shaped} is found by solving a deterministic problem obtained by setting~$\bm \xi = (\bm r, \bm c)$. If the subproblem evaluating~$\Phi(\bm w)$ is not solved to optimality within the time limit, we use the current lower bound of the subproblem to generate a valid integer cut of the form~\eqref{eq:integer-cut} and use the upper bound to update the objective value~$\theta^\star$ in step 4. If the main problem is not solved to optimality within the time limit, we 
report~$\theta^\star$ as the objective value and the current optimal objective value of \newqj{the} main problem~\eqref{eq:mins-w-master} as the lower bound.

\begin{landscape}
\setlength{\tabcolsep}{2pt}
\begin{table}[!ht]
\renewcommand{\arraystretch}{1.5}
\begin{footnotesize}

\begin{tabular}{cccccccccccccccccccc } 
 \hline
 & & & \multicolumn{5}{c}{MINLO} & & \multicolumn{5}{c}{ Decomposition} & & \multicolumn{5}{c}{ Decomposition with Added Cuts}\\
 $N$ & $K$ & $\,$ & Opt($\#$) & Time(s) & Gap & Better($\#$) & Improvement & $\,$ &  Opt($\#$) & Time(s) & Gap & Better($\#$) & Improvement & $\,$ &  Opt($\#$) & Time(s) & Gap & Better($\#$) & Improvement\\ 
 \hline
\multirow{3}{*}{5} & 2 && $ 0 $ & $ 7200.0 $ & $ - $ & $ 2 $ & $ -21.6\% $  && $ \textbf{6} $ & $ 1218.2 $ & $ 0.2\% $ & $ \textbf{5} $ & $ \textbf{30.4\%} $  && $ \textbf{6} $ & $ \textbf{426.0} $ & $ \textbf{0.1\%} $ & $ 4 $ & $ \textbf{30.4\%} $ \\
~ & 3 && $ 0 $ & $ 7200.0 $ & $ - $ & $ 1 $ & $ 14.7\% $  && $ 0 $ & $ 7200.0 $ & $ 5.6\% $ & $ \textbf{7} $ & $ 36.5\% $  && $ \textbf{2} $ & $ \textbf{663.0} $ & $ \textbf{1.1\%} $ & $ 5 $ & $ \textbf{36.9\%} $ \\
~ & 4 && $ 0 $ & $ 7200.0 $ & $ - $ & $ 2 $ & $ 13.4\% $  && $ 0 $ & $ 7200.0 $ & $ 9.9\% $ & $ 1 $ & $ 37.2\% $  && $ \textbf{2} $ & $ \textbf{936.0} $ & $ \textbf{6.1\%} $ & $ \textbf{13} $ & $ \textbf{37.6\%} $ \\
\hline
\multirow{3}{*}{10} & 2 && $ \textbf{1} $ & $ \textbf{4404.5} $ & $ - $ & $ 2 $ & $ 3.6\% $  && $ 0 $ & $ 7200.0 $ & $ 29.5\% $ & $ 3 $ & $ 43.3\% $  && $ 0 $ & $ 7200.0 $ & $ \textbf{20.2\%} $ & $ \textbf{11} $ & $ \textbf{43.8\%} $ \\
~ & 3 && $ \textbf{0} $ & $ \textbf{7200.0} $ & $ - $ & $ 3 $ & $ 38.8\% $  && $ \textbf{0} $ & $ \textbf{7200.0} $ & $ 26.4\% $ & $ 7 $ & $ 51.0\% $  && $ \textbf{0} $ & $ \textbf{7200.0} $ & $ \textbf{20.4\%} $ & $ \textbf{10} $ & $ \textbf{51.3\%} $ \\
~ & 4 && $ \textbf{0} $ & $ \textbf{7200.0} $ & $ - $ & $ 0 $ & $ 9.6\% $  && $ \textbf{0} $ & $ \textbf{7200.0} $ & $ 25.1\% $ & $ 8 $ & $ 53.9\% $  && $ \textbf{0} $ & $ \textbf{7200.0} $ & $ \textbf{20.0\%} $ & $ \textbf{12} $ & $ \textbf{54.2\%} $ \\
\hline
\multirow{3}{*}{15} & 2 && $ \textbf{0} $ & $ \textbf{7200.0} $ & $ - $ & $ 1 $ & $ -3.5\% $  && $ \textbf{0} $ & $ \textbf{7200.0} $ & $ 29.8\% $ & $ 7 $ & $ 36.9\% $  && $ \textbf{0} $ & $ \textbf{7200.0} $ & $ \textbf{28.1\%} $ & $ \textbf{11} $ & $ \textbf{37.1\%} $ \\
~ & 3 && $ \textbf{0} $ & $ \textbf{7200.0} $ & $ - $ & $ 1 $ & $ -20.9\% $  && $ \textbf{0} $ & $ \textbf{7200.0} $ & $ 25.9\% $ & $ 3 $ & $ 44.5\% $  && $ \textbf{0} $ & $ \textbf{7200.0} $ & $ \textbf{24.0\%} $ & $ \textbf{16} $ & $ \textbf{45.6\%} $ \\
~ & 4 && $ \textbf{0} $ & $ \textbf{7200.0} $ & $ - $ & $ 0 $ & $ -41.3\% $  && $ \textbf{0} $ & $ \textbf{7200.0} $ & $ 24.4\% $ & $ 6 $ & $ 47.4\% $  && $ \textbf{0} $ & $ \textbf{7200.0} $ & $ \textbf{22.3\%} $ & $ \textbf{14} $ & $ \textbf{48.7\%} $ \\
\hline
\multirow{3}{*}{20} & 2 && $ \textbf{0} $ & $ \textbf{7200.0} $ & $ - $ & $ 1 $ & $ -26.9\% $  && $ \textbf{0} $ & $ \textbf{7200.0} $ & $ 28.5\% $ & $ 6 $ & $ 40.1\% $  && $ \textbf{0} $ & $ \textbf{7200.0} $ & $ \textbf{27.6\%} $ & $ \textbf{13} $ & $ \textbf{40.6\%} $ \\
~ & 3 && $ \textbf{0} $ & $ \textbf{7200.0} $ & $ - $ & $ 1 $ & $ -19.2\% $  && $ \textbf{0} $ & $ \textbf{7200.0} $ & $ 25.6\% $ & $ 1 $ & $ 45.8\% $  && $ \textbf{0} $ & $ \textbf{7200.0} $ & $ \textbf{23.2\%} $ & $ \textbf{18} $ & $ \textbf{49.3\%} $ \\
~ & 4 && $ \textbf{0} $ & $ \textbf{7200.0} $ & $ - $ & $ 0 $ & $ -64.5\% $  && $ \textbf{0} $ & $ \textbf{7200.0} $ & $ 25.4\% $ & $ 1 $ & $ 46.3\% $  && $ \textbf{0} $ & $ \textbf{7200.0} $ & $ \textbf{21.6\%} $ & $ \textbf{19} $ & $ \textbf{52.4\%} $ \\
\hline
\end{tabular}

\caption{Summary of computational results of the R$\&$D problem with loans. $\mathbf{Decomposition}$ and $\mathbf{ Decomposition\;\;with\;\;Added\;\;Cuts}$ 
refer to Algorithm~\ref{alg:l-shaped} without and with the cuts introduced in Section~\ref{subsec:tighter-cuts}, respectively. $\mathbf{Opt(\#)}$ corresponds to the number of instances solved to optimality, $\mathbf{Time(s)}$ to the average computational time (in seconds) for instances solved to optimality, and $\mathbf{Gap}$ to the average optimality gap for the instances not solved within the time limit. We write '-' when no valid lower bound was found. $\mathbf{Better(\#)}$ denotes the number of instances in each method that achieved a better objective value than the other method. $\mathbf{Improvement}$ denotes the average improvement in the objective value of the $K$-adaptability solution found in the time limit over the static solution found in the time limit. The improvement is calculated as the ratio of the difference between the objective values of the $K$-adaptability solution and the static solution to the objective value of the static solution.}
\label{table:r&d-loan-per}
\end{footnotesize}
\end{table}
\end{landscape}

Table~\ref{table:r&d-loan-per} summarizes the computational results across those instances. From the table, we observe that the improved decomposition algorithm achieves optimality or finds better solutions for more instances than the basic decomposition algorithm and solves the monolithic MINLO problem when $K$ is large. Across all instances, the proposed algorithm has a significantly smaller optimality gap,\newqj{,} which decreases with~$K$. We notice that the MINLO problem is often unable to obtain a valid lower bound across all problem sizes. As the number of policies~$K$ increases, the number of instances for which the improved decomposition algorithm achieves a better objective value increases. For example, when~$N=20$, the number of better solutions found by our approach is~$\{13, \, 18, \, 19\}$ with~$K \in \{2, \,3, \,4\}$, respectively. The improvements in the objective value given by using the improved decomposition algorithm over the static solutions are $\{40.6\%, \, 49.3\%, \, 52.4\%\}$, whereas the improvements obtained by solving the MINLO problem and using the basic decomposition algorithm are $\{-26.9\%, \, -19.2\%, \,-64.5\%\}$ and $\{40.1\%, \, 45.8\%, \, 46.3\%\}$, respectively. In all cases, the improvement in the objective value produced by the decomposition algorithm over the static solutions increases with the number of policies~$K$. This is expected since larger values~$K$ offer more flexibility in the recourse decisions, leading to a better objective value. However, the trend is not clear in MINLO solutions due to computational difficulties, and the improvement can even be negative.

We also present the performance of the decomposition Algorithm \ref{alg:l-shaped} using the Branch-and-Bound Algorithm \ref{alg:bnb} with the cuts introduced in Section \ref{subsec:tighter-cuts}. Table~\ref{table:r&d-loan-bb} summarizes the computational results for the same randomly generated instances and should be read in companion with Table~\ref{table:r&d-per}. From the table, we observe that, across all cases, using Algorithm \ref{alg:bnb} to evaluate~$\Phi(\bm w)$ solves fewer instances to optimality and with better objective values. It also requires longer solving times, has larger gaps, and shows less improvement. Particularly, when the problem size is large (e.g.,\,$N > 10$), the decomposition algorithm using Algorithm \ref{alg:bnb} has difficulty finding high-quality solutions with positive improvement compared to the static solution. This again demonstrates the efficiency and advantages of the proposed Branch-and-Cut Algorithm \ref{alg:bnc}.

\setlength{\tabcolsep}{2pt}

\begin{table}[!htp]
\renewcommand{\arraystretch}{1.5}
\begin{footnotesize}
\centering
\begin{tabular}{cccccccccc } 
 \hline
 $N$ & $K$ & $\,$ & Opt($\#$) & Time(s) & Gap & Better($\#$) & Improvement \\ 
 \hline
\multirow{3}{*}{5} & 2 && $ 0 $ & $ 7200.0 $ & $ 27.1\% $ & $ 0 $ & $ 7.6\% $\\
~ & 3 && $ 0 $ & $ 7200.0 $ & $ 33.1\% $ & $ 0 $ & $ 4.6\% $\\
~ & 4 && $ 0 $ & $ 7200.0 $ & $ 34.0\% $ & $ 0 $ & $ 2.8\% $\\
\hline
\multirow{3}{*}{10} & 2 && $ 0 $ & $ 7200.0 $ & $ 48.6\% $ & $ 0 $ & $ -1.9\% $\\
~ & 3 && $ 0 $ & $ 7200.0 $ & $ 48.8\% $ & $ 0 $ & $ -2.3\% $\\
~ & 4 && $ 0 $ & $ 7200.0 $ & $ 48.3\% $ & $ 0 $ & $ -2.0\% $\\
\hline
\multirow{3}{*}{15} & 2 && $ 0 $ & $ 7200.0 $ & $ 53.6\% $ & $ 0 $ & $ -10.0\% $\\
~ & 3 && $ 0 $ & $ 7200.0 $ & $ 56.2\% $ & $ 0 $ & $ -13.2\% $\\
~ & 4 && $ 0 $ & $ 7200.0 $ & $ 56.1\% $ & $ 0 $ & $ -11.2\% $\\
\hline
\multirow{3}{*}{20} & 2 && $ 0 $ & $ 7200.0 $ & $ 61.8\% $ & $ 0 $ & $ -23.4\% $\\
~ & 3 && $ 0 $ & $ 7200.0 $ & $ 65.9\% $ & $ 0 $ & $ -27.9\% $\\
~ & 4 && $ 0 $ & $ 7200.0 $ & $ 72.3\% $ & $ 0 $ & $ -39.3\% $\\
\hline
\end{tabular}
\caption{Computational results of Algorithm~\ref{alg:l-shaped} using the Branch-and-Bound Algorithm \ref{alg:bnb} with the cuts introduced in Section~\ref{subsec:tighter-cuts} for R$\&$D problem with loans. $\mathbf{Opt(\#)}$ corresponds to the number of instances solved to optimality, $\mathbf{Time(s)}$ to the average computational time (in seconds) for instances solved to optimality, and $\mathbf{Gap}$ to the average optimality gap for the instances not solved within the time limit. $\mathbf{Better(\#)}$ denotes the number of instances in each method that achieved a better objective value than the other method. $\mathbf{Improvement}$ denotes the average improvement in the objective value of the $K$-adaptability solution found in the time limit over the static solution found in the time limit. The improvement is calculated as the ratio of the difference between the objective values of the $K$-adaptability solution and the static solution to the objective value of the static solution.}
\label{table:r&d-loan-bb}
\end{footnotesize}
\end{table}

Next, we examine the suboptimality of the RO solution. We define~$z^\star_\text{RO}$ and~$z^\star_\text{DRO}$ as in Section \ref{subsubsec: b&b-res}. The relative differences between~$z^\star_\text{RO}$ and~$z^\star_\text{DRO}$ are summarized in Table~\ref{table:r&d-loan-subopt}.
\begin{table}[!ht]
\renewcommand{\arraystretch}{1.5}
\begin{tabular}{ ccccc } 
 \hline
$K$ & $N=5$ & $N=10$ & $N=15$ & $N=20$ \\ 
 \hline
$ 1 $ & $15.74\%$ & $13.58\%$ & $20.38\%$ & $19.47\%$\\
$ 2 $ & $42.54\%$ & $49.12\%$ & $48.81\%$ & $50.5\%$\\
$ 3 $ & $43.83\%$ & $44.93\%$ & $49.71\%$ & $51.25\%$\\
$ 4 $ & $43.71\%$ & $43.83\%$ & $46.92\%$ & $45.71\%$\\
\hline
\end{tabular}
\caption{Suboptimality of RO solution to the R$\&$D problem with loans. Each cell represents an average of 20 instances. In each instance, the relative differences are calculated by $\displaystyle \frac{(z^*_\text{DRO} - z^*_\text{RO})}{(z^*_\text{DRO} + z^*_\text{RO})/2}$.}
\label{table:r&d-loan-subopt}
\end{table}
From the table, we see that the average suboptimality of the RO solution among all instances is up to 51.25\%, which shows the importance of incorporating distributional information in the uncertainty set.

We also compare the out-of-sample performance of the RO solution and the DRO solution. We sample 200 distributions that satisfy the moment constraints in~$\mathcal{P}_\text{RD}$. For each instance and distribution, we evaluate the expected net return~$\mathbb E_{\mathbb P}[-{\bm \xi^r}^\top (\bm w + \theta \bm y(\bm \xi)) + (1 + \alpha)x]$ using~$\bm y(\bm \xi)$ as the solution computed by RO and DRO approaches. The details of the sampling procedure can be found in Section \ref{ecsubsub:bb-out-sample}. Table \ref{table:r&d-loan-out-sample} reports the relative difference of the median, the~$95\%$ and~$99\%$ quantiles between the two methods. From the table, we can see that, on the simulated non-worst-case distributions, the DRO solutions outperform the RO solutions by up to~$49\%$.

\begin{table}[!htp]
\renewcommand{\arraystretch}{1.5}
% \begin{center}
\begin{tabular}{ cccccc } 
 \hline
$N$ & $K$ && median & $95\%$ quantile & $99\%$ quantile  \\ 
 \hline
\multirow{4}{*}{5} & $ 1 $ && $ 18.96\% $ & $ 18.66\% $ & $ 18.55\% $ \\
~ & $ 2 $ && $ 39.90\% $ & $ 40.01\% $ & $ 40.07\% $ \\
~ & $ 3 $ && $ 35.93\% $ & $ 35.97\% $ & $ 35.99\% $ \\
~ & $ 4 $ && $ 33.84\% $ & $ 33.92\% $ & $ 33.94\% $ \\
\hline
\multirow{4}{*}{10} & $ 1 $ && $ 17.68\% $ & $ 17.35\% $ & $ 17.27\% $ \\
~ & $ 2 $ && $ 47.68\% $ & $ 47.62\% $ & $ 47.62\% $ \\
~ & $ 3 $ && $ 31.81\% $ & $ 31.82\% $ & $ 31.88\% $ \\
~ & $ 4 $ && $ 38.31\% $ & $ 38.26\% $ & $ 38.28\% $ \\
\hline
\multirow{4}{*}{15} & $ 1 $ && $ 18.26\% $ & $ 18.05\% $ & $ 17.81\% $ \\
~ & $ 2 $ && $ 49.08\% $ & $ 48.98\% $ & $ 48.85\% $ \\
~ & $ 3 $ && $ 34.21\% $ & $ 34.33\% $ & $ 34.24\% $ \\
~ & $ 4 $ && $ 36.88\% $ & $ 36.87\% $ & $ 36.80\% $ \\
\hline
\multirow{4}{*}{20} & $ 1 $ && $ 18.13\% $ & $ 17.90\% $ & $ 17.75\% $ \\
~ & $ 2 $ && $ 49.43\% $ & $ 49.31\% $ & $ 49.25\% $ \\
~ & $ 3 $ && $ 30.77\% $ & $ 30.89\% $ & $ 30.90\% $ \\
~ & $ 4 $ && $ 36.33\% $ & $ 36.30\% $ & $ 36.27\% $ \\
 \hline
\end{tabular}

% \end{center}

\caption{The relative difference between the median, the ~$95\%$ and~$99\%$ quantiles of the expected objective values in the 200 simulated distributions given by the RO and DRO solutions to the R$\&$D problem with loans. Each cell represents an average of 20 instances.}
\label{table:r&d-loan-out-sample}

\end{table}

\end{appendix}

\clearpage

\end{document}